\documentclass[10pt]{amsart}

\usepackage[numbers,square]{natbib}
\usepackage{amsmath,amssymb,amsfonts,amsthm}
\usepackage{mathrsfs}
\usepackage{hyperref}

\hypersetup{hidelinks}
\allowdisplaybreaks[4]

\theoremstyle{plain}
\newtheorem{theorem}{Theorem}
\newtheorem{proposition}{Proposition}
\newtheorem{corollary}{Corollary}
\newtheorem{lemma}{Lemma}
\newtheorem{definition}{Definition}
\newtheorem{remark}{Remark}

\numberwithin{equation}{section}

\def\vint{\mathop{\mathchoice%
		{\setbox0\hbox{$\displaystyle\intop$}\kern 0.22\wd0%
			\vcenter{\hrule width 0.6\wd0}\kern -0.82\wd0}%
		{\setbox0\hbox{$\textstyle\intop$}\kern 0.2\wd0%
			\vcenter{\hrule width 0.6\wd0}\kern -0.8\wd0}%
		{\setbox0\hbox{$\scriptstyle\intop$}\kern 0.2\wd0%
			\vcenter{\hrule width 0.6\wd0}\kern -0.8\wd0}%
		{\setbox0\hbox{$\scriptscriptstyle\intop$}\kern 0.2\wd0%
			\vcenter{\hrule width 0.6\wd0}\kern -0.8\wd0}}%
	\mathopen{}\int}

\title[Singular sets of elliptic equations]
{Bounds of singular sets for elliptic equations in $C^{1,Dini}$ domains with singular potentials}

\author{Zhiwei Wang}
\address{Department of Mathematics, Louisiana State University, Baton Rouge, LA, USA}
\email{zwang30@lsu.edu}

\author{Jiuyi Zhu}
\address{Department of Mathematics, Louisiana State University, Baton Rouge, LA, USA}
\email{zhu@math.lsu.edu}


\begin{document}
	
	\begin{abstract}
		We study the quantitative codimension-two estimate for the singular set in the
		boundary neighborhoods of the solutions of 
		\[
		\Delta u+V(x)u=0\qquad\text{in }\Omega,
		\qquad
		u=0\qquad\text{on }\partial\Omega,
		\]
		where \(\Omega\subset\mathbb R^n\) is a bounded
		\(C^{1,\mathrm{Dini}}\) domain and \(V\in L^p(\Omega)\) for some
		\(p>n\).  
		We first prove
		the explicit upper bound for the  doubling index is given by 
		\(C(n,p,\Omega)(1+\|V\|_{L^p(\Omega)}^{\frac{2p}{3p-2n}})\).
		The analytic input is an interior volume estimate for solutions of
		second order elliptic equations 
		with uniformly elliptic Dini leading coefficients and \(V\in L^p\). Using
		the  quantitative doubling index bound,  boundary flattening, we show an explicit upper bound for singular sets in the neighborhood of the boundary of the \(C^{1,\mathrm{Dini}}\) domain.
	\end{abstract}
	
	\maketitle
	
	\section{Introduction}\label{sec:introduction}
	
	We study the singular sets of solutions to second-order elliptic equations
	with Dini continuous leading coefficients and an $L^p$ potential. Our main
	results are an interior volume estimate under a quantitative doubling
	assumption and an estimate for interior singular points approaching the
	boundary of a $C^{1,\mathrm{Dini}}$ domain.
	
	Let $n\ge2$ and $p>n$. Our first main interior result considers a nontrivial weak solution
	$u\in W^{1,2}(B_2)$ of
	\begin{equation}\label{eq:main-equation}
		\partial_i(a^{ij}(x)\partial_j u)+V(x)u=0\qquad\text{in }B_2.
	\end{equation}
	The matrix $a=(a^{ij})$ is symmetric and uniformly elliptic:
	\begin{equation}\label{eq:ellipticity}
		\lambda|\xi|^2\le a^{ij}(x)\xi_i\xi_j\le\Lambda|\xi|^2,
		\qquad x\in B_2,\quad\xi\in\mathbb R^n,
	\end{equation}
	where $0<\lambda\le\Lambda$. Its modulus is controlled by a
	nonnegative nondecreasing function $\omega$ such that
	\begin{equation}\label{eq:dini-modulus}
		\begin{gathered}
			|a^{ij}(x)-a^{ij}(y)|\le\omega(|x-y|)
			\quad(x,y\in B_2,\ |x-y|\le1),\\
			\omega(1)\le1,\qquad\int_0^1\frac{\omega(t)}t\,dt<\infty.
		\end{gathered}
	\end{equation}
	We assume
	\begin{equation}\label{eq:V-bound}
		\|V\|_{L^p(B_2)}\le M
	\end{equation}
	and the doubling estimate
	\begin{equation}\label{eq:euclidean-doubling-assumption}
		\vint_{B_{2r}(x)}u^2
		\le 4^{\mathcal D}\vint_{B_r(x)}u^2
		\qquad\text{whenever }B_{2r}(x)\subset B_2,
	\end{equation}
	for some \(\mathcal D\ge0\). This multiplicative form remains meaningful
	when the average on the smaller ball vanishes.
	Here $\vint$ denotes the integral average. Dini regularity and $p>n$
	give a continuous gradient. Thus the nodal and singular sets are
	defined pointwise by
	\begin{equation}\label{eq:singular-set-definition}
		Z(u)=\{x\in B_1:u(x)=0\},\qquad
		S(u)=\{x\in B_1:u(x)=0,\ \nabla u(x)=0\}.
	\end{equation}
	
	The codimension-two structure of singular and critical sets was developed
	in \cite{Han1994,HanHardtLin1998,HardtEtAl1999}. Quantitative
	stratification led to effective critical-set estimates for Lipschitz
	coefficients \cite{CheegerNaberValtorta2015} and to sharp
	codimension-two Minkowski bounds \cite{NaberValtorta2017}.
	Huang and Jiang \cite{HuangJiang} obtained singular-set volume estimates
	for H\"older continuous leading coefficients and bounded lower-order
	coefficients under a doubling assumption. We treat a general Dini modulus
	and an $L^p$ zeroth-order term. The doubling assumption remains essential
	for this coefficient class: in dimensions at least three, unique
	continuation can fail even for H\"older coefficients
	\cite{Miller1974,Mandache1996}.
	
	Define the accumulated modulus and its inverse scale by
	\begin{equation}\label{eq:5-Omega-omega}
		\Omega_\omega(t)=\omega(t)+\int_0^t\frac{\omega(s)}s\,ds,
		\qquad 0<t\le1,
	\end{equation}
	and
	\begin{equation}\label{eq:5-rho-omega}
		\rho_\omega(\delta)=
		\sup\{0<t\le1:\Omega_\omega(t)\le\delta\},
		\qquad0<\delta\le1.
	\end{equation}
	The Dini condition ensures that $\rho_\omega(\delta)>0$.
	
	For the constants $c$ and $C$ in the following theorem, put
	\begin{equation}\label{eq:5-final-delta}
		\delta
		=c(1+M)^{-C}
		\exp\bigl[-C(1+\mathcal D)^2\bigr].
	\end{equation}
	
	\begin{theorem}\label{thm:main-volume-estimate}
		Suppose that \eqref{eq:main-equation}--\eqref{eq:V-bound} and
		\eqref{eq:euclidean-doubling-assumption} hold.  There are constants
		\(c\in(0,1)\) and \(C>1\), depending only on
		\(n,p,\lambda,\Lambda\), and \(\omega\), with the following property.
		Then, for every $0<r\le1$,
		\begin{equation}\label{eq:5-final-clean-D-M}
			\begin{split}
				&\left|\{z\in B_1:\operatorname{dist}(z,S(u))<r\}\right|\\
				&\qquad\le C(1+M)^{C}
				\exp\bigl[C(1+\mathcal D)^2\bigr]
				\rho_\omega(\delta)^{-2}r^2.
			\end{split}
		\end{equation}
		Consequently, the upper Minkowski dimension of $S(u)$ is at most $n-2$,
		and its upper $(n-2)$-dimensional Minkowski content is finite.
	\end{theorem}
	
	In particular, under a
	H\"older bound on $\omega$, this gives the following explicit estimate compared to the result in \cite{HuangJiang}.
	
	\begin{corollary}\label{cor:main-holder}
		Under the assumptions of Theorem~\ref{thm:main-volume-estimate}, suppose
		that
		\begin{equation}\label{eq:5-holder-modulus-corollary}
			\omega(t)\le t^\alpha\quad(0<t\le1),
			\qquad 0<\alpha\le1.
		\end{equation}
		Then there is $C=C(n,p,\lambda,\Lambda,\alpha)$ such that
		\begin{equation}\label{eq:5-final-holder-clean}
			\left|\{z\in B_1:\operatorname{dist}(z,S(u))<r\}\right|
			\le C(1+M)^C
			\exp\bigl[C(1+\mathcal D)^2\bigr]r^2,
			\qquad 0<r\le1.
		\end{equation}
	\end{corollary}
	
	We next consider the boundary problem. Boundary unique continuation on
	Dini domains was established in
	\cite{AdolfssonEscauriaza1997,KukavicaNystrom1998}.
	Kenig and Zhao \cite{KenigZhao2022Boundary} proved quantitative
	codimension-two Minkowski bounds and rectifiability for the singular set
	of a harmonic function vanishing on an open portion of a Dini boundary.
	Their set includes boundary points. They subsequently estimated the full
	critical set near a $C^{1,\alpha}$ boundary by flattening and reflection,
	followed by the interior theorem of Huang and Jiang
	\cite{KenigZhao2025}. This reduction points to the role of Dini leading
	coefficients in the corresponding problem on $C^{1,\mathrm{Dini}}$
	domains.
	
	Our boundary result allows an $L^p$ potential. It concerns interior
	singular points uniformly near the boundary. The singular-set condition
	$u=0$ is necessary for the argument: when $V\ne0$, subtracting the value
	of $u$ at a critical point produces an inhomogeneous equation. The global
	Dirichlet condition also matters. A local estimate retains a fixed-scale
	doubling quantity, which we bound here by propagation of mass through
	the bounded connected domain.
	
	Let $\Omega\subset\mathbb R^n$ be a bounded connected
	$C^{1,\mathrm{Dini}}$ domain. Let $p>n$, $V\in L^p(\Omega)$, and let
	$u\in W^{1,2}_0(\Omega)$ be a nontrivial weak solution of
	\begin{equation}\label{eq:intro-boundary-schrodinger}
		\Delta u+V(x)u=0\qquad\text{in }\Omega.
	\end{equation}
	Put again $M=\|V\|_{L^p(\Omega)}$ and
	\begin{equation}\label{eq:intro-boundary-singular-set}
		S_\Omega(u)=\{x\in\Omega:u(x)=0,\ \nabla u(x)=0\}.
	\end{equation}
	Let $\omega_{\partial\Omega}$ denote a modulus of continuity of the
	boundary normal in the boundary charts. For the geometric constants
	$r_0,c,C$ in the following theorem, write
	\begin{equation}\label{eq:intro-boundary-transformed-modulus}
		\widetilde\omega(t)
		=C\omega_{\partial\Omega}(C r_0 t),
		\qquad 0<t\le1,
	\end{equation}
	and
	\begin{equation}\label{eq:intro-boundary-delta-global}
		\begin{split}
			\widetilde{\delta}
			={}&c
			(1+r_0^{2-n/p}M)^{-C}
			\exp\left[-C\left(1+M^{\frac{2p}{3p-2n}}\right)^2\right].
		\end{split}
	\end{equation}
	\begin{theorem}\label{thm:boundary-volume-estimate}
		There are $r_0,c\in(0,1)$ and $C>1$, depending only on $n,p$ and
		the quantitative $C^{1,\mathrm{Dini}}$ character of $\Omega$, such that
		$\widetilde\omega(1)\le1$ and, for every $x_0\in\partial\Omega$ and
		$0<r\le cr_0$,
		\begin{equation}\label{eq:intro-boundary-volume-estimate}
			\begin{split}
				&\left|\left\{x\in\Omega\cap B_{c r_0}(x_0):
				\operatorname{dist}\bigl(x,S_\Omega(u)\cap
				B_{2c r_0}(x_0)\bigr)<r\right\}\right|\\
				&\quad\le
				\frac{Cr_0^{n-2}r^2}
				{\widetilde{\delta}\,\rho_{\widetilde\omega}(\widetilde{\delta})^2}.
			\end{split}
		\end{equation}
	\end{theorem}
	
	Consequently, the upper Minkowski dimension of
	$S_\Omega(u)\cap B_{c r_0}(x_0)$ is at most $n-2$, and its upper
	$(n-2)$-dimensional Minkowski content is finite.
	
	The interior proof uses harmonic approximation and quantitative
	uniqueness of the homogeneous polynomial when the doubling index is
	nearly constant. Cone splitting and a covering argument then give the
	volume estimate as in \cite{HuangJiang,NaberValtorta2017}. The biggest difference lies Dini summability controls the errors accumulated over
	the relevant scales.
	
	For the boundary application, we check the doubling bound estimate using a weighted frequency
	argument, which gives doubling bound \(C(n,p,\Omega)(1+\|V\|_{L^p(\Omega)}^{\frac{2p}{3p-2n}})\).
	Then we flatten the boundary using regularized normal coordinates and
	reflect the equation as in \cite{KenigZhao2022Boundary,KenigZhao2025} to reach problem \eqref{eq:main-equation}, to which the interior theorem applies.
	
	We use the following convention for constants. In the interior argument, a
	structural constant depends only on \(n,p,\lambda,\Lambda\) and the fixed
	Dini modulus \(\omega\) in \eqref{eq:dini-modulus}; in the boundary argument,
	it depends only on \(n,p\) and the quantitative \(C^{1,\mathrm{Dini}}\)
	character of \(\Omega\). In particular, structural constants are independent
	of \(M\) and \(\mathcal D\). Constants denoted by \(c>0\) and \(C>0\), without
	subscripts, are structural and may decrease or increase from line to line. A subscript is retained
	only when a constant is fixed for use in a later definition or iteration. Its
	source is stated when it is introduced and recalled later whenever it is not
	immediate from the context.

    \begin{remark}
       Almost at the same time, a similar upper bound as Theorem \ref{thm:main-volume-estimate} for the critical sets of elliptic equations with Dini continuous leading
coefficients was obtained independently by Jiang and Wang in \cite{JW26}  (See Theorem 1.5) with some different arguments.
    \end{remark}

    	{\bf Acknowledgment:} We used ChatGPT 5.6 Sol to check our calculations and assist with the layout. All underlying ideas and insights are from the authors.
	
	\section{Almost monotone formula}\label{sec:almost-monotonicity}
	
	From this section through Section \ref{sec:covering}, we assume
	\eqref{eq:main-equation}--\eqref{eq:euclidean-doubling-assumption}. We will repeatedly use
	following notations. 
	
	For \(x\in B_1\), let
	\(A_x=a(x)^{1/2}\) be the positive definite square root of \(a(x)\).
	For scalar- or vector-valued functions \(f,g\) with square-integrable
	traces on \(\partial B_1\), set
	\begin{equation}\label{eq:standing-spherical-normalization}
		\langle f,g\rangle=\vint_{\partial B_1}f\cdot g,
		\qquad
		\|f\|^2=\langle f,f\rangle.
	\end{equation}
	Throughout Section \ref{sec:almost-monotonicity} to \ref{sec:singular-covering-volume}, a function is called normalized 
	when \(\|f\|=1\).
	Any normalization
	at a different sphere or ball will be written explicitly.
	
	For a function \(f\), the spherical doubling index is defined by
	\begin{equation}\label{eq:euclidean-spherical-doubling}
		D^f(z,s)=\log_4
		\frac{\vint_{\partial B_{2s}(z)}f^2}
		{\vint_{\partial B_s(z)}f^2}.
	\end{equation}
	This convention applies, in particular, to the harmonic comparison
	functions \(h\) and the rescaled functions \(u_{x,r}\) below. For the
	original solution in the \(x\)-coordinates, the notation \(D^u(x,r)\) is
	reserved for the one defined in
	\eqref{eq:singular-doubling}.
	
	Let \(x\in B_1\) and \(r>0\) satisfy
	\(x+2rA_xB_1\subset B_2\). Whenever the relevant averages below are
	positive, define the normalized rescaling
	\begin{equation}\label{eq:uxr-definition}
		u_{x,r}(y)=\frac{u(x+rA_xy)}
		{\left(\vint_{\partial B_1}u(x+rA_xz)^2\,dz\right)^{1/2}}.
	\end{equation}
	Then \(\|u_{x,r}\|=1\). We define the ball doubling index by
	\begin{equation}\label{eq:ball-average-doubling-index}
		\widetilde D_x(r)
		=\log_4\frac{\vint_{B_2}u(x+rA_xy)^2\,dy}
		{\vint_{B_1}u(x+rA_xy)^2\,dy}
		=\log_4\frac{\vint_{x+2rA_xB_1}u(z)^2\,dz}
		{\vint_{x+rA_xB_1}u(z)^2\,dz}.
	\end{equation}
	The spherical doubling index is
	\begin{equation}\label{eq:singular-doubling}
		D^u(x,r)=\log_4\frac{\vint_{\partial B_2}u(x+rA_xy)^2\,dy}
		{\vint_{\partial B_1}u(x+rA_xy)^2\,dy}.
	\end{equation}
	By \eqref{eq:euclidean-spherical-doubling} and
	\eqref{eq:uxr-definition}, \(D^u(x,r)=D^{u_{x,r}}(0,1)\).
	We remark that, compared to the study of critical set, the rescaling in \eqref{eq:uxr-definition} does not subtract \(u(x)\),
	because it will be used at centers in the zero set of \(u\). Lemma
	\ref{lem:positive-spherical-average} proves that the relevant spherical
	averages are positive under the small-scale condition
	\eqref{eq:spherical-normalization-range}; hence all later uses of
	\eqref{eq:uxr-definition} and \eqref{eq:singular-doubling} are well defined.

	\subsection{Interior estimates and rescaling}
	
	\begin{proposition}\label{prop:C1-interior}
		There are
		\(C,q>0\), depending only on the structural data, such that
		\begin{equation}\label{eq:C1-estimate-normalized}
			\|u\|_{L^\infty(B_{1/2})}+\|\nabla u\|_{L^\infty(B_{1/4})}
			\le C(1+M)^q\|u\|_{L^2(B_1)}.
		\end{equation}
		For the remainder of the paper, fix
		\begin{equation}\label{eq:K-definition}
			K:=C(1+M)^q\ge1.
		\end{equation}
	\end{proposition}
	
	\begin{proof}
		Since \(p>n\), one has \(2-n/p>1\). Fix
		\(\zeta\in C_c^1(B_1)\) such that
		\[
		0\le\zeta\le1,
		\qquad \zeta=1\ \text{on }B_{2/3},
		\qquad |\nabla\zeta|\le4.
		\]
		Fix a structural constant \(C\ge1\) large enough for the Caccioppoli and
		Sobolev estimates associated with this cutoff. Choose a structural
		constant \(0<c\le1/8\) so that
		\(C^2c^{2-n/p}\le1/2\), and set
		\[
		\rho:=c(1+M)^{-1/(2-n/p)}.
		\]
		In particular,
		\[
		0<\rho\le c\le\frac18,
		\qquad C^2c^{2-n/p}\le\frac12.
		\]
		For \(x_0\in B_{1/2}\), one has
		\(B_\rho(x_0)\subset B_1\). Set
		\(v(y)=u(x_0+\rho y)\). On \(B_1\),
		\[
		\partial_i(a^{ij}(x_0+\rho y)\partial_jv)
		+V_{x_0,\rho}(y)v=0,
		\qquad V_{x_0,\rho}(y)=\rho^2V(x_0+\rho y),
		\]
		and the change of variables gives
		\begin{equation}\label{eq:C1-rescaled-potential}
			\|V_{x_0,\rho}\|_{L^p(B_1)}
			\le \rho^{2-n/p}M
			=c^{2-n/p}\frac{M}{1+M}
			\le c^{2-n/p}
			\le\frac1{2C^2}.
		\end{equation}
		Testing the equation by \(\zeta^2v\) and applying H\"older's inequality
		gives
		\[
		\begin{aligned}
			\int_{B_1}\zeta^2|\nabla v|^2
			&\le C\int_{B_1}v^2
			+C\|V_{x_0,\rho}\|_{L^p(B_1)}
			\|\zeta v\|_{L^{2p/(p-1)}(B_1)}^2.
		\end{aligned}
		\]
		Since $p>n$, H\"older inequality and Sobolev inequality yields
		\[
		\|\zeta v\|_{L^{2p/(p-1)}(B_1)}^2
		\le C\left(\int_{B_1}\zeta^2|\nabla v|^2
		+\int_{B_1}v^2\right).
		\]
		Consequently, combining above estimate and \eqref{eq:C1-rescaled-potential} gives
		\begin{align*}
			\int_{B_1}\zeta^2|\nabla v|^2
			&\le C\left(1+\|V_{x_0,\rho}\|_{L^p(B_1)}\right)
			\int_{B_1}v^2\\
			&\quad+C^2\|V_{x_0,\rho}\|_{L^p(B_1)}
			\int_{B_1}\zeta^2|\nabla v|^2\\
			&\le\left(C+\frac1{2C}\right)\int_{B_1}v^2
			+\frac12\int_{B_1}\zeta^2|\nabla v|^2.
		\end{align*}
		After absorption,
		\[
		\int_{B_{2/3}}|\nabla v|^2
		\le3C\int_{B_1}v^2.
		\]
		Since \(p>n\) and the norm in
		\eqref{eq:C1-rescaled-potential} is bounded independently of \(M\),
		the local boundedness estimate gives
		\begin{equation}\label{eq:C1-local-boundedness-caccioppoli}
			\|v\|_{L^\infty(B_{3/4})}
			+\|\nabla v\|_{L^2(B_{2/3})}
			\le C\|v\|_{L^2(B_1)}.
		\end{equation}
		Moreover,
		\[
		\|V_{x_0,\rho}v\|_{L^p(B_{2/3})}
		\le Cc^{2-n/p}\|v\|_{L^2(B_1)}.
		\]
		The rescaled coefficient has modulus \(t\mapsto\omega(\rho t)\), and
		\[
		\int_0^1\frac{\omega(\rho t)}t\,dt
		=\int_0^\rho\frac{\omega(s)}s\,ds
		\le\int_0^1\frac{\omega(s)}s\,ds.
		\]
		Uniform Dini continuity
		implies the partial Dini mean-oscillation hypothesis of
		\cite[Theorem~1.2]{ChoiKimLee2020}. Since we assume
		\(p>n\), \cite[Remarks~1.3--1.4]{ChoiKimLee2020} permits the scalar
		right-hand side in that estimate. Applying it to
		\[
		\operatorname{div}\bigl(a(x_0+\rho\,\cdot)\nabla v\bigr)
		=-V_{x_0,\rho}v,
		\]
		together with \eqref{eq:C1-local-boundedness-caccioppoli}, it yields
		\[
		\|v\|_{L^\infty(B_{1/2})}
		+\|\nabla v\|_{L^\infty(B_{1/4})}
		\le C\|v\|_{L^2(B_1)}.
		\]
		Scaling back gives
		\[
		\|u\|_{L^\infty(B_{\rho/2}(x_0))}
		+\rho\|\nabla u\|_{L^\infty(B_{\rho/4}(x_0))}
		\le C\rho^{-n/2}\|u\|_{L^2(B_\rho(x_0))}
		\le C\rho^{-n/2}\|u\|_{L^2(B_1)}.
		\]
		Taking the supremum over \(x_0\in B_{1/2}\) for the first term and
		over \(x_0\in B_{1/4}\) for the second gives
		\[
		\rho^{-n/2}\le\rho^{-n/2-1}
		=c^{-n/2-1}(1+M)^{(n/2+1)/(2-n/p)}.
		\]
		Since \(c\) is independent of \(M\), its fixed power is absorbed into
		the structural constant. Therefore
		\[
		\|u\|_{L^\infty(B_{1/2})}
		+\|\nabla u\|_{L^\infty(B_{1/4})}
		\le C(1+M)^{(n/2+1)/(2-n/p)}\|u\|_{L^2(B_1)}.
		\]
		Thus \eqref{eq:C1-estimate-normalized} holds with
		\[
		q=\frac{n/2+1}{2-n/p}
		=\frac{p(n+2)}{2(2p-n)}.
		\]
	\end{proof}
	
	Whenever \(x\) and \(r\) are fixed, we use the
	following coefficients on the rescaled domain
	\(\{y\in\mathbb R^n:x+rA_xy\in B_2\}\):
	\begin{equation}\label{eq:rescaled-coefficients-definition}
		\widetilde a(y)=A_x^{-1}a(x+rA_xy)A_x^{-1},
		\qquad
		\widetilde V(y)=r^2V(x+rA_xy).
	\end{equation}
	
	The next lemma verifies that the spherical averages used in the
	normalization and the doubling index do not vanish at the small scales
	considered below.
	
	\begin{lemma}\label{lem:positive-spherical-average}
		There is \(c_1=c_1(n,p,\lambda,\Lambda)>0\), independent of \(M\) and
		\(\mathcal D\), such that the following holds. If \(x\in B_1\) and
		\begin{equation}\label{eq:spherical-normalization-range}
			0<r\le\frac1{16\sqrt\Lambda},
			\qquad M(2r)^{2-n/p}\le c_1,
		\end{equation}
		then, for every \(0<s\le2r\),
		\begin{equation}\label{eq:positive-spherical-average}
			\vint_{\partial B_1}u(x+sA_xy)^2\,dy>0.
		\end{equation}
	\end{lemma}
	
	\begin{proof}
		Fix \(0<s\le2r\). In particular,
		\(s<1/\sqrt\Lambda\), so \(x+sA_xB_1\subset B_2\). Use
		\eqref{eq:rescaled-coefficients-definition} at the center \(x\) and
		radius \(s\).
		Suppose that the average in \eqref{eq:positive-spherical-average}
		vanishes. The trace of \(y\mapsto u(x+sA_xy)\) is zero on
		\(\partial B_1\), so this function belongs to \(W^{1,2}_0(B_1)\). By
		\eqref{eq:rescaled-coefficients-definition}, it satisfies
		\[
		\partial_{y_i}\Bigl(\widetilde a^{ij}(y)
		\partial_{y_j}\bigl[u(x+sA_xy)\bigr]\Bigr)
		+\widetilde V(y)u(x+sA_xy)=0
		\qquad\text{in }B_1.
		\]
		The ellipticity of \(a\), the change of variables, and
		Sobolev--Poincar\'e give
		\begin{align*}
			\frac\lambda\Lambda\int_{B_1}
			|\nabla_y u(x+sA_xy)|^2\,dy
			&\le \int_{B_1}|\widetilde V(y)|u(x+sA_xy)^2\,dy\\
			&\le \|\widetilde V\|_{L^p(B_1)}
			\|u(x+sA_x\,\cdot)\|_{L^{2p/(p-1)}(B_1)}^2\\
			&\le CMs^{2-n/p}\int_{B_1}
			|\nabla_y u(x+sA_xy)|^2\,dy\\
			&\le CM(2r)^{2-n/p}\int_{B_1}
			|\nabla_y u(x+sA_xy)|^2\,dy,
		\end{align*}
		where \(C=C(n,p,\lambda)\). Choose \(c_1\) so that
		\(Cc_1\le\lambda/(2\Lambda)\). The preceding estimate then implies
		\(\nabla_y u(x+sA_xy)=0\) in \(B_1\), and the zero trace gives
		\(u(x+sA_xy)=0\) there. Since
		\[
		B_{\sqrt\lambda s}(x)\subset x+sA_xB_1,
		\]
		we have \(u=0\) on \(B_{\sqrt\lambda s}(x)\), contrary to its nontriviality by strong unique continuation.
	\end{proof}

	\begin{lemma}\label{lem:scaled-coefficients}
		Let \(x\in B_1\), and assume \eqref{eq:spherical-normalization-range}.
		If \(R>0\) satisfies
		\[
		x+rA_xB_R\subset B_2,
		\qquad 2\sqrt\Lambda rR\le1,
		\]
		then the transformed equation is
		\begin{equation}\label{eq:scaled-equation}
			\partial_i(\widetilde a^{ij}\partial_j u_{x,r})
			+\widetilde V u_{x,r}=0
			\qquad\text{in }B_R.
		\end{equation}
		Moreover, \(\widetilde a(0)=I\), and on \(B_R\),
		\begin{align}
			|\widetilde a(y)-\widetilde a(z)|
			&\le C(\lambda,\Lambda)\omega(\sqrt\Lambda r|y-z|),
			\label{eq:scaled-a-modulus}\\
			\|\widetilde V\|_{L^p(B_R)}
			&\le C(n,p,\lambda)Mr^{2-n/p}.
			\label{eq:scaled-V-Lp}
		\end{align}
		The constants in the transformed ellipticity bounds depend only on
		\(\lambda,\Lambda\). Finally, for \(0<s\le1\),
		\begin{equation}\label{eq:singular-scaling-identity}
			D^u(x,rs)=D^{u_{x,r}}(0,s).
		\end{equation}
	\end{lemma}
	
	\begin{proof}
		For
		\(X=x+rA_xy\), the identities
		\[
		\nabla_Xu(X)=\frac{\left(\vint_{\partial B_1}u(x+rA_xz)^2\,dz\right)^{1/2}}rA_x^{-1}\nabla_y u_{x,r}(y),
		\qquad dX=r^n\det(A_x)\,dy
		\]
		inserted into the weak formulation give \eqref{eq:scaled-equation}.
		Moreover,
		\begin{align*}
			|\widetilde a(y)-\widetilde a(z)|
			&\le \|A_x^{-1}\|^2
			\omega(r|A_x(y-z)|)
			\le\lambda^{-1}\omega(\sqrt\Lambda r|y-z|),\\
			\|\widetilde V\|_{L^p(B_R)}^p
			&=r^{2p-n}(\det A_x)^{-1}
			\int_{x+rA_xB_R}|V(X)|^p\,dX
			\le \lambda^{-n/2}M^pr^{2p-n}.
		\end{align*}
		Since \(A_x^2=a(x)\), one has \(\widetilde a(0)=I\); and
		\[
		\frac\lambda\Lambda|\xi|^2
		\le \widetilde a^{ij}(y)\xi_i\xi_j
		\le\frac\Lambda\lambda|\xi|^2.
		\]
		Finally the multiplicative factor \(\left(\vint_{\partial B_1}u(x+rA_xz)^2\,dz\right)^{1/2}\) cancels in both spherical
		averages, and replacing \(y\) by \(sy\) gives
		\[
		D^{u_{x,r}}(0,s)
		=\log_4\frac{\vint_{\partial B_2}u(x+rsA_xy)^2\,dy}
		{\vint_{\partial B_1}u(x+rsA_xy)^2\,dy}
		=D^u(x,rs).
		\]
	\end{proof}
	
	\begin{lemma}\label{lem:doubling-over-balls-bound}
		If \(x\in B_1\) and
		\(B_{2\sqrt\Lambda r}(x)\subset B_2\), then there is
		\(C_0=C_0(n,\lambda,\Lambda)\ge1\) such that
		\begin{equation}\label{eq:ball-doubling-bound}
			\widetilde D_x(r)\le C_0(1+\mathcal D).
		\end{equation}
	\end{lemma}
	
	\begin{proof}
		The spectral bounds for \(A_x\) give
		\[
		B_{\sqrt\lambda s}(x)\subset x+sA_xB_1
		\subset B_{\sqrt\Lambda s}(x),
		\qquad |x+sA_xB_1|=s^n\det(A_x)|B_1|.
		\]
		Consequently,
		\begin{align*}
			\vint_{x+2rA_xB_1}u^2
			&\le\left(\frac\Lambda\lambda\right)^{n/2}
			\vint_{B_{2\sqrt\Lambda r}(x)}u^2,\\
			\vint_{x+rA_xB_1}u^2
			&\ge\left(\frac\lambda\Lambda\right)^{n/2}
			\vint_{B_{\sqrt\lambda r}(x)}u^2.
		\end{align*}
		Let
		\[
		N=\left\lceil\log_2\left(2\sqrt{\Lambda/\lambda}\right)\right\rceil.
		\]
		Put \(s_0=2^{1-N}\sqrt\Lambda r\). Then
		\(\sqrt\lambda r/2<s_0\le\sqrt\lambda r\) and
		\(2^Ns_0=2\sqrt\Lambda r\). Iterating
		\eqref{eq:euclidean-doubling-assumption} along
		\(s_0,2s_0,\ldots,2^Ns_0\), whose largest ball is the ball in the
		hypothesis, and then comparing the volumes of \(B_{s_0}\) and
		\(B_{\sqrt\lambda r}\), gives
		\[
		\frac{\vint_{B_{2\sqrt\Lambda r}(x)}u^2}
		{\vint_{B_{\sqrt\lambda r}(x)}u^2}
		\le4^{N\mathcal D}
		\frac{\vint_{B_{s_0}(x)}u^2}
		{\vint_{B_{\sqrt\lambda r}(x)}u^2}
		\le2^n4^{N\mathcal D}.
		\]
		Here every intermediate ball is contained in
		\(B_{2\sqrt\Lambda r}(x)\subset B_2\). Therefore
		\[
		\widetilde D_x(r)
		\le n\log_4(\Lambda/\lambda)+\frac n2+N\mathcal D
		\le C_0(n,\lambda,\Lambda)(1+\mathcal D).
		\]
	\end{proof}
	
	\subsection{Harmonic approximation}
	
	Let \(K\) be the constant fixed in \eqref{eq:K-definition}, and let
	\(C_0\) be the constant from Lemma
	\ref{lem:doubling-over-balls-bound}, enlarged when necessary below.
	For \(0<\varepsilon\le1/10\), define
	\begin{equation}\label{eq:ha-r-epsilon}
		r_1(\varepsilon)
		=\sup\left\{0<r\le\frac1{16\sqrt\Lambda}:
		\omega(8\sqrt\Lambda r)+Mr^{2-n/p}
		\le c\varepsilon^{(n+2)/2}
		\bigl(K4^{3C_0(1+\mathcal D)}\bigr)^{-(n+2)/2}\right\}.
	\end{equation}
	The Dini condition and \(2-n/p>0\) imply
	\(r_1(\varepsilon)>0\).
	We fix \(c\) in \eqref{eq:ha-r-epsilon} so that
	\(2^{2-n/p}c\le c_1\), where $c_1$ comes from Lemma \ref{lem:positive-spherical-average}. Since \(\varepsilon\le1/10\) and
	\(K4^{3C_0(1+\mathcal D)}\ge1\), every
	\(0<r\le r_1(\varepsilon)\) satisfies
	\[
	0<r\le\frac1{16\sqrt\Lambda},
	\qquad M(2r)^{2-n/p}\le c_1.
	\]
	Hence Lemma \ref{lem:positive-spherical-average} makes every normalization
	and spherical doubling index used below well defined.
	If a defining supremum is not attained, all assertions at that endpoint
	are understood by passage to the limit. This convention will also be used
	for the radii defined later.
	
	\begin{lemma}\label{lem:harmonic-approximation-singular}
		If \(x\in Z(u)\cap B_1\) and \(0<r\le r_1(\varepsilon)\), there
		is a harmonic function \(h\) in \(B_4\) such that
		\begin{equation}\label{eq:ha-conclusion}
			h(0)=0,\qquad \vint_{\partial B_1}h^2=1,
			\qquad \|h-u_{x,r}\|_{L^\infty(B_3)}\le\varepsilon.
		\end{equation}
	\end{lemma}
	
	\begin{proof}
		Normalize first by the ball average, writing
		\(w(y)=u(x+rA_xy)/(\vint_{B_1}u(x+rA_xz)^2\,dz)^{1/2}\).
		The choice of radius makes the transformed equation available on \(B_8\).
		The normalization gives \(\vint_{B_1}w^2=1\). At the scales
		\(r,2r,4r\), Lemma \ref{lem:doubling-over-balls-bound} gives
		\[
		\widetilde D_x(2^jr)
		=\log_4\frac{\vint_{B_{2^{j+1}}}w^2}
		{\vint_{B_{2^j}}w^2}
		\le C_0(1+\mathcal D),
		\qquad j=0,1,2.
		\]
		Thus
		\[
		\vint_{B_8}w^2\le4^{3C_0(1+\mathcal D)}.
		\]
		The scaled form of \eqref{eq:C1-estimate-normalized}, applied on fixed
		balls covering \(B_5\) and contained in \(B_8\), then yields
		\[
		\|w\|_{C^1(B_5)}
		\le CK\left(\vint_{B_8}w^2\right)^{1/2}
		\le CK4^{3C_0(1+\mathcal D)}.
		\]
		Let \(g\) be the harmonic replacement of \(w\) in \(B_4\), i.e,
		$$\begin{cases}
			\Delta g=0,\ \text{in}\ B_4\\
			g=w,\ \text{on}\ \partial B_4.
		\end{cases}$$
		Then
		\(v=w-g\) has zero boundary values and satisfies
		\[
		\Delta v=\operatorname{div}((I-\widetilde a)\nabla w)-\widetilde Vw.
		\]
		Let \(G_4(y,z)\) be the Green function of \(B_4\). For
		\(y\in B_4\), integration by parts in the weak equation gives
		\[
		v(y)=-\int_{B_4}\nabla_zG_4(y,z)\cdot
		(I-\widetilde a(z))\nabla w(z)\,dz
		-\int_{B_4}G_4(y,z)\widetilde V(z)w(z)\,dz.
		\]
		The Green function bounds
		\[
		\sup_{y\in B_4}\int_{B_4}|\nabla_zG_4(y,z)|\,dz\le C,
		\qquad
		\sup_{y\in B_4}\|G_4(y,\cdot)\|_{L^{p'}(B_4)}\le C
		\]
		hold because \(p>n\). Together with
		\eqref{eq:scaled-a-modulus}, \eqref{eq:scaled-V-Lp}, and the preceding
		\(C^1\) bound, they give
		\begin{equation}\label{eq:ha-sup-small}
			\|v\|_{L^\infty(B_4)}
			\le CK4^{3C_0(1+\mathcal D)}
			\bigl(\omega(8\sqrt\Lambda r)+Mr^{2-n/p}\bigr).
		\end{equation}
		Then, the triangle inequality and monotonicity of harmonic spherical means give
		\[
		\left(\vint_{B_1}g^2\right)^{1/2}\ge1-\|v\|_\infty\ge\frac12,
		\qquad
		\left(\vint_{\partial B_1}g^2\right)^{1/2}
		\ge\left(\vint_{B_1}g^2\right)^{1/2}.
		\]
		Also \(w(0)=0\), so \(|g(0)|\le\|v\|_\infty\). Hence
		\[
		\left(\vint_{\partial B_1}w^2\right)^{1/2}\ge\frac14,
		\qquad
		\left(\vint_{\partial B_1}|g-g(0)|^2\right)^{1/2}\ge\frac14
		\]
		after decreasing the structural constant in \eqref{eq:ha-r-epsilon}.
		
		Set
		\[
		h=\frac{g-g(0)}{(\vint_{\partial B_1}|g-g(0)|^2)^{1/2}}.
		\]
		Since \(u_{x,r}=w/(\vint_{\partial B_1}w^2)^{1/2}\), we obtain
		\[
		\|h-u_{x,r}\|_{L^\infty(B_3)}
		\le C\|v\|_\infty\bigl(1+\|w\|_{L^\infty(B_3)}\bigr).
		\]
		In view of the \(C^1\) bound for \(w\), the right-hand side is at most
		\[
		C\bigl(K4^{3C_0(1+\mathcal D)}\bigr)^2
		\bigl(\omega(8\sqrt\Lambda r)+Mr^{2-n/p}\bigr)
		\le Cc\varepsilon^{(n+2)/2}
		\bigl(K4^{3C_0(1+\mathcal D)}\bigr)^{(2-n)/2}
		\le\varepsilon.
		\]
		This proves the assertion after fixing \(c\) sufficiently small.
	\end{proof}
	
	\begin{lemma}\label{lem:spherical-doubling-bound-singular}
		Choose \(0<c\le1/10\) sufficiently small and set
		\(r_2=r_1(c/(1+\mathcal D))\), where \(r_1\) is defined in
		\eqref{eq:ha-r-epsilon}.
		For \(x\in Z(u)\cap B_1\) and \(0<r\le r_2\),
		\begin{equation}\label{eq:lem254-conclusion}
			D^u(x,r)\le C_0(1+\mathcal D).
		\end{equation}
	\end{lemma}
	
	\begin{proof}
		Use the functions \(w,g\), and \(v=w-g\) from the preceding proof with
		\(\varepsilon=c/(1+\mathcal D)\). Choosing \(c\) sufficiently small gives
		\[
		\|v\|_{L^\infty(B_4)}\le\frac18,
		\qquad
		\left(\vint_{\partial B_1}w^2\right)^{1/2}
		\ge\frac14.
		\]
		The two bounds
		\[
		\widetilde D_x(r),\ \widetilde D_x(2r)
		\le C_0(1+\mathcal D)
		\]
		from \eqref{eq:ball-doubling-bound} give
		\[
		\left(\vint_{B_4}w^2\right)^{1/2}
		\le4^{C_0(1+\mathcal D)}.
		\]
		Hence
		\[
		\left(\vint_{B_4}g^2\right)^{1/2}
		\le4^{C_0(1+\mathcal D)}+\|v\|_{L^\infty(B_4)}.
		\]
		Since \(g\) is harmonic on \(B_4\), its interior estimate on \(B_3\)
		gives
		\[
		\left(\vint_{\partial B_2}g^2\right)^{1/2}
		\le\|g\|_{L^\infty(B_3)}
		\le C\left(\vint_{B_4}g^2\right)^{1/2}
		\le C4^{C_0(1+\mathcal D)}.
		\]
		The triangle inequality on \(\partial B_2\) now yields
		\[
		\frac{\vint_{\partial B_2}w^2}{\vint_{\partial B_1}w^2}
		\le16\left(C4^{C_0(1+\mathcal D)}+\|v\|_{L^\infty(B_4)} \right)^2
		\le C4^{2C_0(1+\mathcal D)}.
		\]
		With \(C_0\) fixed at the outset large enough, taking \(\log_4\) proves
		\eqref{eq:lem254-conclusion}.
	\end{proof}
	
	\begin{lemma}\label{lem:comparison-harmonic}
		For \(0<\varepsilon\le1/10\), using \(r_1\) from
		\eqref{eq:ha-r-epsilon}, put
		\begin{equation}\label{eq:harmonic-comparison-radius}
			r_3(\varepsilon)
			=r_1(\varepsilon^{10C_0(1+\mathcal D)}).
		\end{equation}
		At every \(x\in Z(u)\cap B_1\) and \(0<r\le r_3(\varepsilon)\),
		the harmonic function in \eqref{eq:ha-conclusion} can be chosen so that
		\begin{equation}\label{eq:harmonic-doubling-comparison}
			|D^h(0,s)-D^{u_{x,r}}(0,s)|\le\varepsilon,
			\qquad \varepsilon^2\le s\le1.
		\end{equation}
	\end{lemma}
	
	\begin{proof}
		
		Fix \(C_0\) so that, for \(0<\varepsilon\le1/10\),
		\[
		\varepsilon^{10C_0(1+\mathcal D)}
		\le\frac{c}{1+\mathcal D}.
		\]
		Then \(r_3(\varepsilon)\le r_2\), and
		\[
		D^{u_{x,r}}(0,t)=D^u(x,rt)\le C_0(1+\mathcal D),
		\qquad0<t\le1.
		\]
		For \(1\le t\le2\), harmonic approximation and increasing harmonic
		spherical means give
		\[
		\left(\vint_{\partial B_{t}}u_{x,r}^2\right)^{1/2}\ge \left(\vint_{\partial B_{t}}h^2\right)^{1/2}
		-\varepsilon^{10C_0(1+\mathcal D)}\ge\frac12.
		\]
		For \(\varepsilon^2\le s<1\), choose the integer \(j\ge1\) with
		\(1\le2^js<2\). Then \(j\le1+2\log_2(1/\varepsilon)\), and the preceding doubling bound give
		\begin{equation}\label{eq:comparison-height-lower}
			\begin{split}
				\left(\vint_{\partial B_{s}}u_{x,r}^2\right)
				&=\left(\vint_{\partial B_{2^js}}u_{x,r}^2\right)\prod_{k=0}^{j-1}4^{-D^{u_{x,r}}(0,2^ks)}\\
				&\ge4^{-jC_0(1+\mathcal D)}\left(\vint_{\partial B_{2^js}}u_{x,r}^2\right)
				\ge\frac14\,4^{-C_0(1+\mathcal D)}
				\varepsilon^{4C_0(1+\mathcal D)}.
			\end{split}
		\end{equation}
		For \(1\le s\le2\), the earlier estimate
		\(\left(\vint_{\partial B_{s}}u_{x,r}^2\right)^{1/2}\ge1/2\) gives the same lower bound after decreasing its
		constant.
		The triangle inequality now gives
		\[
		|\left(\vint_{\partial B_{s}}h^2\right)^{1/2}-\left(\vint_{\partial B_{s}}u_{x,r}^2\right)^{1/2}|
		\le\varepsilon^{10C_0(1+\mathcal D)}
		\le\frac{\varepsilon}{100}\left(\vint_{\partial B_{s}}u_{x,r}^2\right)^{1/2},
		\]
		where \(C_0\) is fixed sufficiently large. Thus, at both radii,
		the quotient
		\[
		\frac{\left(\vint_{\partial B_s}h^2\right)^{1/2}}
		{\left(\vint_{\partial B_s}u_{x,r}^2\right)^{1/2}}
		\]
		lies in \([1-\varepsilon/100,1+\varepsilon/100]\). Applying this at
		\(s\) and \(2s\), we obtain
		\[
		D^h(0,s)-D^{u_{x,r}}(0,s)
		=2\log_4
		\frac{\left(\vint_{\partial B_{2s}}h^2\right)^{1/2}/\left(\vint_{\partial B_{2s}}u_{x,r}^2\right)^{1/2}}
		{\left(\vint_{\partial B_{s}}h^2\right)^{1/2}/\left(\vint_{\partial B_{s}}u_{x,r}^2\right)^{1/2}},
		\]
		and hence
		\[
		|D^h(0,s)-D^{u_{x,r}}(0,s)|
		\le2\log_4\frac{1+\varepsilon/100}{1-\varepsilon/100}
		\le\varepsilon.
		\]
		This proves \eqref{eq:harmonic-doubling-comparison}.
	\end{proof}
	
	\subsection{Almost monotonicity}
	
	For a harmonic function, let
	\begin{equation}\label{eq:singular-frequency}
		N^h(a,r)=\frac{r\int_{B_r(a)}|\nabla h|^2}
		{\int_{\partial B_r(a)}h^2}.
	\end{equation}
	Orthogonal expansion into homogeneous harmonic polynomials shows that
	\(N^h(a,r)\) is nondecreasing and
	\begin{equation}\label{eq:frequency-doubling-comparison}
		N^h(a,r)\le D^h(a,r)\le N^h(a,2r).
	\end{equation}
	Consequently,
	\begin{equation}\label{eq:harmonic-doubling-separated-scales}
		D^h(a,s)\le D^h(a,r)\qquad\text{whenever }0<2s\le r,
	\end{equation}
	because
	\[
	D^h(a,s)\le N^h(a,2s)\le N^h(a,r)\le D^h(a,r).
	\]
	We use the following standard consequences of this expansion; see
	\cite[Lemmas~2.17 and 2.20]{HuangJiang} and
	\cite[Lemma~3.16 and Theorem~3.19]{NaberValtorta2017}.
	
	\begin{lemma}\label{lem:pinching-integer}
		Let $h$ be harmonic, then there is a  \(\varepsilon_0(n)>0\) with the following properties.
		If \(D^h(0,1)-D^h(0,1/20)\le\varepsilon\le\varepsilon_0\), then there is an
		integer \(d\ge0\) satisfies
		\begin{equation}\label{eq:pinching-conclusion}
			|D^h(0,s)-d|\le3\varepsilon,\qquad 1/10<s<1/2.
		\end{equation}
		If \(D^h(0,1)\le d-\varepsilon\), then
		\begin{equation}\label{eq:drop-conclusion}
			D^h(0,\varepsilon/2)\le d-1+\varepsilon.
		\end{equation}
		In particular, \(D^h(0,1)\le d+1/2\) implies
		\(D^h(0,\varepsilon/2)\le d+\varepsilon\), after decreasing
		\(\varepsilon_0\).
	\end{lemma}
	
	Let \(\varepsilon_0(n)\) be the constant supplied by Lemma
	\ref{lem:pinching-integer}, and set
	\[
	a_0=\min\{1/40,10\varepsilon_0(n)\}.
	\]
	We prove the main theorem in this subsection, namely, the almost monotonicity of doubling index.
	\begin{theorem}\label{thm:almost-monotonicity}
		Let \(0<\varepsilon\le1/10\), and set
		\begin{equation}\label{eq:almost-monotonicity-radius}
			r_4(\varepsilon)=r_3(a_0\varepsilon).
		\end{equation}
		Here \(r_3\) is defined in \eqref{eq:harmonic-comparison-radius}, and
		\(a_0\) is the constant chosen above from Lemma
		\ref{lem:pinching-integer}.
		Then, for \(x\in Z(u)\cap B_1\) and \(0<2s\le r\le r_4(\varepsilon)\),
		\begin{equation}\label{eq:almost-monotonicity}
			D^u(x,s)\le D^u(x,r)+\varepsilon.
		\end{equation}
		At the same radii, \(D^u(x,r)\ge1-\varepsilon\).
	\end{theorem}
	
	\begin{proof}
		Fix \(x\in Z(u)\cap B_1\).
		At every scale \(t\le r_4(\varepsilon)\), Lemma \ref{lem:comparison-harmonic} is available
		with error \(a_0\varepsilon\) on relative radii
		\([(a_0\varepsilon)^2,1]\).
		For the comparison function at scale \(t\), write
		\(h=\sum_{k\ge1}P_k\) in homogeneous harmonic polynomials. Since
		\(h(0)=0\), orthogonality gives
		\[
		4^{D^h(0,1)}
		=\frac{\sum_{k\ge1}4^k\|P_k\|^2}
		{\sum_{k\ge1}\|P_k\|^2}\ge4.
		\]
		Consequently,
		\begin{equation}\label{eq:almost-monotonicity-lower-bound}
			D^u(x,t)\ge D^h(0,1)-a_0\varepsilon
			\ge1-a_0\varepsilon\ge1-\varepsilon.
		\end{equation}
		We first prove that, for every integer \(d\ge0\),
		\begin{equation}\label{eq:integer-upper-bound-propagation}
			D^u(x,t)\le d+\varepsilon/10
			\quad\Longrightarrow\quad
			D^u(x,s)\le d+\varepsilon/2\quad(0<s\le t/2).
		\end{equation}
		Indeed, the harmonic comparison at scale \(t\) satisfies
		\[
		D^h(0,1)\le D^u(x,t)+a_0\varepsilon
		\le d+\varepsilon/10+a_0\varepsilon<d+1/2.
		\]
		Since \(a_0\varepsilon\le\varepsilon_0(n)\), the last assertion of
		Lemma \ref{lem:pinching-integer}, used with parameter
		\(a_0\varepsilon\), gives
		\[
		D^h(0,a_0\varepsilon/2)\le d+a_0\varepsilon,
		\qquad
		D^u(x,a_0\varepsilon t/2)\le d+2a_0\varepsilon.
		\]
		Set \(t_k=(a_0\varepsilon/2)^kt\). We prove inductively that
		\(D^u(x,t_k)\le d+\varepsilon/10\). The case \(k=0\) is the hypothesis.
		If it holds at \(k\), the comparison function \(h_k\) at scale \(t_k\)
		satisfies
		\begin{align*}
			D^{h_k}(0,1)&\le D^u(x,t_k)+a_0\varepsilon<d+1/2,\\
			D^{h_k}(0,a_0\varepsilon/2)&\le d+a_0\varepsilon,\\
			D^u(x,t_{k+1})
			&\le D^{h_k}(0,a_0\varepsilon/2)+a_0\varepsilon
			\le d+2a_0\varepsilon\le d+\varepsilon/10.
		\end{align*}
		Moreover, for every
		\((a_0\varepsilon)^2\le\sigma\le1/2\), comparison and
		\eqref{eq:harmonic-doubling-separated-scales} give
		\[
		D^u(x,\sigma t_k)\le D^{h_k}(0,\sigma)+a_0\varepsilon
		\le D^{h_k}(0,1)+a_0\varepsilon
		\le d+\varepsilon/10+2a_0\varepsilon\le d+\varepsilon/2.
		\]
		Thus the desired bound holds on
		\([ (a_0\varepsilon)^2t_k,t_k/2]\). These intervals cover
		\((0,t/2]\): indeed,
		\[
		\frac{t_{k+1}}2=\frac{a_0\varepsilon}{4}t_k
		\ge (a_0\varepsilon)^2t_k,
		\]
		because \(a_0\varepsilon\le1/4\), and their lower endpoints tend to
		zero. This proves \eqref{eq:integer-upper-bound-propagation}.
		
		For every \(r\le r_4(\varepsilon)\), choose
		\(d=\lfloor D^u(x,r)+\varepsilon/2\rfloor\).
		If \(D^u(x,r)\le d+\varepsilon/10\), we apply
		\eqref{eq:integer-upper-bound-propagation}: since
		\(d\le D^u(x,r)+\varepsilon/2\), this gives
		\(D^u(x,s)\le D^u(x,r)+\varepsilon\).
		Otherwise,
		\[
		d+\varepsilon/10<D^u(x,r)<d+1-\varepsilon/2.
		\]
		Let \(h\) be the harmonic comparison function at scale \(r\). Then
		\[
		D^h(0,1)\le D^u(x,r)+a_0\varepsilon
		\le d+1-a_0\varepsilon,
		\]
		because \(a_0\varepsilon\le\varepsilon/4\). Applying
		\eqref{eq:drop-conclusion} with the integer \(d+1\) and parameter
		\(a_0\varepsilon\), then comparing at relative radius
		\(a_0\varepsilon/2\), gives
		\[
		D^h(0,a_0\varepsilon/2)\le d+a_0\varepsilon,
		\qquad D^u(x,a_0\varepsilon r/2)
		\le d+2a_0\varepsilon\le d+\varepsilon/10.
		\]
		Thus \eqref{eq:integer-upper-bound-propagation}, applied at
		\(a_0\varepsilon r/2\), gives
		\(D^u(x,s)\le d+\varepsilon/2\le D^u(x,r)+\varepsilon\) for
		\(0<s\le a_0\varepsilon r/4\). If
		\(a_0\varepsilon r/4\le s\le r/2\), then
		\[
		(a_0\varepsilon)^2\le a_0\varepsilon/4\le s/r\le1/2,
		\]
		and comparison together with
		\eqref{eq:harmonic-doubling-separated-scales} gives
		\[
		D^u(x,s)\le D^h(0,s/r)+a_0\varepsilon
		\le D^h(0,1)+a_0\varepsilon
		\le D^u(x,r)+2a_0\varepsilon\le D^u(x,r)+\varepsilon.
		\]
		This proves the upper bound for every \(0<2s\le r\); the lower bound
		is \eqref{eq:almost-monotonicity-lower-bound}.
	\end{proof}
	
	\section{Quantitative uniqueness of tangent maps}\label{sec:tangent-uniqueness}
	
	We continue to use the spherical norm
	\eqref{eq:standing-spherical-normalization}. We recall three facts that
	will also be used in the next section. Firstly: if \(P,Q\) are
	homogeneous harmonic polynomials of degree
	\(d\), then
	\begin{equation}\label{eq:hhp-ball-sphere}
		\vint_{B_1}P^2=\frac n{n+2d}\vint_{\partial B_1}P^2.
	\end{equation}
	The other two read like below.
	\begin{lemma}\label{lem:grad-inner-product-hhp}
		For homogeneous harmonic polynomials of the same degree \(d\ge0\),
		\begin{equation}\label{eq:grad-inner-product-hhp}
			\langle\nabla P,\nabla Q\rangle=d(2d+n-2)\langle P,Q\rangle,
			\qquad
			\|P\|_{L^\infty(\partial B_1)}
			\le C(n)(1+d)^{(n-1)/2}\|P\|.
		\end{equation}
		If \(d\ge1\) and \(P\) is orthogonal to the degree-\(d\) polynomials
		invariant along \(\operatorname{span}\{e\}\), where \(e\) is a unit
		vector, then \(\|P\|\le\|\partial_eP\|\).
	\end{lemma}
	These are the standard spherical harmonic identities and estimates; see
	\cite[Lemmas~3.2--3.3 and Proposition~3.9]{NaberValtorta2017}, or
	\cite[Lemmas~3.3, 3.5 and~3.10]{HuangJiang}.
	
	With \(C_0\) denoting the constant from Lemma
	\ref{lem:doubling-over-balls-bound}, enlarged as in Lemmas
	\ref{lem:spherical-doubling-bound-singular} and
	\ref{lem:comparison-harmonic}, and with \(K\) from
	\eqref{eq:K-definition}, define
	\begin{equation}\label{eq:42-AD-definition}
		A_{\mathcal D}=C(n,\lambda,\Lambda)4^{C_0(1+\mathcal D)},
		\qquad d_{\mathcal D}=\lceil C_0(1+\mathcal D)\rceil,
	\end{equation}
	where the constant in the exponent is chosen larger than the constant
	needed in the spherical doubling bound. This allows the exponential
	factors from that bound to be absorbed into \(A_{\mathcal D}\).
	We use the error
	\begin{equation}\label{eq:prop335-Dini-smallness}
		\eta(r):=K A_{\mathcal D}
		\left(\int_0^{2\sqrt\Lambda r}\frac{\omega(s)}s\,ds
		+Mr^{2-n/p}\right).
	\end{equation}
	In particular, \(\eta(r)\to0\) as \(r\downarrow0\).
	We continue to use the radius \(r_4\) defined in
	\eqref{eq:almost-monotonicity-radius}.
	
	The next lemma gives quantitative $C^1$ estimate for the tangent maps $u_{x,r}$ by assuming lower bound of doubling index, which plays key role in the proof of quantitative uniqueness of tangent maps.
	\begin{lemma}\label{lem:growth-estimate-gradient}
		Let \(x\in Z(u)\cap B_1\), \(10t_2\le t_1\le r_4(1/10)\),
		and \(D^u(x,s)\ge\gamma\ge0\) on \([t_2,t_1]\). Then
		\begin{align}
			\vint_{B_t}u_{x,t_1}^2
			&\le A_{\mathcal D}
			\begin{cases}\left(\frac{t_2}{t_1}\right)^{\,2\gamma-9/5}t^{9/5},&t\le \frac{t_2}{t_1},\\
				t^{2\gamma},&\frac{t_2}{t_1}\le t\le2,
			\end{cases}\label{eq:growth-ball}\qquad\\
			\sup_{B_t}|u_{x,t_1}|
			&\le CK A_{\mathcal D}^{1/2}
			\begin{cases}\left(\frac{t_2}{t_1}\right)^{\,\gamma-1}t,&t\le \frac{t_2}{t_1},\\
				t^\gamma,&\frac{t_2}{t_1}\le t\le1,
			\end{cases}\label{eq:growth-sup}\\
			\sup_{B_t}|\nabla u_{x,t_1}|
			&\le CK A_{\mathcal D}^{1/2}
			\begin{cases}\left(\frac{t_2}{t_1}\right)^{\,\gamma-1},&t\le \frac{t_2}{t_1},\\
				t^{\gamma-1},&\frac{t_2}{t_1}\le t\le1.
			\end{cases}\label{eq:growth-gradient}
		\end{align}
	\end{lemma}
	\begin{proof}
		Since
		\(\left(\vint_{\partial B_{1}}u_{x,t_1}^2\right)=1\), Lemma \ref{lem:spherical-doubling-bound-singular} and Lemma \ref{lem:comparison-harmonic} are available at the scale \(t_1\) give, after fixing the
		leading constant in \(A_{\mathcal D}\),
		\begin{equation}\label{eq:growth-sphere-large}
			\left(\vint_{\partial B_{t}}u_{x,t_1}^2\right)\le A_{\mathcal D},\qquad1\le t\le2.
		\end{equation}
		Indeed, since
		\(t_1\le r_4(1/10)=r_3(a_0/10)\), choose the harmonic comparison at
		scale \(t_1\) with parameter \(a_0/10\). Then \(h\) is harmonic in
		\(B_4\), and \eqref{eq:ha-conclusion} gives
		\[
		\|h-u_{x,t_1}\|_{L^\infty(B_3)}
		\le (a_0/10)^{10C_0(1+\mathcal D)}\le1.
		\]
		Since \(h^2\) is subharmonic, its spherical averages are
		nondecreasing. Therefore, for \(1\le t\le2\),
		\begin{align*}
			\left(\vint_{\partial B_t}u_{x,t_1}^2\right)^{1/2}
			&\le\left(\vint_{\partial B_t}h^2\right)^{1/2}+1\\
			&\le\left(\vint_{\partial B_2}h^2\right)^{1/2}+1\\
			&\le\left(\vint_{\partial B_2}u_{x,t_1}^2\right)^{1/2}+2\\
			&=2^{D^{u_{x,t_1}}(0,1)}+2
			=2^{D^u(x,t_1)}+2\\
			&\le C2^{C_0(1+\mathcal D)}.
		\end{align*}
		Here we used \(\vint_{\partial B_1}u_{x,t_1}^2=1\), the scaling
		identity, and Lemma \ref{lem:spherical-doubling-bound-singular}.
		Squaring and enlarging the leading constant in \(A_{\mathcal D}\) proves
		\eqref{eq:growth-sphere-large}.
		
		Let \(\frac{t_2}{t_1}\le s\le1\), and choose \(i\ge0\) so that
		\(1\le2^is<2\). The scaling identity and the hypothesis give
		\(D^{u_{x,t_1}}(0,2^ks)=D^u(x,2^kst_1)\ge\gamma\) for
		\(0\le k<i\). Hence
		\begin{align}
			\left(\vint_{\partial B_{s}}u_{x,t_1}^2\right)
			&=\left(\vint_{\partial B_{2^is}}u_{x,t_1}^2\right)\prod_{k=0}^{i-1}
			4^{-D^{u_{x,t_1}}(0,2^ks)}
			\le A_{\mathcal D}4^{-i\gamma}
			\le A_{\mathcal D}s^{2\gamma}.
			\label{eq:growth-sphere-middle}
		\end{align}
		Theorem \ref{thm:almost-monotonicity}, with error \(1/10\), gives
		\[
		D^{u_{x,t_1}}(0,s)=D^u(x,st_1)\ge\frac9{10},
		\qquad0<s\le1.
		\]
		If \(0<s\le \frac{t_2}{t_1}\), choose \(i\ge0\) with
		\(\frac{t_2}{t_1}\le2^is<\frac{2t_2}{t_1}\). Using the last lower bound and then
		\eqref{eq:growth-sphere-middle},
		\begin{align}
			\left(\vint_{\partial B_{s}}u_{x,t_1}^2\right)
			&\le4^{-9i/10}\left(\vint_{\partial B_{2^is}}u_{x,t_1}^2\right)
			\le C A_{\mathcal D}4^{-9i/10}\left(\frac{t_2}{t_1}\right)^{2\gamma}\notag\\
			&\le C A_{\mathcal D}
			\left(\frac{t_2}{t_1}\right)^{2\gamma-9/5}s^{9/5}.
			\label{eq:growth-sphere-small}
		\end{align}
		
		The polar-coordinate identity
		\begin{equation}\label{eq:growth-radial-identity}
			\vint_{B_t}u_{x,t_1}^2
			=\frac n{t^n}\int_0^t s^{n-1}\left(\vint_{\partial B_{s}}u_{x,t_1}^2\right)\,ds
		\end{equation}
		and \eqref{eq:growth-sphere-small} give, for \(t\le \frac{t_2}{t_1}\),
		\[
		\vint_{B_t}u_{x,t_1}^2
		\le\frac{C A_{\mathcal D}\left(\frac{t_2}{t_1}\right)^{2\gamma-9/5}}{t^n}
		\int_0^t s^{n-1+9/5}\,ds
		\le A_{\mathcal D}\left(\frac{t_2}{t_1}\right)^{2\gamma-9/5}t^{9/5}.
		\]
		If \(\frac{t_2}{t_1}\le t\le2\), split
		\eqref{eq:growth-radial-identity} at \(\frac{t_2}{t_1}\) and use both spherical
		bounds:
		\begin{align*}
			\vint_{B_t}u_{x,t_1}^2
			&\le\frac{C A_{\mathcal D}}{t^n}
			\left(\left(\frac{t_2}{t_1}\right)^{2\gamma-9/5}
			\int_0^{\frac{t_2}{t_1}}s^{n-1+9/5}\,ds
			+\int_{\frac{t_2}{t_1}}^t s^{n-1+2\gamma}\,ds\right)\\
			&\le\frac{C A_{\mathcal D}}{t^n}
			\left(\left(\frac{t_2}{t_1}\right)^{n+2\gamma}+t^{n+2\gamma}\right)
			\le A_{\mathcal D}t^{2\gamma}.
		\end{align*}
		This proves \eqref{eq:growth-ball}.
		
		Let \(\frac{t_2}{t_1}\le t\le1\) and \(y\in B_t\). The scaled form of
		Proposition \ref{prop:C1-interior} on \(B_t(y)\subset B_{2t}\) gives
		\begin{align*}
			|\nabla u_{x,t_1}(y)|
			&\le CKt^{-1}
			\left(\vint_{B_t(y)}u_{x,t_1}^2\right)^{1/2}\\
			&\le CKt^{-1}
			\left(2^n\vint_{B_{2t}}u_{x,t_1}^2\right)^{1/2}
			\le CK A_{\mathcal D}^{1/2}t^{\gamma-1}.
		\end{align*}
		For \(0<t\le \frac{t_2}{t_1}\) and \(y\in B_t\), apply the same estimate on
		\(B_{\frac{t_2}{t_1}}(y)\subset B_{\frac{2t_2}{t_1}}\):
		\[
		|\nabla u_{x,t_1}(y)|
		\le CK\left(\frac{t_2}{t_1}\right)^{-1}
		\left(2^n\vint_{B_{\frac{2t_2}{t_1}}}u_{x,t_1}^2\right)^{1/2}
		\le CK A_{\mathcal D}^{1/2}\left(\frac{t_2}{t_1}\right)^{\gamma-1}.
		\]
		Taking suprema proves \eqref{eq:growth-gradient}. Finally,
		\(u_{x,t_1}(0)=0\), and the line segment from zero to any
		\(y\in B_t\) lies in \(B_t\). Therefore
		\[
		|u_{x,t_1}(y)|
		\le |y|\sup_{B_t}|\nabla u_{x,t_1}|
		\le t\sup_{B_t}|\nabla u_{x,t_1}|,
		\]
		which gives \eqref{eq:growth-sup} in both ranges.
	\end{proof}
	
	For a function \(f\) with nonzero spherical trace, write
	\begin{equation}\label{eq:tangent-spherical-normalization}
		f_{0,t}(z)=\frac{f(tz)}
		{\left(\vint_{\partial B_t}f^2\right)^{1/2}}.
	\end{equation}
	
	The proposition below gives an approximation harmonic function, through which we can relate $u_{x,r}$ to a homogenous harmonic polynomial. We fix a structural constant \(0<\varepsilon_1<1/10\) and recall $r_4(\varepsilon)$ comes from Theorem \ref{thm:almost-monotonicity}, $\eta (r)$ is defined in \eqref{eq:prop335-Dini-smallness}.
	
	\begin{proposition}\label{prop:quantitative-harmonic-approximation-DiniLp}
		Let \(x\in Z(u)\cap B_1\), and suppose
		\[
		10t_2\le t_1\le\min\{r_4(1/10),(10\sqrt\Lambda)^{-1}\},
		\qquad
		\mathop{\mathrm{osc}}_{[t_2,t_1]}D^u(x,\cdot)\le\varepsilon_1,
		\]
		and let \(\gamma=\inf_{[t_2,t_1]}D^u(x,\cdot)\). There is a harmonic function
		\(h\) in \(B_1\), with \(h(0)=0\), such that
		\begin{align}
			|h(y)-u_{x,t_1}(y)|&\le C\eta(t_1)|y|^\gamma
			&&(\frac{t_2}{t_1}\le |y|\le3/4),\label{eq:prop335-DiniLp-annulus}\\
			|h(y)-u_{x,t_1}(y)|&\le C\eta(t_1)\left(\frac{t_2}{t_1}\right)^{\,\gamma-1}|y|
			&&(|y|\le \frac{t_2}{t_1}).\label{eq:prop335-DiniLp-inner}
		\end{align}
	\end{proposition}
	
	\begin{proof}
		Let \(d=\lfloor\gamma\rfloor\). In the coordinates of
		Lemma \ref{lem:scaled-coefficients},
		\[
		\Delta u_{x,t_1}=\operatorname{div}((I-\widetilde a)\nabla u_{x,t_1})
		-\widetilde V u_{x,t_1}.
		\]
		We use the expansion of the fundamental solution at the origin,
		\(\Gamma(y,z)=\sum_{k\ge0}\Gamma_k(z)P_{z,k}(y)\), with
		\(P_{z,k}\) homogeneous harmonic of degree \(k\) and normalized on
		\(\partial B_1\). The elementary kernel bounds used here are
		\begin{align}
			|\Gamma_k(z)|&\le C(4/3)^k|z|^{2-n-k},\nonumber\\
			|\nabla_z(\Gamma_k(z)P_{z,k}(y))|
			&\le C(4/3)^k(1+k)^{(n-1)/2}|z|^{1-n-k}|y|^k,
			\label{eq:Green-coefficient-bound}\\
			\left|\nabla_z^j\left(\Gamma(y,z)-\sum_{k=0}^d
			\Gamma_k(z)P_{z,k}(y)\right)\right|
			&\le C2^{d+1}|y|^{d+1}|z|^{1-n-d-j},\quad j=0,1,
			\label{eq:Green-remainder-bound}
		\end{align}
		when \(|z|\ge2|y|\). For the proof, see for example  \cite[Lemma 3.23]{HuangJiang}; in dimension
		two the coefficient bound is used only for \(k\ge1\), and the
		constant term is treated as \(\Gamma(y,z)-\Gamma(0,z)\).
		
		Define
		\begin{equation}\label{eq:prop335-phi-definition}
			\begin{split}
				\phi(y)={}&\int_{B_1}(\Gamma(y,z)-\Gamma(0,z))\Delta u_{x,t_1}(z)\,dz\\
				&-\sum_{k=1}^d\int_{B_1\setminus B_{\frac{t_2}{t_1}}}
				\Gamma_k(z)P_{z,k}(y)\Delta u_{x,t_1}(z)\,dz,
				\qquad h=u_{x,t_1}-\phi.
			\end{split}
		\end{equation}
		Put
		\[
		F=(I-\widetilde a)\nabla u_{x,t_1},
		\qquad \Delta u_{x,t_1}=\operatorname{div}F-\widetilde V u_{x,t_1}.
		\]
		We remark the meaning of above equation: no second derivative of \(u_{x,t_1}\) is required in
		\eqref{eq:prop335-phi-definition}. On a truncated annulus, an integral
		against \(\operatorname{div}F\) is defined by integration by parts,
		including its boundary terms. Near the kernel poles we first remove
		small balls and then pass to the limit. Since \(F\) is continuous, the
		inner boundary contributions are \(O(\delta)\) in dimension at least
		three and \(O(\delta|\log\delta|)\) in dimension two. Both tend to zero.
		The resulting volume integrals converge by the estimates below.
		Equivalently one can convolve the distributional equation on a slightly
		larger ball and pass to the limit. Since
		\[
		\Delta_y(\Gamma(y,z)-\Gamma(0,z))=\delta_z,
		\qquad \Delta_yP_{z,k}(y)=0,
		\]
		the definition \eqref{eq:prop335-phi-definition} gives
		\(\Delta\phi=\Delta u_{x,t_1}\) in \(B_1\). Hence \(\Delta h=0\)
		there, while \(\phi(0)=u_{x,t_1}(0)=0\) gives \(h(0)=0\).
		
		We first estimate \(\phi(y)\) for \(\frac{t_2}{t_1}\le |y|\le1/4\). Write
		\begin{equation}\label{eq:prop335-decomposition}
			\phi(y)=I_1(y)+I_2(y)+I_3(y),
		\end{equation}
		where
		\begin{align}
			I_1(y)&=\int_{B_{2|y|}}(\Gamma(y,z)-\Gamma(0,z))
			\Delta u_{x,t_1}(z)\,dz,\label{eq:prop335-I1}\\
			I_2(y)&=-\sum_{k=1}^d\int_{\frac{t_2}{t_1}<|z|<2|y|}
			\Gamma_k(z)P_{z,k}(y)\Delta u_{x,t_1}(z)\,dz,
			\label{eq:prop335-I2}\\
			I_3(y)&=\int_{2|y|<|z|<1}
			\left(\Gamma(y,z)-\sum_{k=0}^d\Gamma_k(z)P_{z,k}(y)\right)
			\Delta u_{x,t_1}(z)\,dz.
			\label{eq:prop335-I3}
		\end{align}
		By Lemma \ref{lem:growth-estimate-gradient},
		\begin{equation}\label{eq:prop335-source-bounds}
			\begin{split}
				|F(z)|&\le CK A_{\mathcal D}^{1/2}
				\omega(\sqrt\Lambda t_1|z|)
				\begin{cases}\left(\frac{t_2}{t_1}\right)^{\gamma-1},&|z|\le \frac{t_2}{t_1},\\
					|z|^{\gamma-1},&\frac{t_2}{t_1}\le|z|\le1,
				\end{cases}\\
				\|\widetilde V\|_{L^p(B_1)}&\le CMt_1^{2-n/p}.
			\end{split}
		\end{equation}
		The growth bounds may also be used on each fixed larger ball contained
		in the transformed domain, with the exponential loss included in
		\(A_{\mathcal D}\).
		For \(\frac{t_2}{t_1}\le |y|\le1/4\), applying Lemma
		\ref{lem:growth-estimate-gradient} at radius \(2|y|\) gives
		\begin{equation}\label{eq:prop335-source-sup-at-R}
			\sup_{B_{2|y|}}|u_{x,t_1}|
			\le CK A_{\mathcal D}^{1/2}|y|^\gamma,
			\qquad
			\sup_{B_{2|y|}}|F|
			\le CK A_{\mathcal D}^{1/2}
			\omega(2\sqrt\Lambda t_1|y|)|y|^{\gamma-1}.
		\end{equation}
		
		For \(I_1\), punctured integration by parts gives
		\begin{align*}
			&\int_{B_{2|y|}}(\Gamma(y,z)-\Gamma(0,z))\operatorname{div}F(z)\,dz\\
			&=\int_{\partial B_{2|y|}}(\Gamma(y,z)-\Gamma(0,z))F(z)\cdot\nu\,d\sigma
			-\int_{B_{2|y|}}\nabla_z(\Gamma(y,z)-\Gamma(0,z))\cdot F(z)\,dz.
		\end{align*}
		Here the inner boundary integrals about \(0\) and \(y\) vanish as in
		the construction of \(h\). The kernel bounds
		\[
		\int_{B_{2|y|}}(|y-z|^{1-n}+|z|^{1-n})\,dz\le C|y|,
		\qquad
		|\Gamma(y,z)-\Gamma(0,z)|\le C|y|^{2-n}\quad(|z|=2|y|)
		\]
		and \eqref{eq:prop335-source-sup-at-R} imply
		\begin{align*}
			&\left|\int_{B_{2|y|}}(\Gamma(y,z)-\Gamma(0,z))
			\operatorname{div}F(z)\,dz\right|\\
			&\quad\le C\bigl(|y|^{2-n}|\partial B_{2|y|}|+|y|\bigr)
			\sup_{B_{2|y|}}|F|\\
			&\quad\le CK A_{\mathcal D}^{1/2}|y|^\gamma
			\omega(2\sqrt\Lambda t_1|y|).
		\end{align*}
		With \(p'=p/(p-1)\), the potential term is bounded by
		\[
		C\|\widetilde V\|_{L^p(B_{2|y|})}
		\sup_{B_{2|y|}}|u_{x,t_1}|
		\left(\int_{B_{2|y|}}(|y-z|^{(2-n)p'}+|z|^{(2-n)p'})\,dz\right)^{1/p'}
		\le CK A_{\mathcal D}Mt_1^{2-n/p}|y|^{\gamma+2-n/p}.
		\]
		Here \((n-2)p'<n\). Consequently,
		\begin{equation}\label{eq:prop335-I1-final}
			|I_1(y)|\le CK A_{\mathcal D}|y|^\gamma
			\left[\omega(2\sqrt\Lambda t_1|y|)+Mt_1^{2-n/p}|y|^{2-n/p}\right].
		\end{equation}
		In dimension two, make the change of variables \(z=|y|\zeta\) in
		\(\Gamma(y,z)-\Gamma(0,z)\). The two \(\log |y|\) terms cancel, and the
		remaining logarithmic kernels have bounded \(L^{p'}\) norm on the
		fixed rescaled ball. This gives the same power \(|y|^{2-n/p}\).
		
		We give the annular calculation for \(I_2\), since the boundary
		terms are important for Dini coefficients. For every \(1\le k\le d\),
		\begin{align}
			&\int_{\frac{t_2}{t_1}<|z|<2|y|}\Gamma_k(z)P_{z,k}(y)\operatorname{div}F(z)\,dz
			\nonumber\\
			&=\int_{\partial B_{2|y|}}\Gamma_k(z)P_{z,k}(y)F(z)\cdot\nu\,d\sigma
			-\int_{\partial B_{\frac{t_2}{t_1}}}\Gamma_k(z)P_{z,k}(y)F(z)\cdot\nu\,d\sigma
			\nonumber\\
			&\quad-\int_{\frac{t_2}{t_1}<|z|<2|y|}
			\nabla_z(\Gamma_k(z)P_{z,k}(y))\cdot F(z)\,dz,
			\label{eq:step4-I2a-decomposition}
		\end{align}
		where \(\nu=z/|z|\) on both spheres. Equations
		\eqref{eq:Green-coefficient-bound} and \eqref{eq:prop335-source-bounds}
		bound the sum of these terms by
		\begin{align}
			CK A_{\mathcal D}^{1/2}
			\sum_{k=1}^d(4/3)^k(1+k)^{(n-1)/2}|y|^k
			\Bigl[&(2|y|)^{\gamma-k}\omega(2\sqrt\Lambda t_1|y|)
			+\left(\frac{t_2}{t_1}\right)^{\gamma-k}\omega(\sqrt\Lambda t_2)\nonumber\\
			&+\int_{\frac{t_2}{t_1}}^{2|y|}\omega(\sqrt\Lambda t_1s)s^{\gamma-k-1}\,ds\Bigr].
			\label{eq:step4-J23-raw}
		\end{align}
		Since \(k\le d=\lfloor\gamma\rfloor\),
		\[
		|y|^k\int_{\frac{t_2}{t_1}}^{2|y|}\omega(\sqrt\Lambda t_1s)s^{\gamma-k-1}\,ds
		\le2^{\gamma-k}|y|^\gamma
		\int_0^{2\sqrt\Lambda t_1|y|}\frac{\omega(s)}s\,ds.
		\]
		For every \(0<a\le\sqrt\Lambda t_1\), monotonicity gives
		\begin{equation}\label{eq:Dini-pointwise-control}
			\omega(a)\le\frac1{\log2}\int_a^{2a}\frac{\omega(s)}s\,ds
			\le\frac1{\log2}\int_0^{2\sqrt\Lambda t_1}\frac{\omega(s)}s\,ds.
		\end{equation}
		Since \(\frac{t_2}{t_1}\le |y|\), \(\gamma-k\ge0\), and \(|y|\le1/4\),
		\begin{align*}
			|y|^k(2|y|)^{\gamma-k}\omega(2\sqrt\Lambda t_1|y|)
			&\le2^{\gamma-k}|y|^\gamma
			\int_0^{2\sqrt\Lambda t_1}\frac{\omega(s)}s\,ds,\\
			|y|^k\left(\frac{t_2}{t_1}\right)^{\gamma-k}\omega(\sqrt\Lambda t_2)
			&\le |y|^\gamma
			\int_0^{2\sqrt\Lambda t_1}\frac{\omega(s)}s\,ds.
		\end{align*}
		Thus both boundary terms in \eqref{eq:step4-I2a-decomposition} have
		the same Dini bound as the volume term. Moreover,
		\[
		\sum_{k=1}^d(4/3)^k(1+k)^{(n-1)/2}2^{\gamma-k}
		\le C(n)2^\gamma\le C A_{\mathcal D}^{1/2}.
		\]
		Thus the divergence contribution to \(I_2\) is at most
		\(CK A_{\mathcal D}|y|^\gamma\int_0^{2\sqrt\Lambda t_1}\omega(s)\,ds/s\).
		For its potential contribution, H\"older's inequality gives
		\begin{equation}\label{eq:step4-I2b-radial}
			\left(\int_{\frac{t_2}{t_1}<|z|<2|y|}|z|^{(\gamma+2-n-k)p'}\,dz\right)^{1/p'}
			\le C(n,p)2^{\gamma-k}|y|^{\gamma+2-k-n/p},
		\end{equation}
		because \(\gamma+2-k-n/p\ge2-n/p>0\). Multiplying by \(|y|^k\),
		using \eqref{eq:prop335-source-bounds}, and summing in \(k\) yields
		\begin{equation}\label{eq:step4-I2-final}
			|I_2(y)|\le CK A_{\mathcal D}|y|^\gamma
			\left[\int_0^{2\sqrt\Lambda t_1}\frac{\omega(s)}s\,ds
			+Mt_1^{2-n/p}|y|^{2-n/p}\right].
		\end{equation}
		
		For \(I_3\), apply \eqref{eq:Green-remainder-bound} and integrate
		by parts on \(B_1\setminus B_{2|y|}\). The two boundary terms and the
		volume term are bounded by
		\begin{align}
			C2^{d+1}K A_{\mathcal D}^{1/2}|y|^{d+1}
			\Bigl[&\omega(\sqrt\Lambda t_1)
			+(2|y|)^{\gamma-d-1}\omega(2\sqrt\Lambda t_1|y|)\nonumber\\
			&+\int_{2|y|}^1\omega(\sqrt\Lambda t_1s)s^{\gamma-d-2}\,ds\Bigr].
			\label{eq:step5-I3a-decomposition}
		\end{align}
		Since \(0<d+1-\gamma\le1\),
		\[
		|y|^{d+1}\int_{2|y|}^1\omega(\sqrt\Lambda t_1s)s^{\gamma-d-2}\,ds
		\le |y|^\gamma\int_0^{2\sqrt\Lambda t_1}\frac{\omega(s)}s\,ds.
		\]
		For the two boundary terms, \eqref{eq:Dini-pointwise-control} and
		\(|y|\le1/4\) give
		\begin{align*}
			|y|^{d+1}\omega(\sqrt\Lambda t_1)
			&\le |y|^\gamma\int_0^{2\sqrt\Lambda t_1}\frac{\omega(s)}s\,ds,\\
			|y|^{d+1}(2|y|)^{\gamma-d-1}
			\omega(2\sqrt\Lambda t_1|y|)
			&\le |y|^\gamma\int_0^{2\sqrt\Lambda t_1}\frac{\omega(s)}s\,ds.
		\end{align*}
		Hence the coefficient estimate is uniform as \(\gamma\) approaches an
		integer.
		The potential term is bounded by
		\begin{align}
			C2^{d+1}K A_{\mathcal D}^{1/2}|y|^{d+1}
			\|\widetilde V\|_{L^p(B_1)}
			\left(\int_{2|y|}^1
			s^{p'(\gamma+1-d-n)+n-1}\,ds\right)^{1/p'}.
			\label{eq:I3b-holder}
		\end{align}
		The exponent satisfies
		\[
		p'(\gamma+1-d-n)+n-1
		=p'(\gamma+1-d-n/p)-1>-1,
		\]
		with a uniform gap, since \(\gamma-d\ge0\) and \(p>n\).
		The radial factor is bounded by \(C(n,p)\). As \(|y|^{d+1}\le |y|^\gamma\),
		we conclude that
		\begin{equation}\label{eq:step5-I3-final}
			|I_3(y)|\le C\eta(t_1)|y|^\gamma.
		\end{equation}
		Together with \eqref{eq:prop335-I1-final} and \eqref{eq:step4-I2-final},
		this proves the annular estimate for \(|y|\le1/4\).
		For \(1/4\le |y|\le3/4\), estimate the two terms in
		\eqref{eq:prop335-phi-definition} directly. Punctured integration by
		parts in the first term and the uniform kernel bounds
		\[
		\sup_{1/4\le|y|\le3/4}
		\int_{B_1}\bigl(|\nabla_z\Gamma(y,z)|
		+|\nabla_z\Gamma(0,z)|\bigr)\,dz\le C,
		\]
		\[
		\sup_{1/4\le|y|\le3/4}
		\left(\|\Gamma(y,\cdot)-\Gamma(0,\cdot)\|_{L^{p'}(B_1)}
		+\sup_{|z|=1}|\Gamma(y,z)-\Gamma(0,z)|\right)\le C
		\]
		give
		\[
		\left|\int_{B_1}(\Gamma(y,z)-\Gamma(0,z))
		\Delta u_{x,t_1}(z)\,dz\right|
		\le CK A_{\mathcal D}^{1/2}
		[\omega(\sqrt\Lambda t_1)+Mt_1^{2-n/p}].
		\]
		For the finite polynomial correction, integration by parts and
		\eqref{eq:Green-coefficient-bound} give the coefficient contribution
		at most
		\begin{align*}
			CK A_{\mathcal D}^{1/2}
			\sum_{k=1}^d(4/3)^k(1+k)^{(n-1)/2}|y|^k
			\Bigl[&\omega(\sqrt\Lambda t_1)
			+\left(\frac{t_2}{t_1}\right)^{\gamma-k}\omega(\sqrt\Lambda t_2)\\
			&+\int_{\frac{t_2}{t_1}}^1\omega(\sqrt\Lambda t_1s)
			s^{\gamma-k-1}\,ds\Bigr].
		\end{align*}
		Its potential contribution is at most
		\[
		CK A_{\mathcal D}^{1/2}Mt_1^{2-n/p}
		\sum_{k=1}^d(4/3)^k(1+k)^{(n-1)/2}|y|^k
		\left(\int_{\frac{t_2}{t_1}}^1
		s^{p'(\gamma+2-n-k)+n-1}\,ds\right)^{1/p'}.
		\]
		Here \(\gamma-k\ge0\),
		\(p'(\gamma+2-n-k)+n-1
		=p'(\gamma+2-k-n/p)-1>-1\), and \(4|y|/3\le1\).
		Using \eqref{eq:Dini-pointwise-control} in the boundary terms, both
		displays are bounded by
		\[
		CK A_{\mathcal D}^{1/2}4^{C\gamma}
		\left(\int_0^{2\sqrt\Lambda t_1}\frac{\omega(s)}s\,ds
		+Mt_1^{2-n/p}\right).
		\]
		Finally, \(|y|^\gamma\ge4^{-\gamma}\), and the leading exponent in
		\(A_{\mathcal D}\) was fixed so that
		\(A_{\mathcal D}^{1/2}4^{(C+1)\gamma}\le C A_{\mathcal D}\).
		Thus \(|\phi(y)|\le C\eta(t_1)|y|^\gamma\).
		
		It remains to prove the linear bound in the inner ball. If
		\(\frac{t_2}{2t_1}\le |y|\le \frac{t_2}{t_1}\), then \(\frac{t_2}{t_1}\le2|y|\), so
		\eqref{eq:prop335-decomposition} is still valid. The estimate for \(I_1\)
		uses the growth bounds at radius \(2|y|\). In the inner boundary terms of
		\(I_2\),
		\[
		|y|^k\left(\frac{t_2}{t_1}\right)^{\gamma-k}\le2^{\gamma-k}|y|^\gamma,
		\qquad1\le k\le d.
		\]
		The resulting factor \(2^\gamma\) is included in \(A_{\mathcal D}\), and
		the preceding estimates give
		\[
		|\phi(y)|\le C\eta(t_1)|y|^\gamma
		\le C\eta(t_1)\left(\frac{t_2}{t_1}\right)^{\gamma-1}|y|,
		\]
		where the last inequality follows from \(\gamma\ge1-\varepsilon_1\) and
		\((|y|t_1/t_2)^{\gamma-1}\le2^{\varepsilon_1}\).
		For \(0<|y|<\frac{t_2}{2t_1}\), the correct decomposition is
		\begin{equation}\label{eq:step7-correct-decomposition}
			\phi(y)=J_1(y)+J_2(y)+J_3(y),
		\end{equation}
		where
		\begin{align*}
			J_1(y)&=\int_{|z|<2|y|}(\Gamma(y,z)-\Gamma(0,z))
			\Delta u_{x,t_1}(z)\,dz,\\
			J_2(y)&=\int_{2|y|<|z|<\frac{t_2}{t_1}}(\Gamma(y,z)-\Gamma(0,z))
			\Delta u_{x,t_1}(z)\,dz,\\
			J_3(y)&=\int_{\frac{t_2}{t_1}<|z|<1}
			\left(\Gamma(y,z)-\sum_{k=0}^d\Gamma_k(z)P_{z,k}(y)\right)
			\Delta u_{x,t_1}(z)\,dz.
		\end{align*}
		On \(B_{2|y|}\), Lemma \ref{lem:growth-estimate-gradient} gives
		\(\sup|\nabla u_{x,t_1}|\le CK A_{\mathcal D}^{1/2}\left(\frac{t_2}{t_1}\right)^{\gamma-1}\)
		and \(|u_{x,t_1}(z)|\le CK A_{\mathcal D}^{1/2}\left(\frac{t_2}{t_1}\right)^{\gamma-1}|z|\).
		The punctured integration-by-parts identity for \(I_1\), now with these
		bounds, gives
		\begin{align*}
			|J_1(y)|
			&\le C |y|\sup_{B_{2|y|}}|F|
			+C\|\widetilde V\|_{L^p(B_{2|y|})}
			\sup_{B_{2|y|}}|u_{x,t_1}|\,|y|^{2-n/p}\\
			&\le CK A_{\mathcal D}^{1/2}\left(\frac{t_2}{t_1}\right)^{\gamma-1}|y|
			[\omega(2\sqrt\Lambda t_1|y|)+Mt_1^{2-n/p}|y|^{2-n/p}].
		\end{align*}
		On the middle annulus the mean-value theorem gives
		\begin{equation}\label{eq:step7-kernel-middle}
			\begin{split}
				|\Gamma(y,z)-\Gamma(0,z)|&\le C|y||z|^{1-n},\\
				|\nabla_z(\Gamma(y,z)-\Gamma(0,z))|&\le C|y||z|^{-n}.
			\end{split}
		\end{equation}
		Integration by parts, including both boundary spheres, bounds the
		coefficient contribution to \(J_2\) by
		\[
		CK A_{\mathcal D}\left(\frac{t_2}{t_1}\right)^{\gamma-1}|y|
		\left[\omega(2\sqrt\Lambda t_1|y|)+\omega(\sqrt\Lambda t_2)
		+\int_{2|y|}^{\frac{t_2}{t_1}}\frac{\omega(\sqrt\Lambda t_1s)}s\,ds\right].
		\]
		For the potential contribution, the extra factor \(|z|\) in the
		inner bound for \(u_{x,t_1}\) gives
		\[
		CK A_{\mathcal D}^{1/2}\left(\frac{t_2}{t_1}\right)^{\gamma-1}|y|
		\|\widetilde V\|_p
		\left(\int_{B_{\frac{t_2}{t_1}}}|z|^{(2-n)p'}\,dz\right)^{1/p'}
		\le CK A_{\mathcal D}Mt_1^{2-n/p}\left(\frac{t_2}{t_1}\right)^{\gamma-1}|y|.
		\]
		Both estimates are bounded by \(C\eta(t_1)\left(\frac{t_2}{t_1}\right)^{\gamma-1}|y|\).
		For \(J_3\), put
		\[
		\mathcal R(y,z)=\Gamma(y,z)-\sum_{k=0}^d
		\Gamma_k(z)P_{z,k}(y).
		\]
		Since \(|y|<\frac{t_2}{2t_1}\), \eqref{eq:Green-remainder-bound} gives, for
		\(|z|\ge \frac{t_2}{t_1}\),
		\[
		|\mathcal R(y,z)|\le C2^{d+1}|y|^{d+1}|z|^{1-n-d},
		\qquad
		|\nabla_z\mathcal R(y,z)|\le C2^{d+1}|y|^{d+1}|z|^{-n-d}.
		\]
		Integration by parts on \(B_1\setminus B_{\frac{t_2}{t_1}}\), including both
		boundary spheres, and \eqref{eq:prop335-source-bounds} yield
		\begin{align*}
			&\left|\int_{B_1\setminus B_{\frac{t_2}{t_1}}}
			\mathcal R(y,z)\operatorname{div}F(z)\,dz\right|\\
			&\quad\le C2^{d+1}K A_{\mathcal D}^{1/2}|y|^{d+1}
			\Bigl[\omega(\sqrt\Lambda t_1)
			+\left(\frac{t_2}{t_1}\right)^{\gamma-d-1}\omega(\sqrt\Lambda t_2)\\
			&\hspace{49mm}
			+\int_{\frac{t_2}{t_1}}^1\omega(\sqrt\Lambda t_1s)
			s^{\gamma-d-2}\,ds\Bigr].
		\end{align*}
		Since \(0<d+1-\gamma\le1\),
		\[
		\int_{\frac{t_2}{t_1}}^1\omega(\sqrt\Lambda t_1s)s^{-1-(d+1-\gamma)}\,ds
		\le \left(\frac{t_2}{t_1}\right)^{-(d+1-\gamma)}
		\int_0^{2\sqrt\Lambda t_1}\frac{\omega(s)}s\,ds.
		\]
		Equation \eqref{eq:Dini-pointwise-control} bounds the two boundary terms.
		Since \(2^{d+1}A_{\mathcal D}^{1/2}\le C A_{\mathcal D}\), the
		coefficient contribution is at most
		\[
		CK A_{\mathcal D}|y|^{d+1}
		\left[1+\left(\frac{t_2}{t_1}\right)^{\gamma-d-1}\right]
		\int_0^{2\sqrt\Lambda t_1}\frac{\omega(s)}s\,ds
		\]
		For the potential contribution, H\"older's inequality gives
		\begin{align*}
			&\left|\int_{B_1\setminus B_{\frac{t_2}{t_1}}}
			\mathcal R(y,z)\widetilde V(z)u_{x,t_1}(z)\,dz\right|\\
			&\quad\le C2^{d+1}K A_{\mathcal D}^{1/2}|y|^{d+1}
			Mt_1^{2-n/p}
			\left(\int_{\frac{t_2}{t_1}}^1
			s^{p'(\gamma+1-n-d)+n-1}\,ds\right)^{1/p'}\\
			&\quad\le CK A_{\mathcal D}Mt_1^{2-n/p}|y|^{d+1},
		\end{align*}
		because
		\(p'(\gamma+1-n-d)+n-1
		=p'(\gamma+1-d-n/p)-1>-1\).
		The inequalities
		\[
		|y|^{d+1}\left(\frac{t_2}{t_1}\right)^{\gamma-d-1}\le |y|\left(\frac{t_2}{t_1}\right)^{\gamma-1},
		\qquad |y|^{d+1}\le |y|\left(\frac{t_2}{t_1}\right)^d
		\le |y|\left(\frac{t_2}{t_1}\right)^{\gamma-1}
		\]
		therefore give the same linear bound. This proves
		\eqref{eq:prop335-DiniLp-inner} in all of \(B_{\frac{t_2}{t_1}}\).
	\end{proof}
	
	Let \(c_2>0\) be a sufficiently small structural constant. Now we give more approximation relations as a corollary of Proposition \ref{prop:quantitative-harmonic-approximation-DiniLp}.
	
	\begin{corollary}\label{cor:relative-harmonic-approximation-DiniLp}
		Under the hypotheses and notation of Proposition
		\ref{prop:quantitative-harmonic-approximation-DiniLp}, let \(h\) be
		the harmonic function constructed there. If
		\(0<q_0\le \frac{t_2}{t_1}\le1/10\), \(0<\varepsilon\le1/100\), and
		\begin{equation}\label{eq:prop335-DiniLp-r0-condition}
			A_{\mathcal D}^{1/2}q_0^{-\varepsilon_1}\eta(t_1)
			\le c_2\varepsilon,
		\end{equation}
		then, for \(\frac{t_2}{t_1}\le t\le1/2\),
		\begin{align}
			\|h-u_{x,t_1}\|_{L^2(\partial B_t,\vint)}
			&\le\varepsilon \left(\vint_{\partial B_{t}}u_{x,t_1}^2\right)^{1/2},
			\label{eq:prop335-DiniLp-paper-form}\\
			(1-\varepsilon)\left(\vint_{\partial B_{t}}u_{x,t_1}^2\right)^{1/2}
			&\le \left(\vint_{\partial B_{t}}h^2\right)^{1/2}\le(1+\varepsilon)\left(\vint_{\partial B_{t}}u_{x,t_1}^2\right)^{1/2},
			\label{eq:prop335-DiniLp-spherical-norm-comparison}\\
			\vint_{\partial B_1}|u_{x,tt_1}-h_{0,t}|^2&\le4\varepsilon^2.
			\label{eq:prop335-DiniLp-rescaling-comparison}
		\end{align}
		Moreover,
		\begin{equation}\label{eq:prop335-DiniLp-doubling}
			|D^h(0,t)-D^u(x,tt_1)|\le10\varepsilon
			\qquad(\frac{t_2}{t_1}\le t\le1/4).
		\end{equation}
	\end{corollary}
	
	\begin{proof}
		For \(\frac{t_2}{t_1}\le t\le1/2\), choose \(j\) with \(1\le2^jt<2\).
		The harmonic approximation at scale \(t_1\) supplies a normalized
		harmonic function \(g\) with
		\(\|g-u_{x,t_1}\|_{L^\infty(B_2)}\le1/4\). Since harmonic spherical
		means are nondecreasing,
		\[
		\left(\vint_{\partial B_{s}}u_{x,t_1}^2\right)^{1/2}
		\ge \left(\vint_{\partial B_{s}}g^2\right)^{1/2}-\frac14
		\ge \left(\vint_{\partial B_{1}}g^2\right)^{1/2}-\frac14
		\ge\frac12,
		\qquad1\le s\le2.
		\]
		For \(0\le k<j\), the corresponding radius \(2^ktt_1\) belongs to
		\([t_2,t_1]\), and therefore
		\(D^{u_{x,t_1}}(0,2^kt)\le\gamma+\varepsilon_1\). Thus
		\begin{equation}\label{eq:lower-bound-normalized-spherical-average}
			\begin{split}
				\left(\vint_{\partial B_{t}}u_{x,t_1}^2\right)
				&=\left(\vint_{\partial B_{2^jt}}u_{x,t_1}^2\right)
				\prod_{k=0}^{j-1}4^{-D^{u_{x,t_1}}(0,2^kt)}\\
				&\ge\frac14\,4^{-j(\gamma+\varepsilon_1)},
			\end{split}
		\end{equation}
		and \(2^{-j}>t/2\) gives
		\[
		\left(\vint_{\partial B_{t}}u_{x,t_1}^2\right)^{1/2}
		\ge c t^{\gamma+\varepsilon_1}
		\ge cA_{\mathcal D}^{-1/2}t^{\gamma+\varepsilon_1}.
		\]
		Consequently,
		\begin{equation}\label{eq:step8-relative-ratio}
			\frac{\|h-u_{x,t_1}\|_{L^2(\partial B_t,\vint)}}
			{\left(\vint_{\partial B_{t}}u_{x,t_1}^2\right)^{1/2}}
			\le C A_{\mathcal D}^{1/2}\eta(t_1)t^{-\varepsilon_1}
			\le C A_{\mathcal D}^{1/2}q_0^{-\varepsilon_1}\eta(t_1)
			\le\varepsilon.
		\end{equation}
		The triangle inequality proves the spherical norm comparison.
		For nonzero vectors in a Hilbert space we have
		\begin{equation}\label{eq:center-comparison-normalization-inequality}
			\left\|\frac f{\|f\|}-\frac g{\|g\|}\right\|
			\le\frac{2\|f-g\|}{\|g\|}.
		\end{equation}
		Applying this inequality with the spherical norm
		of \(h\) and \(u_{x,t_1}\), and using the identity
		\[
		u_{x,tt_1}(\xi)=\frac{u_{x,t_1}(t\xi)}{\left(\vint_{\partial B_{t}}u_{x,t_1}^2\right)^{1/2}},
		\]
		proves \eqref{eq:prop335-DiniLp-rescaling-comparison}. At both \(t\)
		and \(2t\), the spherical norms differ by factors in
		\([1-\varepsilon,1+\varepsilon]\). Hence, for \(t\le1/4\),
		\[
		|D^h(0,t)-D^u(x,tt_1)|
		\le2\log_4\frac{1+\varepsilon}{1-\varepsilon}\le10\varepsilon.
		\]
		This proves all the relative comparisons.
	\end{proof}
	
	With Corollary \ref{cor:relative-harmonic-approximation-DiniLp} in hand, we are going to prove the main result in this subsection. Before that, we cite the following quantitative uniqueness of tangent maps defined in \eqref{eq:tangent-spherical-normalization} for harmonic function.
	\begin{proposition}\label{prop:uniform-symmetry-under-pinching}
		Let \(h\) be harmonic on \(B_{2t_1}\), and let
		\(0<20t_2\le t_1\). If
		\(D^h(0,t_1)-D^h(0,t_2)\le\varepsilon\), with
		\(\varepsilon\le\varepsilon_0(n)\), there is a degree \(d\ge0\)
		such that the normalized degree-\(d\) term \(P\) of the expansion of
		\(h\) at zero satisfies
		\[
		\vint_{\partial B_1}|h_{0,t}-P|^2\le C\varepsilon
		\qquad(3t_2\le t\le t_1/3).
		\]
	\end{proposition}
	This is \cite[Proposition 2.26 and Remark 2.27]{HuangJiang}, obtained by concentrating
	the spherical harmonic coefficients in one degree. The same proof
	includes degree zero when the constant term is retained in the
	normalization, as we do here.
	
	Now we are ready to prove the quantitative uniqueness of tangent map.
	\begin{theorem}\label{thm:uniform-symmetry-pinched}
		Let \(x\in Z(u)\cap B_1\), \(0<q_0\le t_2/t_1\le10^{-2}\), and
		\(t_1\le\min\{r_4(1/10),(10\sqrt\Lambda)^{-1}\}\).
		There are structural constants \(\varepsilon_2,c_3>0\)
		such that, if \(\varepsilon\le\varepsilon_2\),
		\begin{equation}\label{eq:symmetry-smallness}
			A_{\mathcal D}^{1/2}q_0^{-\varepsilon_1}\eta(t_1)
			\le c_3\varepsilon,
			\qquad |D^u(x,s)-D^u(x,t_1)|\le\varepsilon\quad(t_2\le s\le t_1),
		\end{equation}
		then one normalized homogeneous harmonic polynomial \(P\) satisfies
		\begin{equation}\label{eq:symmetry-conclusion}
			\vint_{\partial B_1}|u_{x,t}-P|^2\le C\varepsilon
			\qquad(4t_2\le t\le t_1/12).
		\end{equation}
		The polynomial is the normalized homogeneous term of the harmonic
		function in Proposition \ref{prop:quantitative-harmonic-approximation-DiniLp}.
	\end{theorem}
	\begin{proof}
		Choose \(c_3\le c_2/100\) and \(\varepsilon_2\le\varepsilon_1/2\), where
		\(c_2\) is the constant fixed before Corollary
		\ref{cor:relative-harmonic-approximation-DiniLp} and
		\(\varepsilon_1\) is fixed before Proposition
		\ref{prop:quantitative-harmonic-approximation-DiniLp}. The hypothesis
		gives
		\[
		A_{\mathcal D}^{1/2}q_0^{-\varepsilon_1}\eta(t_1)
		\le c_3\varepsilon=c_2\frac{c_3}{c_2}\varepsilon,
		\qquad
		\mathop{\mathrm{osc}}_{[t_2,t_1]}D^u(x,\cdot)
		\le2\varepsilon\le\varepsilon_1.
		\]
		Thus Proposition \ref{prop:quantitative-harmonic-approximation-DiniLp}
		and Corollary \ref{cor:relative-harmonic-approximation-DiniLp} apply with
		comparison error \(c_3\varepsilon/c_2\). The same harmonic function \(h\) satisfies
		\[
		|D^h(0,s)-D^u(x,t_1s)|\le10\frac{c_3}{c_2}\varepsilon
		\qquad(t_2/t_1\le s\le1/4).
		\]
		For \(t_2/t_1\le s_1,s_2\le1/4\),
		\[
		|D^h(0,s_1)-D^h(0,s_2)|
		\le20\frac{c_3}{c_2}\varepsilon
		+|D^u(x,t_1s_1)-D^u(x,t_1s_2)|
		\le\left(2+20\frac{c_3}{c_2}\right)\varepsilon.
		\]
		Decrease \(\varepsilon_2\) so that this is at most
		\(\varepsilon_0(n)\), and apply Proposition
		\ref{prop:uniform-symmetry-under-pinching} with outer radius \(1/4\)
		and inner radius \(t_2/t_1\); the ratio is at least \(25\).
		One normalized homogeneous term \(P\) then satisfies
		\[
		\|h_{0,s}-P\|_{L^2(\partial B_1,\vint)}
		\le C\sqrt\varepsilon
		\qquad(3t_2/t_1\le s\le1/12).
		\]
		The relative comparison with \(u_{x,t_1}\) holds on
		\([t_2/t_1,1/2]\). For
		\(4t_2/t_1\le s\le1/12\),
		\[
		\|u_{x,t_1s}-P\|^2
		\le2\|u_{x,t_1s}-h_{0,s}\|^2
		+2\|h_{0,s}-P\|^2
		\le C\varepsilon.
		\]
		This is \eqref{eq:symmetry-conclusion}. The homogeneous term is chosen from
		the expansion of this one function \(h\), so it is independent of
		the intermediate radius.
	\end{proof}
	
	\begin{lemma}\label{lem:integer-pinching-drop}
		There is a sufficiently small structural constant \(c>0\) such that,
		for every sufficiently small structural \(\varepsilon>0\) and
		\(x\in Z(u)\cap B_1\), the following two assertions hold. 
		\begin{itemize}
			\item 	If
			\(q_0\le t_2/t_1\le1/400\),
			\(t_1\le\min\{r_4(\varepsilon/100),(10\sqrt\Lambda)^{-1}\}\),
			\(A_{\mathcal D}^{1/2}q_0^{-\varepsilon_1}\eta(t_1)
			\le c\varepsilon\), and
			\(|D^u(x,t_2)-D^u(x,t_1)|\le\varepsilon\), then a positive integer
			\(d\) satisfies
			\begin{equation}\label{eq:degree-drop-integer-conclusion}
				|D^u(x,s)-d|\le C\varepsilon\qquad(4t_2\le s\le t_1/16).
			\end{equation}
			\item Let \(d\) be a positive integer. If
			\(t_1\le r_4(\varepsilon/100)\), Lemma
			\ref{lem:comparison-harmonic} is valid at scale \(t_1\) with error
			\(\varepsilon/100\), and \(D^u(x,t_1)\le d-\varepsilon\), then
			\begin{equation}\label{eq:degree-drop-conclusion}
				D^u(x,s)\le d-1+\varepsilon\qquad(0<s\le\varepsilon t_1/8).
			\end{equation}
		\end{itemize}
	\end{lemma}
	\begin{proof}
		Apply Proposition
		\ref{prop:quantitative-harmonic-approximation-DiniLp} and Corollary
		\ref{cor:relative-harmonic-approximation-DiniLp} with outer scale
		\(t_1/2\), inner scale \(2t_2\), and comparison error
		\(\varepsilon/100\). Their scale ratio is \(4t_2/t_1\), and
		\[
		4q_0\le\frac{4t_2}{t_1}\le\frac1{100},
		\qquad
		A_{\mathcal D}^{1/2}\left(\frac{4t_2}{t_1}\right)^{-\varepsilon_1}\eta(t_1/2)
		\le A_{\mathcal D}^{1/2}q_0^{-\varepsilon_1}\eta(t_1)
		\le c\varepsilon.
		\]
		Moreover, for \(2t_2\le s\le t_1/2\), Theorem
		\ref{thm:almost-monotonicity} gives
		\[
		D^u(x,t_2)-\frac\varepsilon{100}
		\le D^u(x,s)
		\le D^u(x,t_1)+\frac\varepsilon{100},
		\]
		so the oscillation on this interval is at most
		\(2\varepsilon+2\varepsilon/100\). Taking the structural upper bound
		on \(\varepsilon\) sufficiently small and decreasing \(c\) if necessary
		verifies both
		hypotheses. The resulting harmonic function satisfies
		\begin{equation}\label{eq:degree-drop-h-u-comparison}
			\left|D^h(0,\tau)-D^u\left(x,\frac{\tau t_1}{2}\right)\right|
			\le\frac{\varepsilon}{10}
			\quad\left(\frac{4t_2}{t_1}\le\tau\le\frac14\right).
		\end{equation}
		Hence \(D^h(0,\cdot)\) varies by at most \(C\varepsilon\) on this
		interval. Apply Lemma \ref{lem:pinching-integer} to
		\(y\mapsto h((80t_2/t_1)y)\). Its radii \(1/20\) and \(1\)
		correspond to \(4t_2/t_1\) and \(80t_2/t_1\), both in the comparison
		interval. It gives an integer \(d\ge0\) with
		\[
		|D^h(0,\tau)-d|\le C\varepsilon
		\quad\left(\frac{8t_2}{t_1}\le\tau\le
		\frac{40t_2}{t_1}\right).
		\]
		The preceding oscillation estimate extends the same bound to
		\([8t_2/t_1,1/8]\). Taking \(\tau=2s/t_1\) in
		\eqref{eq:degree-drop-h-u-comparison} proves
		\eqref{eq:degree-drop-integer-conclusion}. The lower bound
		\(D^u\ge1-\varepsilon/100\) excludes degree zero.
		
		For the second assertion, let \(h\) be the harmonic function supplied
		by Lemma \ref{lem:comparison-harmonic} at scale \(t_1\). Then
		\[
		|D^h(0,\tau)-D^u(x,\tau t_1)|\le\frac{\varepsilon}{100}
		\qquad\left((\varepsilon/100)^2\le\tau\le1\right).
		\]
		In particular,
		\[
		D^h(0,1)\le D^u(x,t_1)+\frac{\varepsilon}{100}
		\le d-\frac{99\varepsilon}{100}
		\le d-\frac\varepsilon2.
		\]
		Applying \eqref{eq:drop-conclusion} with parameter \(\varepsilon/2\)
		gives
		\[
		D^h(0,\varepsilon/4)\le d-1+\frac\varepsilon2.
		\]
		Since \((\varepsilon/100)^2\le\varepsilon/4\le1\), comparison at that
		radius gives
		\[
		D^u(x,\varepsilon t_1/4)\le d-1+51\varepsilon/100.
		\]
		If \(0<s\le\varepsilon t_1/8\), then
		\[
		2s\le\frac{\varepsilon t_1}{4}\le t_1
		\le r_4(\varepsilon/100).
		\]
		Thus Theorem \ref{thm:almost-monotonicity}, with error
		\(\varepsilon/100\), small radius \(s\), and outer radius
		\(\varepsilon t_1/4\), gives
		\[
		D^u(x,s)
		\le D^u(x,\varepsilon t_1/4)+\frac\varepsilon{100}
		\le d-1+\frac{52\varepsilon}{100}
		\le d-1+\varepsilon.
		\]
		This proves \eqref{eq:degree-drop-conclusion}.
	\end{proof}
	
	\begin{proposition}\label{prop:low-doubling-regularity}
		There are structural constants \(0<q_0\le10^{-4}\),
		\(\varepsilon_3,c_5>0\) such that
		\(x\in Z(u)\cap B_1\) is regular $\left(|\nabla u(x)|>0\right)$ whenever following five conditions hold: 
		\begin{equation}\label{eq:regularity-scale-smallness}
			\begin{gathered}
				0<\varepsilon\le\varepsilon_3,\qquad
				r\le r_4(\varepsilon),\qquad
				r/2\le\min\{r_4(1/10),(10\sqrt\Lambda)^{-1}\},\\
				A_{\mathcal D}^{1/2}q_0^{-\varepsilon_1}\eta(r/2)
				\le c_5\varepsilon,\qquad D^u(x,r)\le1+\varepsilon.
			\end{gathered}
		\end{equation}
	\end{proposition}
	\begin{proof}
		Almost monotonicity gives
		\begin{equation}\label{eq:regularity-near-one}
			1-\varepsilon\le D^u(x,s)\le1+2\varepsilon
			\qquad(q_0r/2\le s\le r/2).
		\end{equation}
		Let \(h\) be the harmonic function in Proposition
		\ref{prop:quantitative-harmonic-approximation-DiniLp}, and apply
		Corollary \ref{cor:relative-harmonic-approximation-DiniLp} with
		comparison error \(c_5\varepsilon/c_2\). This is permitted because
		\[
		A_{\mathcal D}^{1/2}q_0^{-\varepsilon_1}\eta(r/2)
		\le c_5\varepsilon=c_2\frac{c_5}{c_2}\varepsilon
		\]
		after choosing \(c_5\le c_2/100\), and
		\(\mathop{\mathrm{osc}}_{[q_0r/2,r/2]}D^u(x,\cdot)\le3\varepsilon
		\le\varepsilon_1\). Its estimate
		inside the smaller radius is
		\begin{equation}\label{eq:regularity-inner-comparison}
			|h(y)-u_{x,r/2}(y)|
			\le C\eta(r/2)q_0^{\gamma-1}|y|
			\qquad(|y|\le q_0),
			\qquad\gamma\ge1-\varepsilon_1.
		\end{equation}
		Both functions vanish at zero and are differentiable there.
		Taking difference quotients along any unit direction yields
		\begin{equation}\label{eq:regularity-gradient-comparison}
			|\nabla h(0)-\nabla u_{x,r/2}(0)|
			\le C\eta(r/2)q_0^{-\varepsilon_1}
			\le Cc_5\varepsilon.
		\end{equation}
		
		Take \(q_0\le1/100\). For
		\(q_0\le s\le1/4\), the relative doubling comparison and
		\eqref{eq:regularity-near-one} give
		\[
		|D^h(0,s)-1|
		\le |D^h(0,s)-D^u(x,sr/2)|
		+|D^u(x,sr/2)-1|
		\le10\frac{c_5}{c_2}\varepsilon+2\varepsilon.
		\]
		Apply Proposition \ref{prop:uniform-symmetry-under-pinching} with inner
		radius \(q_0\) and outer radius \(1/4\). Its homogeneous degree is the
		integer within \(C\varepsilon\) of \(D^h\) on this interval; the displayed
		bound forces that integer to be one when \(\varepsilon_3\) is small.
		There is consequently a normalized linear polynomial \(P\) such that
		\[
		\|h_{0,1/16}-P\|_{L^2(\partial B_1,\vint)}
		\le C\sqrt\varepsilon.
		\]
		The harmonic function \(h_{0,1/16}-P\) satisfies
		\[
		|\nabla h_{0,1/16}(0)-\nabla P|
		\le C\|h_{0,1/16}-P\|_{L^2(\partial B_1,\vint)}
		\le C\sqrt\varepsilon.
		\]
		Since a normalized linear polynomial has \(|\nabla P|=\sqrt n\),
		\begin{equation}\label{eq:regularity-normalized-gradient-lower}
			|\nabla h_{0,1/16}(0)|\ge\sqrt n-C\sqrt\varepsilon
			\ge\sqrt n/2.
		\end{equation}
		It remains to compare this normalized gradient with \(\nabla h(0)\).
		Using
		\[
		\vint_{\partial B_{2s}}u_{x,r/2}^2
		=4^{D^u(x,sr/2)}\vint_{\partial B_s}u_{x,r/2}^2
		\]
		at \(s=2^j/16\), \(j=0,1,2,3\), gives
		\begin{equation}\label{eq:regularity-spherical-average-lower}
			\begin{split}
				\left(\vint_{\partial B_{1/16}}u_{x,r/2}^2\right)
				&=\left(\vint_{\partial B_{1}}u_{x,r/2}^2\right)
				4^{-\sum_{j=0}^3D^u(x,2^jr/32)}\\
				&\ge4^{-4(1+2\varepsilon)}\left(\vint_{\partial B_{1}}u_{x,r/2}^2\right)
				=4^{-4(1+2\varepsilon)}.
			\end{split}
		\end{equation}
		The relative approximation implies
		\(\left(\vint_{\partial B_{1/16}}h^2\right)^{1/2}\ge(1-C\varepsilon)
		\left(\vint_{\partial B_{1/16}}u_{x,r/2}^2\right)^{1/2}\ge c>0\). Thus
		\begin{equation}\label{eq:regularity-grad-h-lower}
			|\nabla h(0)|
			=16\left(\vint_{\partial B_{1/16}}h^2\right)^{1/2}
			|\nabla h_{0,1/16}(0)|\ge c>0.
		\end{equation}
		Choose \(\varepsilon_3\) and \(c_5\)
		small relative to this lower bound. Equation
		\eqref{eq:regularity-gradient-comparison} then implies
		\(\nabla u_{x,r/2}(0)\ne0\). Since
		\[
		\nabla u_{x,r/2}(0)
		=\frac{\frac r2A_x^T\nabla u(x)}
		{\left(\vint_{\partial B_1}u(x+\frac r2A_xz)^2\,dz\right)^{1/2}},
		\]
		the invertibility of \(A_x\) proves \(\nabla u(x)\ne0\).
	\end{proof}
	
	Let \(q_0,\varepsilon_3,c_5\) be the constants supplied by Proposition
	\ref{prop:low-doubling-regularity}; let \(\varepsilon_1\) be the constant
	fixed before Proposition
	\ref{prop:quantitative-harmonic-approximation-DiniLp}; let \(r_4\) be
	defined in \eqref{eq:almost-monotonicity-radius}; and let
	\(A_{\mathcal D}\) and \(\eta\) be defined in
	\eqref{eq:42-AD-definition} and \eqref{eq:prop335-Dini-smallness}.
	For later use define
	\begin{equation}\label{eq:regularity-radius-definition}
		\begin{split}
			r_5=\sup\Bigl\{0<R\le(10\sqrt\Lambda)^{-1}:\;&
			R\le r_4(\varepsilon_3/10),\quad
			R/2\le r_4(1/10),\\
			&A_{\mathcal D}^{1/2}q_0^{-\varepsilon_1}\eta(R/2)
			\le c_5\varepsilon_3/10\Bigr\}.
		\end{split}
	\end{equation}
	This is positive by the Dini condition, and
	\eqref{eq:regularity-scale-smallness} holds with
	\(\varepsilon=\varepsilon_3/10\) at every smaller scale.

	\section{Cone splitting}\label{sec:cone-splitting}
	
	For \(d\in\{0,1,\ldots\}\), denote the $d$-homogenous harmonic polynomial space as
	\[
	\mathcal P_d
	=\left\{P\in\mathbb R[y_1,\ldots,y_n]:
	\Delta P=0,\quad P(ty)=t^dP(y)
	\text{ for every \(t>0\) and \(y\in\mathbb R^n\)}\right\}.
	\]
	And recall that the word normalized has the meaning fixed in
	\eqref{eq:standing-spherical-normalization}.
	We first give the definition for symmetry of functions used in this section.
	\begin{definition}\label{def:quantitative-symmetry}
		\begin{itemize}
			\item Let \(V\subset\mathbb R^n\) be a
			\(k\)-dimensional linear subspace, and \(f(x)\) be a function. We
			say that \(f(x)\) is
			\(k\)-symmetric with respect to \(V\) if it is invariant under every
			translation in \(V\), namely,
			\begin{equation}\label{eq:polynomial-k-symmetry}
				f(y+v)=f(y)\qquad
				\text{for every \(y\in\mathbb R^n\) and \(v\in V\)}.
			\end{equation}
			We call \(f(x)\) \(k\)-symmetric if such a \(k\)-dimensional subspace
			exists.
			
			\item Let \(x\in Z(u)\cap B_1\), and let \(r>0\) be such that
			\(u_{x,r}\) in \eqref{eq:uxr-definition} is defined. For
			\(\varepsilon\ge0\), we say that \(u\) is
			\((k,\varepsilon,r,x)\)-symmetric with respect to \(V\) if there are
			\(d\in\{0,1,\ldots\}\) and a normalized polynomial
			\(P\in\mathcal P_d\), \(k\)-symmetric with respect to \(V\), such that
			\begin{equation}\label{eq:quantitative-symmetry-definition}
				\vint_{\partial B_1}|u_{x,r}-P|^2\le\varepsilon.
			\end{equation}
			\item Let \(I\subset(0,\infty)\) be an interval such that \(u_{x,t}\) is
			defined for every \(t\in I\). We say
			that the \((k,\varepsilon)\)-symmetry of \(u\) at \(x\) with respect
			to \(V\) is uniform on \(I\) if one degree \(d\) and one normalized
			polynomial \(P\in\mathcal P_d\), \(k\)-symmetric with respect to
			\(V\), satisfy
			\begin{equation}\label{eq:uniform-quantitative-symmetry-definition}
				\vint_{\partial B_1}|u_{x,t}-P|^2\le\varepsilon
				\qquad\text{for every \(t\in I\)}.
			\end{equation}
		\end{itemize}
		
	\end{definition}
	
	We next recall the harmonic pinching estimates, saying two different pinched points give an almost translation invariant direction.
	
	\begin{lemma}\label{lem:two-center-coefficient-estimate}
		Let \(d\ge2\), let \(h\) be harmonic in \(B_{20d}\), and let
		\(x\in B_1\). Suppose
		\[
		|D^h(a,s)-d|\le\varepsilon
		\quad(a=0,x,\quad1/10\le s\le9d),
		\qquad \varepsilon\le\varepsilon_0(n).
		\]
		If \(P_{a,d}\) is the degree-\(d\) term of the expansion at \(a\),
		then \(P_{0,d}\ne0\) and
		\begin{equation}\label{eq:Pxd-P0d-estimate}
			\|P_{x,d}-P_{0,d}\|\le C\sqrt\varepsilon|x|\|P_{0,d}\|,
			\qquad
			\|x\cdot\nabla P_{0,d}\|\le C\sqrt\varepsilon\|\nabla P_{0,d}\|.
		\end{equation}
	\end{lemma}
	\begin{proof}
		This is the coefficient calculation in
		\cite[Lemma~4.5(1)--(2)]{HuangJiang}; see also
		\cite[Lemma~3.22(1), (3)]{NaberValtorta2017}. We include the
		estimates needed below. Write
		\[
		h(a+y)=\sum_{m=0}^{\infty}P_{a,m}(y),
		\qquad a=0,x.
		\]
		Orthogonality of the homogeneous harmonic terms gives
		\[
		N^h(a,s)=
		\frac{\sum_{m\ge0}m s^{2m}\|P_{a,m}\|^2}
		{\sum_{m\ge0}s^{2m}\|P_{a,m}\|^2}.
		\]
		By \eqref{eq:frequency-doubling-comparison},
		\(N^h(a,1/5)\ge d-\varepsilon\) and
		\(N^h(a,9d)\le d+\varepsilon\).
		After expanding the two frequency quotients and using
		\(\varepsilon\le1/2\), we obtain
		\begin{align*}
			\sum_{m<d}5^{2(d-m)}\|P_{a,m}\|^2
			&\le2\varepsilon\|P_{a,d}\|^2\\
			&\quad+2\sup_{k\ge1}(k+1)(45d)^{-2k}
			\sum_{m>d}(9d)^{2(m-d)}\|P_{a,m}\|^2,\\
			\sum_{m>d}(9d)^{2(m-d)}\|P_{a,m}\|^2
			&\le2\varepsilon\|P_{a,d}\|^2\\
			&\quad+2\sup_{k\ge1}(k+1)(45d)^{-2k}
			\sum_{m<d}5^{2(d-m)}\|P_{a,m}\|^2.
		\end{align*}
		The coefficient
		\(2\sup_{k\ge1}(k+1)(45d)^{-2k}\) is smaller than \(1/2\) for
		\(d\ge2\).
		Adding and absorbing proves
		\begin{equation}\label{eq:two-center-spectral-concentration}
			\sum_{m<d}5^{2(d-m)}\|P_{a,m}\|^2
			+\sum_{m>d}(9d)^{2(m-d)}\|P_{a,m}\|^2
			\le C\varepsilon\|P_{a,d}\|^2,
			\qquad a=0,x.
		\end{equation}
		In particular, \(P_{a,d}\ne0\); otherwise
		\eqref{eq:two-center-spectral-concentration} would force every
		homogeneous term \(P_{a,m}\) to vanish.
		
		For \(Q_{d+k}\in\mathcal P_{d+k}\), we use $[Q_{d+k}]_d$ to represent its homogeneous component of degree \(d\) in Taylor expansion. Taylor's formula and
		Lemma \ref{lem:grad-inner-product-hhp} give
		\begin{align}
			[Q_{d+k}(x+\cdot)]_d
			&=\frac{(x\cdot\nabla)^kQ_{d+k}}{k!}=\sum_{|\beta|=k} \frac{x^\beta}{\beta!}\partial^{\beta}Q_{d+k}(\cdot),\notag\\
			\|[Q_{d+k}(x+\cdot)]_d\|
			&\le C(n)2^{k/2}\binom{d+k}{k}
			\left(\frac{d+k}{d}\right)^{(n-2)/4}
			|x|^k\|Q_{d+k}\|.\label{eq:two-center-translation-bound}
		\end{align}
		For the summation in \eqref{eq:two-center-translation-bound}, we have an upper bound independent of degree:
		\begin{equation}\label{eq:two-center-translation-sum}
			\sum_{k\ge1}
			\left[
			\frac{2^{k/2}}{(9d)^k}\binom{d+k}{k}
			\left(1+\frac{k}{d}\right)^{(n-2)/4}
			\right]^2\le C(n).
		\end{equation}
		To see this, use
		\((1+k/d)^{(n-2)/4}\le C(n)(3/2)^k\) and
		\[
		\sum_{k\ge1}\binom{d+k}{k}
		\left(\frac{\sqrt2}{6d}\right)^k
		=\left(1-\frac{\sqrt2}{6d}\right)^{-d-1}-1.
		\]
		The last expression is bounded uniformly for \(d\ge2\), and the
		square sum is bounded by the square of this positive series.
		Since
		\[
		P_{x,d}-P_{0,d}
		=\sum_{k\ge1}[P_{0,d+k}(x+\cdot)]_d,
		\]
		Cauchy--Schwarz, \(|x|^k\le|x|\),
		\eqref{eq:two-center-translation-sum}, and
		\eqref{eq:two-center-spectral-concentration} yield
		\[
		\|P_{x,d}-P_{0,d}\|
		\le C|x|
		\left(\sum_{m>d}(9d)^{2(m-d)}\|P_{0,m}\|^2\right)^{1/2}
		\le C\sqrt\varepsilon|x|\|P_{0,d}\|.
		\]
		Consequently,
		\(\|P_{x,d}\|\le(1+C\sqrt\varepsilon)\|P_{0,d}\|\).
		
		For the next lower coefficient one has the identity
		\begin{equation}\label{eq:two-center-degree-d-minus-one}
			P_{x,d-1}=P_{0,d-1}+x\cdot\nabla P_{0,d}
			+\sum_{k\ge1}[P_{0,d+k}(x+\cdot)]_{d-1}.
		\end{equation}
		The same derivative calculation, now with \(k+1\) derivatives, uses
		\(\binom{d+k}{k+1}\le d\binom{d+k}{k}\) and gives
		\[
		\sum_{k\ge1}
		\left[
		\frac{C(n)2^{(k+1)/2}}{(9d)^k}
		\binom{d+k}{k+1}
		\left(\frac{d+k}{d-1}\right)^{(n-2)/4}
		\right]^2\le C(n)d^2.
		\]
		Thus the final sum in
		\eqref{eq:two-center-degree-d-minus-one} has norm at most
		\(Cd\sqrt\varepsilon\|P_{0,d}\|\). Moreover,
		\eqref{eq:two-center-spectral-concentration} at the two centers gives
		\[
		\|P_{0,d-1}\|\le C\sqrt\varepsilon\|P_{0,d}\|,
		\qquad
		\|P_{x,d-1}\|\le C\sqrt\varepsilon\|P_{x,d}\|.
		\]
		It follows from \eqref{eq:two-center-degree-d-minus-one} that
		\[
		\|x\cdot\nabla P_{0,d}\|
		\le Cd\sqrt\varepsilon\|P_{0,d}\|.
		\]
		Finally,
		\(\|\nabla P_{0,d}\|^2=d(2d+n-2)\|P_{0,d}\|^2\), so the
		last estimate is the second inequality in
		\eqref{eq:Pxd-P0d-estimate}.
	\end{proof}
	
	\begin{lemma}\label{lem:almost-invariant-polynomial}
		Let \(P\in\mathcal P_d\setminus\{0\}\), \(d\ge1\), and let
		\(1\le k\le n-2\) and
		\(e_1,\ldots,e_k\) be orthonormal. If
		\[
		\|\partial_{e_i}P\|^2\le\delta\|\nabla P\|^2\quad(1\le i\le k),
		\qquad\delta\le c(n,k)d^{-4k},
		\]
		then there is \(P_1\in\mathcal P_d\), invariant along
		\(\operatorname{span}\{e_1,\ldots,e_k\}\),
		such that
		\[
		\|P-P_1\|\le C(n,k)d^k\sqrt\delta\|P\|.
		\]
	\end{lemma}
	\begin{proof}
		The projection argument follows
		\cite[Proposition~3.9 and Lemma~3.11]{NaberValtorta2017}; see also
		\cite[Lemma~4.15]{HuangJiang}. Let \(Q_j\) be the orthogonal
		projection of \(P\) onto the degree-\(d\) harmonic polynomials
		invariant along \(e_1,\ldots,e_j\), with \(Q_0=P\). The cited
		one-direction estimate, applied in the orthogonal complement of
		\(\operatorname{span}\{e_1,\ldots,e_{j-1}\}\), gives
		\[
		\|Q_{j-1}-Q_j\|\le C(n)\|\partial_{e_j}Q_{j-1}\|.
		\]
		Since
		\(\|\nabla F\|^2=d(2d+n-2)\|F\|^2\) for \(F\in\mathcal P_d\),
		\[
		\|P-Q_j\|
		\le(1+Cd)\|P-Q_{j-1}\|+Cd\sqrt\delta\|P\|.
		\]
		Iteration yields
		\(\|P-Q_k\|\le C(n,k)d^k\sqrt\delta\|P\|\).
		The smallness assumption makes this at most \(\|P\|/2\), so
		\(Q_k\ne0\). Take \(P_1=Q_k\).
	\end{proof}
	\begin{remark}\label{rem:normalized-almost-invariant-polynomial}
		If the bound in Lemma \ref{lem:almost-invariant-polynomial} is at most
		\(\|P\|/2\), then \(\|P_1\|\ge\|P\|/2\). Hence
		\eqref{eq:center-comparison-normalization-inequality}, applied with
		\(f=P_1\) and \(g=P\), and the lemma give
		\[
		\left\|\frac P{\|P\|}-\frac{P_1}{\|P_1\|}\right\|^2
		\le\frac{4\|P-P_1\|^2}{\|P\|^2}
		\le C(n,k)d^{2k}\delta.
		\]
	\end{remark}
	
	Fix a structural constant \(c_6>0\) sufficiently small for all the
	requirements below.
	Here \(\varepsilon_1\) is the constant fixed before Proposition
	\ref{prop:quantitative-harmonic-approximation-DiniLp}; \(r_4\) and
	\(r_5\) are defined in \eqref{eq:almost-monotonicity-radius} and
	\eqref{eq:regularity-radius-definition}; and
	\(A_{\mathcal D},d_{\mathcal D},\eta\) are defined in
	\eqref{eq:42-AD-definition}--\eqref{eq:prop335-Dini-smallness}.
	
	For \(0<\varepsilon\le\varepsilon_1/4\), let
	\(r_6(\varepsilon)\) be the supremum of
	the radii satisfying
	\begin{equation}\label{eq:42-rcs-definition}
		\begin{aligned}
			r&\le\frac{\min\{r_4(1/10),r_5,
				(10\sqrt\Lambda)^{-1}\}}
			{2^{12}(1+\lambda^{-1/2})(1+d_{\mathcal D})},\\
			A_{\mathcal D}^{C(n)}100^{d_{\mathcal D}+1}(1+d_{\mathcal D})^{C(n)}
			&\eta\bigl(2^{12}(1+\lambda^{-1/2})(1+d_{\mathcal D})r\bigr)
			\le c_6\varepsilon.
		\end{aligned}
	\end{equation}
	The two fixed exponents are chosen to dominate the powers below.
	This is a positive radius by the Dini condition. Its dependence on
	\(M\) and \(\mathcal D\) is understood as for \(r_1,\ldots,r_5\).
	
	\begin{proposition}\label{prop:nearby-center-comparison}
		Let \(x\in Z(u)\cap B_1\), \(0<r\le r_6(\varepsilon)\),
		and \(2\le d\le d_{\mathcal D}\). Suppose 
		\[
		|D^u(x,s)-d|\le\varepsilon
		\quad\left(r/100\le s\le
		2^{12}(1+\lambda^{-1/2})(1+d_{\mathcal D})r\right).
		\]
		Put \(R=20(1+\lambda^{-1/2})r\). We can find a single harmonic function
		\(h\) on \(B_{20d_{\mathcal D}}\) such that, for every
		\(y\in Z(u)\cap \overline B_r(x)\), with
		\(z=A_x^{-1}(y-x)/R\),
		\begin{align}
			|D^u(y,\sigma R)-D^h(z,\sigma)|&\le C\varepsilon,
			\label{eq:cone-splitting-long-doubling-comparison}\\
			\vint_{\partial B_1}|u_{y,\sigma R}-h_{z,\sigma}|^2
			&\le C\varepsilon^2,
			\label{eq:cone-splitting-long-blowup-comparison}
		\end{align}
		for \(1/10\le\sigma\le9d_{\mathcal D}\). Here
		\(h_{z,\sigma}(\xi)=h(z+\sigma\xi)/
		(\vint_{\partial B_1}h(z+\sigma\zeta)^2\,d\zeta)^{1/2}\).
	\end{proposition}
	
	\begin{proof}
		Choose the least integer \(J\) such that \(2^J\ge80d_{\mathcal D}\),
		and set
		\[
		R_0=2^JR,\qquad r_j=2^{-j}R_0
		\quad(1\le j\le J).
		\]
		Then \(r_J=R\) and
		\[
		80d_{\mathcal D}\le2^J<160d_{\mathcal D},\qquad
		R_0<2^{12}(1+\lambda^{-1/2})(1+d_{\mathcal D})r.
		\]
		The definition of \(r_6(\varepsilon)\) therefore makes \(R_0\)
		admissible in Proposition
		\ref{prop:quantitative-harmonic-approximation-DiniLp}.
		On \([r/100,R_0]\), the oscillation of \(D^u(x,\cdot)\) is at
		most \(2\varepsilon\le\varepsilon_1\). Let
		\[
		\gamma=\inf_{r/100\le t\le R_0}D^u(x,t).
		\]
		Apply Proposition \ref{prop:quantitative-harmonic-approximation-DiniLp}
		on this interval and denote its harmonic function by \(h_0\).
		Thus \(\gamma\ge d-\varepsilon>1\), and
		\begin{equation}\label{eq:cone-pointwise-error}
			|h_0(y)-u_{x,R_0}(y)|
			\le C\eta(R_0)
			\begin{cases}
				\left(\frac{r}{100R_0}\right)^{\gamma-1}|y|,& |y|\le\frac{r}{100R_0},\\
				|y|^\gamma,&\frac{r}{100R_0}\le |y|\le3/4.
			\end{cases}
		\end{equation}
		Define
		\[
		h_j(\xi)=\frac{h_0(2^{-j}\xi)}{\left(\vint_{\partial B_{2^{-j}}}h_0^2\right)^{1/2}}
		\quad(1\le j\le J).
		\]
		The denominator is nonzero by the following calculation. Since
		\(\left(\vint_{\partial B_{1}}u_{x,R_0}^2\right)=1\), the definition of the doubling index gives
		\begin{equation}\label{eq:cone-height-exact}
			\left(\vint_{\partial B_{2^{-j}}}u_{x,R_0}^2\right)^{1/2}
			=2^{-\sum_{l=1}^jD^u(x,2^{-l}R_0)}
			\ge 2^{-j(d+\varepsilon)}.
		\end{equation}
		Moreover \(\frac{r}{100R_0}\le2^{-J}\le2^{-j}\). Hence
		\eqref{eq:cone-pointwise-error}, both outside and inside \(B_{r/(100R_0)}\),
		gives
		\[
		\sup_{B_{2/3}}
		|h_0(2^{-j}\xi)-u_{x,R_0}(2^{-j}\xi)|
		\le C\eta(R_0)2^{-j\gamma}.
		\]
		Using \(\gamma\ge d-\varepsilon\) in
		\eqref{eq:cone-height-exact}, we find
		\begin{equation}\label{eq:cone-relative-outer-error}
			\frac{\sup_{B_{2/3}}
				|h_0(2^{-j}\cdot)-u_{x,R_0}(2^{-j}\cdot)|}
			{\left(\vint_{\partial B_{2^{-j}}}u_{x,R_0}^2\right)^{1/2}}
			\le C\eta(R_0)2^{2j\varepsilon}
			\le C(160d_{\mathcal D})^{2\varepsilon}\eta(R_0).
		\end{equation}
		The triangle inequality on \(\partial B_{2^{-j}}\) gives the
		same bound for
		\[
		\left|
		\frac{\left(\vint_{\partial B_{2^{-j}}}h_0^2\right)^{1/2}}
		{\left(\vint_{\partial B_{2^{-j}}}u_{x,R_0}^2\right)^{1/2}}-1
		\right|.
		\]
		After enlarging \(C\) in the preceding estimates, the choice of
		\(c_6\), the second condition in \eqref{eq:42-rcs-definition} implies
		\begin{equation}\label{eq:cone-master-error}
			C(160d_{\mathcal D})^{2\varepsilon}\eta(R_0)\le cA_{\mathcal D}^{-C(n)}
			100^{-d_{\mathcal D}-1}(1+d_{\mathcal D})^{-C(n)}
			\varepsilon.
		\end{equation}
		Combing \eqref{eq:cone-height-exact} and \eqref{eq:cone-master-error} shows $h_j(\xi)$ is well defined.
		
		Using the identity
		\[
		\begin{split}
			u_{x,r_j}-h_j={}&\frac{u_{x,R_0}(2^{-j}\cdot)-h_0(2^{-j}\cdot)}
			{\left(\vint_{\partial B_{2^{-j}}}u_{x,R_0}^2\right)^{1/2}}\\
			&+h_j\left(\frac{\left(\vint_{\partial B_{2^{-j}}}h_0^2\right)^{1/2}}
			{\left(\vint_{\partial B_{2^{-j}}}u_{x,R_0}^2\right)^{1/2}}-1\right)
		\end{split}
		\]
		Since \(\left(\vint_{\partial B_{1}}h_j^2\right)=1\), harmonic interior estimates give
		\[
		\|h_j\|_{L^\infty(B_{2/3})}\le C(n).
		\]
		We therefore obtain
		\begin{equation}\label{eq:common-harmonic-fixed-scale-error}
			\|u_{x,r_j}-h_j\|_{L^\infty(B_{2/3})}
			\le C(160d_{\mathcal D})^{2\varepsilon}\eta(R_0).
		\end{equation}
		
		Fix \(y\in Z(u)\cap\overline B_r(x)\), and write
		\[
		z_j=A_x^{-1}\frac{y-x}{r_j}.
		\]
		Then
		\begin{equation}\label{eq:nearby-center-common-variable}
			x+r_jA_x(z_j+sA_x^{-1}A_y\xi)=y+sr_jA_y\xi.
		\end{equation}
		Moreover,
		\[
		|z_j|\le
		\frac{\lambda^{-1/2}r}{20(1+\lambda^{-1/2})r}
		\le\frac1{20}.
		\]
		We record the matrix comparison:
		\[
		|A_x^{-1}A_y-I|\le C(\lambda)\omega(|y-x|)
		\le C\int_0^{2\sqrt\Lambda R_0}\frac{\omega(t)}t\,dt.
		\]
		Here we used monotonicity of \(\omega\) and
		\[
		\omega(r)\le(\log2)^{-1}\int_r^{2r}\frac{\omega(t)}t\,dt.
		\]
		Lemma \ref{lem:growth-estimate-gradient} gives
		\(\sup_{B_{2/3}}|\nabla u_{x,r_j}|\le CK A_{\mathcal D}^{1/2}\).
		By the choice of \(c_6\), \(|A_x^{-1}A_y-I|\le1/10\).
		For \(1/10\le s\le1/2\), the segment from
		\(z_j+s\xi\) to \(z_j+sA_x^{-1}A_y\xi\) lies in \(B_{2/3}\), since
		\[
		|z_j|+s(1+|A_x^{-1}A_y-I|)
		\le\frac1{20}+\frac{11}{20}=\frac35.
		\]
		The fundamental theorem of calculus and the definition of \(\eta\)
		therefore give
		\begin{equation}\label{eq:nearby-center-metric-error}
			\sup_{\partial B_1}
			|u_{x,r_j}(z_j+sA_x^{-1}A_y\xi)-u_{x,r_j}(z_j+s\xi)|
			\le C\eta(R_0)\le C(160d_{\mathcal D})^{2\varepsilon}\eta(R_0).
		\end{equation}
		
		It remains to obtain a lower bound for the nearby normalization.
		Three applications of the pinching assumption at \(x\) give
		\[
		\left(\vint_{\partial B_{1/8}}u_{x,r_j}^2\right)^{1/2}
		=2^{-D^u(x,r_j/8)-D^u(x,r_j/4)-D^u(x,r_j/2)}
		\ge8^{-d-\varepsilon}\ge8^{-d_{\mathcal D}-1}.
		\]
		Inequalities \eqref{eq:common-harmonic-fixed-scale-error} and
		\eqref{eq:cone-master-error} therefore imply
		\[
		\left(\vint_{\partial B_{1/8}}h_j^2\right)^{1/2}
		\ge\frac12\,8^{-d_{\mathcal D}-1}.
		\]
		Since harmonic spherical
		means are increasing,
		\[
		\left(\vint_{B_{1/4}}h_j^2\right)\ge c\,8^{-2d_{\mathcal D}-2}.
		\]
		The inclusions
		\[
		B_{1/4}\subset B_{3/10}(z_j)
		\subset B_{3/5}(z_j)\subset B_{2/3}
		\]
		and \(\sup_{B_{2/3}}|h_j|\le C\) show that
		\[
		\left(\vint_{B_{3/10}(z_j)}h_j^2\right)\ge c\,8^{-2d_{\mathcal D}-2},
		\qquad \left(\vint_{B_{3/5}(z_j)}h_j^2\right)\le C.
		\]
		Logarithmic convexity of harmonic ball means at radii
		\[
		\frac1{10},\qquad\frac3{10},\qquad\frac35
		\]
		gives
		\[
		\left(\vint_{B_{3/10}(z_j)}h_j^2\right)
		\le \left(\vint_{B_{1/10}(z_j)}h_j^2\right)^{\log2/\log6}
		\left(\vint_{B_{3/5}(z_j)}h_j^2\right)^{1-\log2/\log6}.
		\]
		Solving for the first factor and using the relation between
		\(d_{\mathcal D}\) and \(A_{\mathcal D}\), we obtain
		\[
		\left(\vint_{B_{1/10}(z_j)}h_j^2\right)\ge cA_{\mathcal D}^{-C(n)}.
		\]
		Hence
		\begin{equation}\label{eq:nearby-center-relative-height}
			\left(\vint_{\partial B_1}
			h_j(z_j+s\xi)^2\,d\xi\right)^{1/2}
			\ge c A_{\mathcal D}^{-C(n)}
			\qquad(1/10\le s\le1/2).
		\end{equation}
		
		For fixed \(j\) and \(s\in[1/10,1/2]\), inequalities \eqref{eq:common-harmonic-fixed-scale-error},
		\eqref{eq:nearby-center-metric-error}, and
		\eqref{eq:nearby-center-relative-height} give
		\[
		\frac{
			\left(\vint_{\partial B_1}
			|u_{x,r_j}(z_j+sA_x^{-1}A_y\xi)-h_j(z_j+s\xi)|^2\,d\xi\right)^{1/2}}
		{\left(\vint_{\partial B_1}h_j(z_j+s\xi)^2\,d\xi\right)^{1/2}}
		\le CA_{\mathcal D}^{C(n)}(160d_{\mathcal D})^{2\varepsilon}\eta(R_0)
		\le\frac{\varepsilon}{100}.
		\]
		By \eqref{eq:nearby-center-common-variable},
		\[
		u_{y,sr_j}(\xi)=
		\frac{u_{x,r_j}(z_j+sA_x^{-1}A_y\xi)}
		{\left(\vint_{\partial B_1}
			u_{x,r_j}(z_j+sA_x^{-1}A_y\zeta)^2\,d\zeta\right)^{1/2}},
		\qquad
		(h_j)_{z_j,s}(\xi)=
		\frac{h_j(z_j+s\xi)}
		{\left(\vint_{\partial B_1}
			h_j(z_j+s\zeta)^2\,d\zeta\right)^{1/2}}.
		\]
		The normalization inequality
		\eqref{eq:center-comparison-normalization-inequality} yields
		\begin{equation}\label{eq:cone-stage-blowup}
			\vint_{\partial B_1}
			|u_{y,sr_j}-(h_j)_{z_j,s}|^2\le C\varepsilon^2.
		\end{equation}
		For \(1/10\le s\le1/4\), the same relative estimate at \(s\)
		and \(2s\) gives
		\begin{equation}\label{eq:cone-stage-doubling}
			|D^u(y,sr_j)-D^{h_j}(z_j,s)|
			\le2\log_4\frac{1+\varepsilon/100}{1-\varepsilon/100}
			\le C\varepsilon.
		\end{equation}
		
		Finally take \(h=h_J\). Rescaling the preceding comparisons at
		\(j=1,\ldots,J\) produces no additional error. Indeed,
		\(r_j=2^{J-j}R\) and \(z_j=2^{j-J}z_J\), so the definitions give
		the identities
		\[
		D^{h_j}(z_j,s)=D^h(z_J,2^{J-j}s),\qquad
		(h_j)_{z_j,s}=h_{z_J,2^{J-j}s},\qquad
		u_{y,sr_j}=u_{y,2^{J-j}sR}.
		\]
		The intervals overlap and satisfy
		\[
		\bigcup_{j=1}^J[2^{J-j}/10,2^{J-j}/4]
		=[1/10,2^J/8]\supset[1/10,9d_{\mathcal D}].
		\]
		Equations \eqref{eq:cone-stage-blowup} and
		\eqref{eq:cone-stage-doubling} therefore prove the two stated
		comparisons. Since \(h_J\) is harmonic in \(B_{2^J}\) and
		\(2^J\ge80d_{\mathcal D}\), it is harmonic on
		\(B_{20d_{\mathcal D}}\).
	\end{proof}
	
	The next lemma records the explicit polynomial error in the
	cone-splitting argument of \cite[Proposition~4.18]{HuangJiang}.
	The estimate below is a quantitative consequence of that argument,
	not the stated error bound of the cited proposition.
	\begin{lemma}\label{lem:quantitative-cone-splitting}
		Let \(h\) be harmonic in \(B_{20d}\), where \(d\ge2\), and let
		\(1\le k\le n-2\) and \(0<\tau\le1\). Suppose that
		\(v_1,\ldots,v_k\in B_1\) satisfy
		\[
		\operatorname{dist}\bigl(v_i,
		\operatorname{span}(v_1,\ldots,v_{i-1})\bigr)\ge\tau
		\quad(1\le i\le k),
		\]
		and
		\[
		|D^h(a,s)-d|\le\varepsilon
		\quad(a=0,v_1,\ldots,v_k,\quad 1/10\le s\le9d),
		\qquad
		\varepsilon\le c(n,k)\tau^{2k}d^{-4k}.
		\]
		Then \(P=P_{0,d}\ne0\), and there is a nonzero
		\(Q\in\mathcal P_d\), invariant along
		\(\operatorname{span}(v_1,\ldots,v_k)\), such that
		\begin{equation}\label{eq:cone-invariant-polynomial-estimate}
			\|P-Q\|\le C(n,k)d^k\tau^{-k}\sqrt\varepsilon\|P\|,
			\qquad
			\left\|\frac P{\|P\|}-\frac Q{\|Q\|}\right\|^2
			\le C(n,k)d^{2k}\tau^{-2k}\varepsilon.
		\end{equation}
	\end{lemma}
	\begin{proof}
		By Lemma \ref{lem:two-center-coefficient-estimate},
		\(\|v_i\cdot\nabla P\|
		\le C\sqrt\varepsilon\|\nabla P\|\).
		Put \(w_i=v_i/|v_i|\). Since \(|v_i|\ge\tau\),
		\(\|\partial_{w_i}P\|
		\le C\tau^{-1}\sqrt\varepsilon\|\nabla P\|\).
		Gram--Schmidt gives orthonormal vectors
		\(e_i=\sum_{l\le i}b_{il}w_l\), with
		\(\sum_{l\le i}|b_{il}|\le C(n,k)\tau^{1-k}\). Hence
		\[
		\|\partial_{e_i}P\|^2
		\le C(n,k)\tau^{-2k}\varepsilon\|\nabla P\|^2.
		\]
		The smallness of \(\varepsilon\) allows us to apply Lemma
		\ref{lem:almost-invariant-polynomial} and Remark
		\ref{rem:normalized-almost-invariant-polynomial} with
		\(\delta=C(n,k)\tau^{-2k}\varepsilon\), proving
		\eqref{eq:cone-invariant-polynomial-estimate}.
	\end{proof}
	
	For an affine set \(L\subset\mathbb R^n\), write
	\[
	B_s(L)=\{z\in\mathbb R^n:\operatorname{dist}(z,L)<s\}.
	\]
	If \(V\) is a linear subspace, put
	\[
	B_s(V)=B_1\cap\{z\in\mathbb R^n:\operatorname{dist}(z,V)<s\}.
	\]
	
	Fix a structural constant \(c_7>0\) sufficiently small for all the
	requirements below.
	Here \(\varepsilon_3\) is the constant from Proposition
	\ref{prop:low-doubling-regularity}, \(C_0\) is the constant fixed after
	Lemma \ref{lem:doubling-over-balls-bound}, and \(d_{\mathcal D}\) is
	defined in \eqref{eq:42-AD-definition}.
	
	For \(0<\tau<c(n,\lambda,\Lambda)\), put
	\begin{equation}\label{eq:confinement-epsilon-choice}
		\varepsilon_4(\tau)
		=c_7\varepsilon_3^2
		(\tau/3)^{2C_0(1+\mathcal D)}
		\le c_7\varepsilon_3^2
		(\tau/3)^{2d_{\mathcal D}-2}.
	\end{equation}
	
	\begin{lemma}\label{lem:singular-confinement}
		Let \(0<\tau<c(n,\lambda,\Lambda)\),
		\(x\in Z(u)\cap B_1\), and
		\(r\le r_6(\varepsilon)\). Suppose that \(u\) is
		\((n-2,\varepsilon,3r,x)\)-symmetric with respect to an
		\((n-2)\)-plane \(V\), in the sense of Definition
		\ref{def:quantitative-symmetry}, with corresponding polynomial \(P\)
		of degree at most \(d_{\mathcal D}\). If
		\(\varepsilon\le\varepsilon_4(\tau)\), then
		\begin{equation}\label{eq:confinement-conclusion-scaled}
			S(u)\cap(x+rA_xB_1)\subset x+rA_xB_\tau(V).
		\end{equation}
		In particular,
		\begin{equation}\label{eq:confinement-conclusion-euclidean}
			S(u)\cap B_{\sqrt\lambda r}(x)
			\subset B_{\sqrt\Lambda\tau r}(x+rA_xV).
		\end{equation}
	\end{lemma}
	
	\begin{proof}
		Since \(3r\le r_4(1/10)\),
		Theorem \ref{thm:almost-monotonicity} gives
		\(D^u(x,s)\ge9/10\) for \(s\le3r\). Lemma
		\ref{lem:growth-estimate-gradient}, followed by the scaled interior
		estimate, yields
		\begin{equation}\label{eq:confinement-C1-bound}
			\vint_{B_2}u_{x,3r}^2\le A_{\mathcal D},\qquad
			\sup_{B_{5/4}}(|u_{x,3r}|+|\nabla u_{x,3r}|)
			\le CK A_{\mathcal D}^{1/2}.
		\end{equation}
		In the coordinates of Lemma \ref{lem:scaled-coefficients},
		\[
		\Delta u_{x,3r}=\partial_iF^i-\widetilde V u_{x,3r},\qquad
		F^i=(\delta^{ij}-\widetilde a^{ij})\partial_j u_{x,3r}.
		\]
		Equations \eqref{eq:scaled-a-modulus},
		\eqref{eq:scaled-V-Lp}, and
		\eqref{eq:confinement-C1-bound} give
		\begin{equation}\label{eq:confinement-source-bounds}
			\|F\|_{L^\infty(B_{5/4})}
			\le CK A_{\mathcal D}^{1/2}\omega(4\sqrt\Lambda r),
			\qquad
			\|\widetilde V\|_{L^p(B_{5/4})}
			\le CMr^{2-n/p}.
		\end{equation}
		
		Choose \(\chi\in C_c^\infty(B_{5/4})\), with
		\(\chi=1\) on \(B_{9/8}\), and let \(\Gamma\) be the fundamental
		solution of the Laplacian. Define
		\[
		\phi(y)=\int\Gamma(y-z)\chi(z)
		\bigl(\partial_iF^i(z)-\widetilde V(z)u_{x,3r}(z)\bigr)\,dz.
		\]
		Integration by parts in the first term gives
		\begin{align*}
			\phi(y)
			={}&-\int\partial_{z_i}\bigl(\Gamma(y-z)\chi(z)\bigr)
			F^i(z)\,dz\\
			&-\int\Gamma(y-z)\chi(z)\widetilde V(z)u_{x,3r}(z)\,dz.
		\end{align*}
		The gradient of \(\Gamma\) is locally integrable. Since \(p>n\),
		\[
		\sup_{y\in B_{9/8}}
		\|\Gamma(y-\cdot)\|_{L^{p/(p-1)}(B_{5/4})}\le C(n,p).
		\]
		Thus \eqref{eq:confinement-C1-bound} and
		\eqref{eq:confinement-source-bounds} imply
		\[
		\|\phi\|_{L^\infty(B_{9/8})}
		\le CK A_{\mathcal D}^{1/2}
		\bigl[\omega(4\sqrt\Lambda r)+Mr^{2-n/p}\bigr]
		\le C\eta(3r).
		\]
		In the last inequality we used
		\[
		\omega(4\sqrt\Lambda r)
		\le\frac1{\log(3/2)}
		\int_{4\sqrt\Lambda r}^{6\sqrt\Lambda r}
		\frac{\omega(t)}t\,dt.
		\]
		The function \(h=u_{x,3r}-\phi\) is harmonic in \(B_{9/8}\). By
		\eqref{eq:42-rcs-definition}, by the choice of \(c_6\),
		\begin{equation}\label{eq:confinement-harmonic-error}
			\|u_{x,3r}-h\|_{L^\infty(B_{9/8})}\le\varepsilon.
		\end{equation}
		
		The assumed approximation on \(\partial B_1\) and
		\eqref{eq:confinement-harmonic-error} give
		\begin{equation}\label{eq:confinement-h-P-boundary}
			\left(\vint_{\partial B_1}|h-P|^2\right)^{1/2}
			\le\sqrt\varepsilon+\varepsilon\le2\sqrt\varepsilon.
		\end{equation}
		Degree zero is impossible. Indeed, in that case \(P\) is a constant
		of absolute value one, while the mean-value property and
		\eqref{eq:confinement-h-P-boundary} give
		\(|h(0)-P|\le2\sqrt\varepsilon\), and
		\eqref{eq:confinement-harmonic-error} gives
		\(|h(0)|\le\varepsilon\). Hence
		\[
		1\le m:=\deg P\le d_{\mathcal D}.
		\]
		Harmonic interior estimates applied to \(h-P\), and then to \(h\),
		yield
		\begin{equation}\label{eq:confinement-C2-bounds}
			\|h-P\|_{C^2(B_{1/2})}\le C\sqrt\varepsilon,\qquad
			\|D^2h\|_{L^\infty(B_{1/2})}\le C.
		\end{equation}
		
		After a rotation, \(V=\{z_1=z_2=0\}\). Since \(P\) depends only
		on the two perpendicular variables, it has the form
		\[
		P(\varrho,\theta)
		=\varrho^m(a\cos m\theta+b\sin m\theta).
		\]
		Its normalization gives
		\[
		1=\vint_{\partial B_1}P^2
		=\frac{a^2+b^2}{2}
		\vint_{\partial B_1}\varrho^{2m}
		\le\frac{a^2+b^2}{2}.
		\]
		Therefore \(a^2+b^2\ge2\), and
		\[
		|\nabla P|=m\sqrt{a^2+b^2}\,\varrho^{m-1}.
		\]
		If \(z\in B_{1/3}\) and
		\(\operatorname{dist}(z,V)\ge\tau/3\), then
		\begin{equation}\label{eq:gradient-lower-P-away-V}
			|\nabla P(z)|\ge\sqrt2\left(\frac{\tau}{3}\right)^{d_{\mathcal D}-1}.
		\end{equation}
		Because
		\[
		\varepsilon\le\varepsilon_4(\tau)
		\le c_7\varepsilon_3^2\left(\frac{\tau}{3}\right)^{2d_{\mathcal D}-2},
		\]
		by the choice of \(c_7\) in \eqref{eq:confinement-epsilon-choice}
		and using \eqref{eq:confinement-C2-bounds} gives
		\begin{equation}\label{eq:confinement-h-gradient-lower}
			|\nabla h(z)|\ge\frac1{\sqrt2}\left(\frac{\tau}{3}\right)^{d_{\mathcal D}-1}.
		\end{equation}
		
		Suppose that a singular point violates
		\eqref{eq:confinement-conclusion-scaled}. Write
		\[
		y=x+rA_x\xi,\qquad |\xi|<1,\qquad
		\operatorname{dist}(\xi,V)\ge\tau.
		\]
		Then \(\xi/3\in B_{1/3}\),
		\(\operatorname{dist}(\xi/3,V)\ge\tau/3\), and \(u_{x,3r}(\xi/3)=0\).
		The symmetry and ellipticity bounds for \(A_x,A_y\) give
		\begin{equation}\label{eq:confinement-matrix-bounds}
			\|A_x^{-1}A_y\|\le\sqrt{\frac\Lambda\lambda},\qquad
			|A_yA_x^{-1}v|\ge\sqrt{\frac\lambda\Lambda}|v|
			\quad(v\in\mathbb R^n).
		\end{equation}
		Choose a sufficiently small constant \(c=c(n,\lambda,\Lambda)>0\) and put
		\[
		\rho=c\varepsilon_3\left(\frac{\tau}{3}\right)^{d_{\mathcal D}-1}.
		\]
		Equations \eqref{eq:confinement-C2-bounds},
		\eqref{eq:confinement-h-gradient-lower}, and
		\eqref{eq:confinement-matrix-bounds} then give
		\begin{align}
			2\rho\|A_x^{-1}A_y\|&\le\frac1{12},\label{eq:confinement-rho-domain}\\
			2\rho\|A_x^{-1}A_y\|^2\|D^2h\|_{L^\infty(B_{1/2})}
			&\le\frac{\varepsilon_3}{200\sqrt n}
			|A_yA_x^{-1}\nabla h(\xi/3)|,\label{eq:confinement-rho-remainder}\\
			\rho|A_yA_x^{-1}\nabla h(\xi/3)|
			&\ge c\varepsilon_3\left(\frac{\tau}{3}\right)^{2d_{\mathcal D}-2}.
			\label{eq:confinement-linear-lower}
		\end{align}
		Also, \(u_{x,3r}(\xi/3)=0\) and
		\eqref{eq:confinement-harmonic-error} imply
		\(|h(\xi/3)|\le\varepsilon\). By the choice of \(c_7\), the bound
		\(\varepsilon\le c_7\varepsilon_3^2\left(\frac{\tau}{3}\right)^{2d_{\mathcal D}-2}\) and
		\eqref{eq:confinement-linear-lower} give
		\begin{equation}\label{eq:confinement-constant-error}
			2\|u_{x,3r}-h\|_{L^\infty(B_{9/8})}+|h(\xi/3)|
			\le\frac{\varepsilon_3}{200\sqrt n}
			\rho|A_yA_x^{-1}\nabla h(\xi/3)|.
		\end{equation}
		
		For \(t=1,2\), Taylor's formula is valid on
		\(\xi/3+t\rho A_x^{-1}A_y\partial B_1\subset B_{1/2}\), by
		\eqref{eq:confinement-rho-domain}, and gives
		\[
		h(\xi/3+t\rho A_x^{-1}A_y\zeta)
		=h(\xi/3)+t\rho A_yA_x^{-1}\nabla h(\xi/3)\cdot\zeta+R_t(\zeta),
		\]
		where
		\begin{equation}\label{eq:confinement-Taylor-remainder}
			|R_t(\zeta)|
			\le\frac{t^2\rho^2\|A_x^{-1}A_y\|^2}{2}
			\|D^2h\|_{L^\infty(B_{1/2})}.
		\end{equation}
		For \(t=1,2\), one has
		\[
		\left(\vint_{\partial B_1}
		|t\rho A_yA_x^{-1}\nabla h(\xi/3)\cdot\zeta|^2\,d\zeta\right)^{1/2}
		=\frac{t\rho|A_yA_x^{-1}\nabla h(\xi/3)|}{\sqrt n},
		\]
		the reverse triangle inequality, together with
		\eqref{eq:confinement-rho-remainder}--
		\eqref{eq:confinement-Taylor-remainder}, yields
		\begin{equation}\label{eq:confinement-spherical-average-estimate}
			\left|
			\left(\vint_{\partial B_1}
			u_{x,3r}(\xi/3+t\rho A_x^{-1}A_y\zeta)^2\,d\zeta\right)^{1/2}
			-\frac{t\rho|A_yA_x^{-1}\nabla h(\xi/3)|}{\sqrt n}
			\right|
			\le\frac{\varepsilon_3}{50}
			\frac{t\rho|A_yA_x^{-1}\nabla h(\xi/3)|}{\sqrt n}.
		\end{equation}
		
		The change of center is:
		\[
		x+3rA_x(\xi/3+t\rho A_x^{-1}A_y\zeta)=y+3tr\rho A_y\zeta.
		\]
		The normalization in \(u_{x,3r}\) cancels from the
		doubling quotient. Thus
		\[
		D^u(y,3r\rho)=\log_4
		\frac{\vint_{\partial B_1}
			u_{x,3r}(\xi/3+2\rho A_x^{-1}A_y\zeta)^2\,d\zeta}
		{\vint_{\partial B_1}
			u_{x,3r}(\xi/3+\rho A_x^{-1}A_y\zeta)^2\,d\zeta}.
		\]
		Using \eqref{eq:confinement-spherical-average-estimate}, we obtain
		\[
		D^u(y,3r\rho)
		\le1+2\log_4
		\frac{1+\varepsilon_3/50}{1-\varepsilon_3/50}
		\le1+\frac{\varepsilon_3}{10}.
		\]
		Moreover \(3r\rho\le3r\le r_5\), where \(r_5\) is defined in
		\eqref{eq:regularity-radius-definition}. All hypotheses of
		Proposition \ref{prop:low-doubling-regularity} hold at this smaller
		scale with error \(\varepsilon_3/10\). That proposition gives
		\(\nabla u(y)\ne0\), contradicting \(y\in S(u)\). This proves
		\eqref{eq:confinement-conclusion-scaled}.
		
		Finally, if
		\(y\in S(u)\cap B_{\sqrt\lambda r}(x)\), then
		\[
		\xi=A_x^{-1}\frac{y-x}{r}
		\quad\text{satisfies}\quad |\xi|<1.
		\]
		The scaled inclusion gives
		\(\operatorname{dist}(\xi,V)<\tau\). If
		\(\xi=v+e\), with \(v\in V\), \(e\perp V\), then
		\[
		\operatorname{dist}(y,x+rA_xV)
		\le r|A_xe|
		\le\sqrt\Lambda r|e|
		<\sqrt\Lambda\tau r.
		\]
		This is \eqref{eq:confinement-conclusion-euclidean}.
	\end{proof}

	\section{Proof of Theorem \ref{thm:main-volume-estimate}}
	\label{sec:covering}\label{sec:singular-covering-volume}
	
	We retain \(K\) from \eqref{eq:K-definition};
	\(C_0\) from Lemma \ref{lem:doubling-over-balls-bound}, with the
	enlargements fixed in the harmonic-approximation subsection; and
	\(A_{\mathcal D},d_{\mathcal D},\eta\) from
	\eqref{eq:42-AD-definition}--\eqref{eq:prop335-Dini-smallness};
	\(\varepsilon_1\) from the choice preceding Proposition
	\ref{prop:quantitative-harmonic-approximation-DiniLp}; \(c_3\) from
	Theorem \ref{thm:uniform-symmetry-pinched}; and
	\(q_0,\varepsilon_3,c_5\) from Proposition
	\ref{prop:low-doubling-regularity}. The radii \(r_2,r_4,r_6\) are fixed
	in Lemma \ref{lem:spherical-doubling-bound-singular},
	\eqref{eq:almost-monotonicity-radius}, and
	\eqref{eq:42-rcs-definition}, respectively. All constants denoted by
	\(c\) and \(C\) in this section depend only on
	\(n,p,\lambda,\Lambda,\omega\), unless a smaller dependence is indicated.
	The constants may decrease or increase, respectively, at subsequent
	occurrences. We continue to use the spherical inner product and norm
	from \eqref{eq:standing-spherical-normalization}.
	
	Put
	\begin{equation}\label{eq:5-Cpin}
		C_3=2^{12}(1+\lambda^{-1/2})(1+d_{\mathcal D}).
	\end{equation}
	Choose \(0<\tau<c(n,\lambda,\Lambda)\) in the range required by Lemma
	\ref{lem:singular-confinement}, and put
	\begin{equation}\label{eq:5-tau-geo}
		\tau_1
		=\frac{\tau}{10^8(1+\lambda^{-1/2})(1+\sqrt{\Lambda/\lambda})}.
	\end{equation}
	We fix
	\begin{equation}\label{eq:5-epsilon-choice}
		\varepsilon
		=c(1+d_{\mathcal D})^{-8n}
		\left(\frac{\tau_1}{3}\right)^{
			4C_0(n,\lambda,\Lambda)(1+\mathcal D)}.
	\end{equation}
	Here \(c\) is chosen so that Lemma \ref{lem:integer-pinching-drop} and
	Proposition \ref{prop:low-doubling-regularity} apply,
	\(C\varepsilon^{1/2}\le\varepsilon_4(\tau_1)\), with
	\(\varepsilon_4\) defined in \eqref{eq:confinement-epsilon-choice}, and
	the polynomial estimates below hold.
	In particular,
	\begin{equation}\label{eq:5-epsilon-simple-bound}
		\varepsilon^{-1}\le C\exp\bigl(C(1+\mathcal D)\bigr).
	\end{equation}
	The polynomial factor in \eqref{eq:5-epsilon-choice} also ensures
	\(\varepsilon\le c\tau^8(1+d_{\mathcal D})^{-8n}\).
	
	We first fix the radius required for the degree drop. Let \(r_7\)
	be the supremum of the radii
	\(0<R\le\min\{r_4(\varepsilon/100),(16\sqrt\Lambda)^{-1}\}\)
	satisfying
	\begin{equation}\label{eq:5-rdrop}
		\omega(8\sqrt\Lambda R)+MR^{2-n/p}
		\le c\left(\frac{\varepsilon}{100}\right)^{
			5(n+2)C_0(n,\lambda,\Lambda)(1+\mathcal D)}
		\bigl(K4^{3C_0(1+\mathcal D)}\bigr)^{-(n+2)/2}.
	\end{equation}
	At these radii the hypotheses of the definite-drop conclusion in
	Lemma \ref{lem:integer-pinching-drop} hold. The Dini condition and
	\(2-n/p>0\) ensure \(r_7>0\).
	
	We choose \(r_8>0\) so small that the spherical doubling
	bound, the definite-drop conclusion of Lemma
	\ref{lem:integer-pinching-drop}, and the regularity and comparison
	results of the preceding section apply at every scale at most
	\(4\cdot10^8 C_3r_8\). More precisely, we require
	\begin{equation}\label{eq:5-rmaster}
		\begin{aligned}
			&4\cdot10^8 C_3r_8
			\le\min\{r_2,r_6(\varepsilon),r_7,1\},\\
			&C A_{\mathcal D}^{C}\,
			\omega(4\cdot10^8 C_3r_8)
			\le \frac{\tau\varepsilon^{1/2}}{10^8}.
		\end{aligned}
	\end{equation}
	Here \(r_2\) is the radius from Lemma
	\ref{lem:spherical-doubling-bound-singular}, \(r_6(\varepsilon)\) is
	defined in \eqref{eq:42-rcs-definition}, and \(r_7\) is defined in
	\eqref{eq:5-rdrop}.
	The fixed power of \(A_{\mathcal D}\) is chosen to dominate the powers in
	the estimates below. We use the same convention in the definition
	of \(r_6\). These choices are possible by the Dini condition.
	At the end of the proof we give an explicit admissible value of
	\(r_8\).
	
	Finally, set
	\begin{equation}\label{eq:5-theta}
		\theta=\frac{\varepsilon}{10^4C_3},\qquad
		\mathcal A=C(n,\lambda,\Lambda)(10^6C_3)^n\theta^{-2},
	\end{equation}
	where the last constant is sufficiently large. Then
	\begin{equation}\label{eq:5-A5-exponential-D}
		0<\theta<\frac17,\qquad
		4C_3\theta=\frac{\varepsilon}{2500},
		\qquad 1\le\mathcal A\le \exp\bigl(C(1+\mathcal D)\bigr).
	\end{equation}
	
	We need to choose the approximating polynomials consistently as the scale
	varies, which contributes the Corollary \ref{cor:5-canonical-coherence} below. Before that, we need the following lemma to quantify the error between two adjacent scales through the Dini modulus and $M:=\|V\|_{L^p}$. Basically, this is the essential difference compared to the H\"older case in \cite{HuangJiang}.
	
	\begin{lemma}\label{lem:5-canonical-polynomial-center}
		Let \(2\le d\le d_{\mathcal D}\), \(0<\sigma\le4r_8\),
		\(x_i\in S(u)\), and
		\(0<\rho_i\le \sigma/(10^5 C_3)\).
		Assume
		\begin{equation}\label{eq:5-canonical-center-pinching}
			d-\varepsilon<D^u(x_i,s)\le d+\varepsilon
			\qquad(2\rho_i\le s\le C_3\sigma).
		\end{equation}
		For every
		\begin{equation}\label{eq:5-canonical-scale-range}
			200\rho_i\le R\le400(1+\lambda^{-1/2})\sigma
		\end{equation}
		there is a normalized \(P_{i,R}\in\mathcal P_d\) such that
		\begin{equation}\label{eq:5-canonical-scale-approximation}
			\|u_{x_i,R/8}-P_{i,R}\|^2\le C\varepsilon.
		\end{equation}
		When \(R\) and \(R/2\) both satisfy \eqref{eq:5-canonical-scale-range}, we have
		\begin{equation}\label{eq:5-canonical-adjacent-polynomials}
			\|P_{i,R}-P_{i,R/2}\|
			\le C K A_{\mathcal D}^{C}2^{d_{\mathcal D}+1}
			\bigl(\omega(2\sqrt\Lambda R)+MR^{2-n/p}\bigr).
		\end{equation}
	\end{lemma}
	
	\begin{proof}
		Let \(H_{i,R}\) be harmonic in \(B_1\) with boundary value
		\(u_{x_i,R}\).  Write
		\(H_{i,R}=\sum_{m\ge0}H_{i,R,m}\), where
		\(H_{i,R,m}\in\mathcal P_m\).
		The scale assumptions give
		\[
		\frac R{100}\ge2\rho_i,
		\qquad
		R\le400(1+\lambda^{-1/2})\sigma<C_3\sigma,
		\]
		so the pinching hypothesis is available on \([R/100,R]\). Applying
		Lemma \ref{lem:growth-estimate-gradient} with
		\(t_1=R\), \(t_2=R/100\), and
		\(\gamma=d-\varepsilon\), and using
		\(\vint_{\partial B_1}u_{x_i,R}^2=1\), gives
		\[
		\|u_{x_i,R}\|_{L^\infty(B_1)}
		+\|\nabla u_{x_i,R}\|_{L^\infty(B_1)}
		\le C K A_{\mathcal D}^{1/2}.
		\]
		Consequently,
		\begin{align*}
			\|u_{x_i,R}\|_{H^1(B_1)}
			+\|u_{x_i,R}\|_{L^{2p/(p-2)}(B_1)}
			&\le C\left(
			\|u_{x_i,R}\|_{L^\infty(B_1)}
			+\|\nabla u_{x_i,R}\|_{L^\infty(B_1)}\right)\\
			&\le C K A_{\mathcal D}^{C}.
		\end{align*}
		Set \(e=u_{x_i,R}-H_{i,R}\in H^1_0(B_1)\). Since the harmonic
		replacement has zero weak Dirichlet pairing with \(e\), the scaled equation
		gives
		\[
		\int_{B_1}|\nabla e|^2
		=\int_{B_1}(I-\widetilde a)\nabla u_{x_i,R}\cdot\nabla e
		+\int_{B_1}\widetilde V u_{x_i,R}e.
		\]
		Consequently, H\"older's inequality and Poincar\'e's inequality imply
		\begin{align*}
			\|\nabla e\|_2^2
			&\le\|I-\widetilde a\|_\infty
			\|\nabla u_{x_i,R}\|_2\|\nabla e\|_2
			+\|\widetilde V\|_p
			\|u_{x_i,R}\|_{2p/(p-2)}\|e\|_2\\
			&\le C K A_{\mathcal D}^{C}
			\bigl(\omega(2\sqrt\Lambda R)+MR^{2-n/p}\bigr)
			\|\nabla e\|_2.
		\end{align*}
		All the norms in this display are over \(B_1\). Dividing by
		\(\|\nabla e\|_2\), unless it is zero, gives
		\begin{equation}\label{eq:5-canonical-replacement-energy}
			\|e\|_{H^1(B_1)}
			\le C K A_{\mathcal D}^{C}
			\bigl(\omega(2\sqrt\Lambda R)+MR^{2-n/p}\bigr).
		\end{equation}
		The trace theorem gives the same bound for
		\((\vint_{\partial B_t}e^2)^{1/2}\) at the four fixed radii
		\[
		t\in\{1/100,1/50,1/8,1/2\}.
		\]
		The lower spherical bound under pinching, \eqref{eq:lower-bound-normalized-spherical-average},
		gives
		\begin{equation}\label{eq:5-canonical-sphere-lower}
			\left(\vint_{\partial B_t}u_{x_i,R}^2\right)^{1/2}
			\ge c A_{\mathcal D}^{-1/2}
			t^{d+\varepsilon+\varepsilon_1}
			\qquad(1/100\le t\le1/2).
		\end{equation}
		For the four fixed radii, the last lower bound is at least
		\(cA_{\mathcal D}^{-1/2}100^{-(d_{\mathcal D}+1)}\).
		Combining \eqref{eq:5-canonical-replacement-energy} with this lower
		bound yields, at each of the four radii,
		\begin{align*}
			\frac{(\vint_{\partial B_t}e^2)^{1/2}}
			{(\vint_{\partial B_t}u_{x_i,R}^2)^{1/2}}
			&\le C K A_{\mathcal D}^{C+1/2}100^{d_{\mathcal D}+1}
			\bigl(\omega(2\sqrt\Lambda R)+MR^{2-n/p}\bigr)\\
			&\le\frac{\varepsilon}{100},
		\end{align*}
		where the final inequality is included in the smallness conditions in
		\eqref{eq:5-rmaster} and in the definition of \(r_6\).
		The triangle inequality gives
		\begin{equation}\label{eq:5-canonical-relative-spheres}
			\left|\frac{(\vint_{\partial B_t}H_{i,R}^2)^{1/2}}
			{(\vint_{\partial B_t}u_{x_i,R}^2)^{1/2}}-1\right|
			\le\frac{\varepsilon}{100}
			\quad\left(t\in\left\{\frac1{100},\frac1{50},\frac18,\frac12\right\}\right).
		\end{equation}
		At \(t=1\) the traces agree. Equation
		\eqref{eq:5-canonical-relative-spheres} gives
		\[
		1-\frac{\varepsilon}{100}
		\le\frac{(\vint_{\partial B_t}H_{i,R}^2)^{1/2}}
		{(\vint_{\partial B_t}u_{x_i,R}^2)^{1/2}}
		\le1+\frac{\varepsilon}{100}
		\]
		at \(t=1/100,1/50,1/2\), while this quotient equals one at
		\(t=1\). Thus
		\[
		D^{H_{i,R}}(0,t)-D^{u_{x_i,R}}(0,t)
		=2\log_4
		\frac{
			(\vint_{\partial B_{2t}}H_{i,R}^2)^{1/2}/
			(\vint_{\partial B_{2t}}u_{x_i,R}^2)^{1/2}}
		{(\vint_{\partial B_t}H_{i,R}^2)^{1/2}/
			(\vint_{\partial B_t}u_{x_i,R}^2)^{1/2}},
		\]
		which transfers the pinching at \(1/100\) and \(1/2\):
		\[
		|D^{H_{i,R}}(0,1/100)-d|
		+|D^{H_{i,R}}(0,1/2)-d|\le C\varepsilon.
		\]
		Proposition \ref{prop:uniform-symmetry-under-pinching} therefore applies,
		with degree \(d\), since the selected integer is within
		\(C\varepsilon<1/4\) of \(d\). Its approximating polynomial is
		\begin{equation}\label{eq:5-canonical-scale-polynomial}
			P_{i,R}=\frac{H_{i,R,d}}{\|H_{i,R,d}\|}.
		\end{equation}
		At radius \(1/8\), the uniform symmetry estimate and orthogonality
		of the homogeneous terms give
		\begin{align*}
			C\varepsilon
			&\ge\left\|
			\frac{H_{i,R}(\cdot/8)}{(\vint_{\partial B_{1/8}}H_{i,R}^2)^{1/2}}
			-P_{i,R}\right\|^2\\
			&=2-2\frac{8^{-d}\|H_{i,R,d}\|}
			{(\vint_{\partial B_{1/8}}H_{i,R}^2)^{1/2}}.
		\end{align*}
		Using \eqref{eq:5-canonical-relative-spheres} at \(t=1/8\) and
		\eqref{eq:5-canonical-sphere-lower}, we obtain
		\[
		\|H_{i,R,d}\|
		\ge c\,8^d\left(\vint_{\partial B_{1/8}}u_{x_i,R}^2\right)^{1/2}
		\ge cA_{\mathcal D}^{-1/2}8^{-(\varepsilon+\varepsilon_1)}.
		\]
		The degree-dependent factors cancel. Since
		\(\varepsilon+\varepsilon_1<1/5\), this proves
		\begin{equation}\label{eq:5-canonical-degree-lower}
			\|H_{i,R,d}\|\ge cA_{\mathcal D}^{-1/2}.
		\end{equation}
		In particular the degree-\(d\) component is nonzero. The exact rescaling
		identity and the normalization inequality give
		\begin{align*}
			&\left\|u_{x_i,R/8}
			-\frac{H_{i,R}(\cdot/8)}
			{(\vint_{\partial B_{1/8}}H_{i,R}^2)^{1/2}}\right\|\\
			&\qquad\le
			2\frac{(\vint_{\partial B_{1/8}}|u_{x_i,R}-H_{i,R}|^2)^{1/2}}
			{(\vint_{\partial B_{1/8}}u_{x_i,R}^2)^{1/2}}
			\le\frac{\varepsilon}{50}.
		\end{align*}
		Together with the uniform symmetry estimate at radius \(1/8\), this yields
		\[
		\|u_{x_i,R/8}-P_{i,R}\|
		\le\frac{\varepsilon}{50}+C\varepsilon^{1/2}.
		\]
		Squaring proves \eqref{eq:5-canonical-scale-approximation}.
		
		For adjacent scales define
		\[
		\widehat H_{i,R}(\xi)
		=\frac{H_{i,R}(\xi/2)}
		{(\vint_{\partial B_{1/2}}H_{i,R}^2)^{1/2}}.
		\]
		The boundary value of \(H_{i,R/2}\) is
		\(u_{x_i,R}(\xi/2)/(\vint_{\partial B_{1/2}}u_{x_i,R}^2)^{1/2}\).
		Thus the normalization inequality gives
		\begin{align*}
			\|\widehat H_{i,R}-H_{i,R/2}\|
			&\le2\frac{(\vint_{\partial B_{1/2}}
				|H_{i,R}-u_{x_i,R}|^2)^{1/2}}
			{(\vint_{\partial B_{1/2}}u_{x_i,R}^2)^{1/2}}\\
			&\le C K A_{\mathcal D}^{C}2^{d_{\mathcal D}+1}
			\bigl(\omega(2\sqrt\Lambda R)+MR^{2-n/p}\bigr),
		\end{align*}
		where \eqref{eq:5-canonical-replacement-energy} bounds the numerator and
		\eqref{eq:5-canonical-sphere-lower} bounds the denominator.
		Let \(\Pi_d\) denote orthogonal projection onto \(\mathcal P_d\).
		The exact identities
		\[
		\Pi_d\widehat H_{i,R}
		=\frac{2^{-d}H_{i,R,d}}
		{(\vint_{\partial B_{1/2}}H_{i,R}^2)^{1/2}},
		\qquad \Pi_dH_{i,R/2}=H_{i,R/2,d}
		\]
		and contractivity give
		\begin{equation}\label{eq:5-canonical-adjacent-degree-components}
			\left\|\frac{2^{-d}H_{i,R,d}}
			{(\vint_{\partial B_{1/2}}H_{i,R}^2)^{1/2}}
			-H_{i,R/2,d}\right\|
			\le\|\widehat H_{i,R}-H_{i,R/2}\|.
		\end{equation}
		Both components are nonzero. The first is a positive multiple of
		\(P_{i,R}\), while the normalized second component is \(P_{i,R/2}\).
		The normalization inequality and \eqref{eq:5-canonical-degree-lower},
		used at \(R/2\), imply
		\begin{equation}\label{eq:5-canonical-normalization-final}
			\begin{aligned}
				\|P_{i,R}-P_{i,R/2}\|
				&\le2\frac{
					\|\Pi_d\widehat H_{i,R}-H_{i,R/2,d}\|}
				{\|H_{i,R/2,d}\|}\\
				&\le C A_{\mathcal D}^{1/2}
				\|\widehat H_{i,R}-H_{i,R/2}\|.
			\end{aligned}
		\end{equation}
		Absorbing this additional fixed power of \(A_{\mathcal D}\) proves
		\eqref{eq:5-canonical-adjacent-polynomials}.
	\end{proof}
	
	For the dyadic scales used below, put
	\begin{equation}\label{eq:5-canonical-dyadic-scales}
		R_{i,m}=2^m(200\rho_i),\qquad P_i=P_{i,R_{i,0}}.
	\end{equation}
	
	\begin{corollary}\label{cor:5-canonical-coherence}
		Under the hypotheses of Lemma \ref{lem:5-canonical-polynomial-center},
		\begin{equation}\label{eq:5-canonical-coherence}
			\|P_{i,R_{i,m}}-P_i\|
			\le C A_{\mathcal D}^{C}2^{d_{\mathcal D}+1}
			\eta(2R_{i,m})\le\frac{\varepsilon}{100}
		\end{equation}
		whenever \(R_{i,m}\) lies in \eqref{eq:5-canonical-scale-range}.
	\end{corollary}
	
	\begin{proof}
		Summation over the dyadic radii gives
		\[
		\sum_{\ell=1}^m\!
		\bigl(\omega(2\sqrt\Lambda R_{i,\ell})
		+MR_{i,\ell}^{2-n/p}\bigr)
		\le C\left(\int_0^{4\sqrt\Lambda R_{i,m}}
		\frac{\omega(t)}t\,dt
		+MR_{i,m}^{2-n/p}\right).
		\]
		Indeed, monotonicity of \(\omega\) gives
		\[
		\omega(2\sqrt\Lambda R_{i,\ell})
		\le\frac1{\log2}\int_{2\sqrt\Lambda R_{i,\ell}}^{4\sqrt\Lambda R_{i,\ell}}
		\frac{\omega(t)}t\,dt,
		\]
		and these intervals have disjoint interiors. Also,
		\[
		\sum_{\ell=1}^m MR_{i,\ell}^{2-n/p}
		=MR_{i,m}^{2-n/p}
		\sum_{q=0}^{m-1}2^{-q(2-n/p)}
		\le CMR_{i,m}^{2-n/p}.
		\]
		Thus the adjacent-scale errors are summable with no factor depending on
		the number of scales, and
		\begin{align*}
			\|P_{i,R_{i,m}}-P_i\|
			&\le C K A_{\mathcal D}^{C}2^{d_{\mathcal D}+1}
			\left(\int_0^{4\sqrt\Lambda R_{i,m}}\frac{\omega(t)}t\,dt
			+MR_{i,m}^{2-n/p}\right)\\
			&\le C A_{\mathcal D}^{C}2^{d_{\mathcal D}+1}
			\eta(2R_{i,m})\le\frac{\varepsilon}{100}.
		\end{align*}
		This is \eqref{eq:5-canonical-coherence}.
	\end{proof}
	
	The next two lemmas give the basic Lipschitz character of singular set, which are the main ingredients for volume estimate.
	\begin{lemma}\label{lem:5-coherent-packing}
		Let \(2\le d\le d_{\mathcal D}\), \(0<\sigma\le4r_8\), and
		let \(\{(x_i,\rho_i)\}_{i\in I}\) be a finite or countable family with
		\[
		x_i\in S(u)\cap B_{\sigma/2}(x_0),\qquad
		0<\rho_i\le\frac{\sigma}{10^5C_3}.
		\]
		Assume that \(B_{200\rho_i}(x_i)\) are pairwise disjoint and that
		\eqref{eq:5-canonical-center-pinching} holds at every \(x_i\).
		There is a subspace \(W\), with \(\dim W\le n-2\), such that
		\begin{equation}\label{eq:5-fixed-cone}
			\operatorname{dist}(x_j-x_i,W)\le10^{-12}|x_j-x_i|
			\qquad(i,j\in I).
		\end{equation}
		The centers lie in an affine translate of a graph over \(W\) with
		Lipschitz constant at most \(10^{-6}\), and
		\begin{equation}\label{eq:5-coherent-packing-bound}
			\sum_i\rho_i^{n-2}\le C(n,\lambda,\Lambda)\sigma^{n-2}.
		\end{equation}
	\end{lemma}
	
	\begin{proof}
		If \(I=\varnothing\), take \(W=\{0\}\); the graph and cone assertions are
		vacuous and \(\sum_i\rho_i^{n-2}=0\). Hence assume \(I\ne\varnothing\).
		Fix \(i_0\in I\),
		let \(P_i\) be supplied by Lemma \ref{lem:5-canonical-polynomial-center},
		and put \(P=P_{i_0}\).
		We first compare the polynomials at two distinct centers. Put
		\(L=|x_i-x_j|\). Disjointness gives
		\[
		200(\rho_i+\rho_j)\le L\le\sigma.
		\]
		The pinching at both centers therefore holds throughout
		\([L/100,C_3L]\). Apply the harmonic comparison
		\eqref{eq:cone-splitting-long-doubling-comparison}--
		\eqref{eq:cone-splitting-long-blowup-comparison}, with base scale
		\(R=20(1+\lambda^{-1/2})L\). The two harmonic centers are \(0\) and
		\[
		z=\frac{A_{x_i}^{-1}(x_j-x_i)}{R}.
		\]
		Their doubling indices are within \(C\varepsilon\) of \(d\) for
		\(1/10\le s\le9d\), since
		\[
		[R/10,9dR]\subset[L/100,C_3L].
		\]
		Let \(Q_i,Q_j\) be the normalized degree-\(d\) parts at these centers.
		Lemma \ref{lem:two-center-coefficient-estimate} gives
		\begin{equation}\label{eq:5-pair-harmonic-estimates}
			\|Q_i-Q_j\|\le C\varepsilon^{1/2},\qquad
			\|z\cdot\nabla Q_i\|
			\le C\varepsilon^{1/2}\|\nabla Q_i\|.
		\end{equation}
		
		For each center choose a dyadic radius \(R_{i,m}\in[8R,16R)\).
		This choice is possible because \(8R\ge200\rho_i\), and it is admissible
		since
		\[
		R_{i,m}<320(1+\lambda^{-1/2})L
		\le400(1+\lambda^{-1/2})\sigma.
		\]
		At the scale \(R_{i,m}/8\), the normalized harmonic rescaling
		at the first center is evaluated at a radius in \([1,2)\).
		Equations \eqref{eq:5-canonical-coherence},
		\eqref{eq:5-canonical-scale-approximation}, and
		\eqref{eq:cone-splitting-long-blowup-comparison}, together with
		Proposition~\ref{prop:uniform-symmetry-under-pinching}, give
		\begin{align*}
			\|P_i-Q_i\|
			&\le\|P_i-P_{i,R_{i,m}}\|
			+\|P_{i,R_{i,m}}-u_{x_i,R_{i,m}/8}\|\\
			&\quad+\|u_{x_i,R_{i,m}/8}-h_{0,R_{i,m}/(8R)}\|
			+\|h_{0,R_{i,m}/(8R)}-Q_i\|\\
			&\le \frac{\varepsilon}{100}
			+C\varepsilon^{1/2}+C\varepsilon+C\varepsilon^{1/2}
			\le C\varepsilon^{1/2}.
		\end{align*}
		The same four estimates, with \(R_{j,\ell}\in[8R,16R)\) and
		harmonic center \(z\), give
		\begin{align*}
			\|P_j-Q_j\|
			&\le\|P_j-P_{j,R_{j,\ell}}\|
			+\|P_{j,R_{j,\ell}}-u_{x_j,R_{j,\ell}/8}\|\\
			&\quad+\|u_{x_j,R_{j,\ell}/8}
			-h_{z,R_{j,\ell}/(8R)}\|
			+\|h_{z,R_{j,\ell}/(8R)}-Q_j\|\\
			&\le\frac{\varepsilon}{100}
			+C\varepsilon^{1/2}+C\varepsilon+C\varepsilon^{1/2}
			\le C\varepsilon^{1/2}.
		\end{align*}
		All norms and comparisons use
		the same variable on \(\partial B_1\). The gradient identity for
		homogeneous harmonic polynomials implies
		\[
		\|\nabla(P_i-Q_i)\|
		=\sqrt{d(2d+n-2)}\,\|P_i-Q_i\|,
		\qquad \|\nabla P_i\|=\sqrt{d(2d+n-2)}.
		\]
		Since \(|z|\) is bounded above and below by positive constants depending
		only on ellipticity, \eqref{eq:5-pair-harmonic-estimates} yields
		\begin{equation}\label{eq:5-polynomial-coherence-with-P}
			\|P_i-P_j\|\le C\varepsilon^{1/2},\qquad
			\left\|\left(A_{x_i}^{-1}\frac{x_j-x_i}{L}\right)
			\cdot\nabla P_i\right\|
			\le C\varepsilon^{1/2}\|\nabla P_i\|.
		\end{equation}
		The Lipschitz dependence of the positive matrix square root and its
		inverse on uniformly elliptic matrices gives
		\[
		\|A_{x_i}^{-1}-A_{x_{i_0}}^{-1}\|
		\le C(\lambda,\Lambda)\omega(2\sigma).
		\]
		To pass to the fixed polynomial and fixed metric, the three errors are
		\begin{align}
			\left\|\left(A_{x_{i_0}}^{-1}\frac{x_j-x_i}{L}\right)
			\cdot\nabla P\right\|
			&\le\left\|\left(A_{x_i}^{-1}\frac{x_j-x_i}{L}\right)
			\cdot\nabla P_i\right\|\notag\\
			&\quad+\left\|\left(A_{x_{i_0}}^{-1}\frac{x_j-x_i}{L}\right)
			\cdot\nabla(P-P_i)\right\|\notag\\
			&\quad+\left\|\left((A_{x_{i_0}}^{-1}-A_{x_i}^{-1})
			\frac{x_j-x_i}{L}\right)
			\cdot\nabla P_i\right\|\notag\\
			&\le C\bigl(\varepsilon^{1/2}+\omega(2\sigma)\bigr)
			\|\nabla P\|.
			\label{eq:5-fixed-direction-three-errors}
		\end{align}
		Here the second term uses \eqref{eq:5-polynomial-coherence-with-P}
		with the pair \(i,i_0\). Since \(2\sigma\le8r_8\),
		\eqref{eq:5-rmaster} gives
		\[
		C A_{\mathcal D}^{C}\omega(2\sigma)
		\le\frac{\tau\varepsilon^{1/2}}{10^8}
		\le\varepsilon^{1/2}.
		\]
		Ellipticity gives
		\(L\le\sqrt\Lambda|A_{x_{i_0}}^{-1}(x_j-x_i)|\). Multiplying
		\eqref{eq:5-fixed-direction-three-errors} by \(L\), we obtain
		\begin{equation}\label{eq:5-pairwise-fixed-polynomial-preliminary}
			\|(A_{x_{i_0}}^{-1}(x_j-x_i))\cdot\nabla P\|
			\le C\varepsilon^{1/2}\|\nabla P\|
			|A_{x_{i_0}}^{-1}(x_j-x_i)|.
		\end{equation}
		
		Set
		\begin{equation}\label{eq:5-gamma-definition}
			\gamma=\frac{\varepsilon^{1/4}}{10^4(1+d_{\mathcal D})^{n-2}}.
		\end{equation}
		The choice of \(\varepsilon\) ensures, for \(2\le d\le d_{\mathcal D}\),
		\begin{equation}\label{eq:5-gamma-consequences}
			\begin{gathered}
				4\gamma^2\le c(n)d^{-4(n-1)},\qquad
				\gamma^2<\frac1{2n(n-1)},\\
				C(n)d^{n-1}\gamma\le\tfrac14,\\
				C\varepsilon^{1/2}
				\le10^{-12}\sqrt{\lambda/\Lambda}\,\gamma.
			\end{gathered}
		\end{equation}
		Consider the positive semidefinite matrix
		\[
		(G_P)_{\alpha\beta}
		=\vint_{\partial B_1}\partial_\alpha P\,\partial_\beta P.
		\]
		For every vector \(v\),
		\[
		v^TG_Pv=\|v\cdot\nabla P\|^2,
		\qquad \operatorname{tr}G_P=\|\nabla P\|^2.
		\]
		Let \(e_\alpha\) be an orthonormal eigenbasis with eigenvalues
		\(\mu_\alpha\ge0\), and let \(W(P)\) be the span of the eigenvectors
		with \(\mu_\alpha\le\gamma^2\|\nabla P\|^2\). Then
		\[
		\|v\cdot\nabla P\|^2
		=\sum_{\alpha=1}^n\mu_\alpha(v\cdot e_\alpha)^2.
		\]
		Discarding the summands corresponding to \(W(P)\) gives
		\begin{equation}\label{eq:5-Gram-distance-lower}
			\|v\cdot\nabla P\|^2
			\ge\gamma^2\|\nabla P\|^2
			\operatorname{dist}(v,W(P))^2.
		\end{equation}
		We have \(\dim W(P)\le n-2\). Indeed, if
		\(e_1,\ldots,e_{n-1}\in W(P)\), then
		\[
		\|e_\alpha\cdot\nabla P\|^2
		=\mu_\alpha\le\gamma^2\|\nabla P\|^2
		\qquad(1\le\alpha\le n-1).
		\]
		By \eqref{eq:5-gamma-consequences},
		\cite[Lemma~3.23]{NaberValtorta2017} then implies \(d=1\),
		contrary to \(d\ge2\).
		
		Define \(W=A_{x_{i_0}}W(P)\). Equations
		\eqref{eq:5-pairwise-fixed-polynomial-preliminary}--
		\eqref{eq:5-Gram-distance-lower} give
		\[
		\operatorname{dist}\bigl(A_{x_{i_0}}^{-1}(x_j-x_i),W(P)\bigr)
		\le C\varepsilon^{1/2}\gamma^{-1}
		|A_{x_{i_0}}^{-1}(x_j-x_i)|
		\le10^{-12}\sqrt{\lambda/\Lambda}\,
		|A_{x_{i_0}}^{-1}(x_j-x_i)|.
		\]
		Since
		\(A_{x_{i_0}}\pi_{W(P)}A_{x_{i_0}}^{-1}(x_j-x_i)\in W\),
		ellipticity therefore gives
		\begin{align*}
			\operatorname{dist}(x_j-x_i,W)
			&\le\left|A_{x_{i_0}}\left(
			A_{x_{i_0}}^{-1}(x_j-x_i)
			-\pi_{W(P)}A_{x_{i_0}}^{-1}(x_j-x_i)\right)\right|\\
			&\le\sqrt\Lambda\,
			\operatorname{dist}\bigl(A_{x_{i_0}}^{-1}(x_j-x_i),W(P)\bigr)\\
			&\le10^{-12}\sqrt\lambda\,
			|A_{x_{i_0}}^{-1}(x_j-x_i)|
			\le10^{-12}|x_j-x_i|,
		\end{align*}
		which is \eqref{eq:5-fixed-cone}.
		
		For \(i\ne j\), orthogonal projection onto \(W\) satisfies
		\[
		|\pi_W(x_j-x_i)|\ge\sqrt{1-10^{-24}}\,|x_j-x_i|,
		\qquad
		|\pi_{W^\perp}(x_j-x_i)|\le10^{-12}|x_j-x_i|.
		\]
		It is therefore injective on the centers, and
		\[
		\frac{|\pi_{W^\perp}(x_j-x_i)|}{|\pi_W(x_j-x_i)|}
		\le\frac{10^{-12}}{\sqrt{1-10^{-24}}}<10^{-11}.
		\]
		Thus the map sending \(\pi_Wx_i\) to \(\pi_{W^\perp}x_i\) is
		Lipschitz with constant at most \(10^{-11}\), and extension theorem
		extends it to \(W\).
		If \(k=\dim W\ge1\), the balls
		\(B^W_{50\rho_i}(\pi_Wx_i)\) are disjoint and contained in a ball of
		radius \(C\sigma\). Comparing their \(k\)-dimensional volumes gives
		\[
		c(k)50^k\sum_i\rho_i^k
		\le C(k)\sigma^k,
		\]
		and therefore
		\[
		\sum_i\rho_i^{n-2}
		\le \sigma^{n-2-k}\sum_i\rho_i^k
		\le C(n)\sigma^{n-2}.
		\]
		If \(k=0\), injectivity implies that there is at most one center, and
		the same estimate follows from \(\rho_i\le\sigma\).
	\end{proof}
	
	\begin{lemma}\label{lem:5-coherent-confinement}
		Assume the hypotheses of Lemma \ref{lem:5-coherent-packing}, and use
		the polynomials \(P_i\) and subspace \(W=A_{x_{i_0}}W(P)\) constructed
		in its proof. Suppose \(\dim W=n-2\). If \(y\in S(u)\),
		\begin{equation}\label{eq:5-coherent-exterior-point}
			10^3\rho_i\le L=|y-x_i|\le\sigma,
		\end{equation}
		and
		\begin{equation}\label{eq:5-coherent-base-pinching-for-y}
			d-\varepsilon<D^u(x_i,s)\le d+\varepsilon
			\qquad(L/100\le s\le C_3L),
		\end{equation}
		then
		\begin{equation}\label{eq:5-coherent-confinement}
			\operatorname{dist}(y-x_i,W)\le\frac{\tau}{10^4}L.
		\end{equation}
	\end{lemma}
	
	\begin{proof}
		Put \(s_L=\frac{20}{3}(1+\lambda^{-1/2})L\) and \(r_L=s_L/3\).
		For \(L/100\le s\le C_3L\),
		\[
		|D^u(x_i,s)-D^u(x_i,C_3L)|\le2\varepsilon.
		\]
		For this application,
		\[
		\frac{L/100}{C_3L}=\frac1{100C_3}\le10^{-2}.
		\]
		The definition of \(r_6(\varepsilon)\), with its fixed power chosen to
		dominate \((100C_3)^{\varepsilon_1}\), and
		\eqref{eq:5-rmaster} give
		\begin{align*}
			C_3L&\le4C_3r_8
			\le\min\{r_4(1/10),(10\sqrt\Lambda)^{-1}\},\\
			A_{\mathcal D}^{1/2}(100C_3)^{\varepsilon_1}\eta(C_3L)
			&\le c_3(2\varepsilon).
		\end{align*}
		Here \(c_3\) is the constant supplied by Theorem
		\ref{thm:uniform-symmetry-pinched}.
		Theorem \ref{thm:uniform-symmetry-pinched}, used with parameter
		\(2\varepsilon\), therefore supplies a normalized
		\(\widetilde P_i\in\mathcal P_d\) such that
		\begin{equation}\label{eq:5-confinement-uniform-polynomial}
			\|u_{x_i,t}-\widetilde P_i\|^2\le C\varepsilon
			\qquad(L/25\le t\le C_3L/12).
		\end{equation}
		Compare this polynomial with \(P_i\) at a dyadic scale
		\[
		160(1+\lambda^{-1/2})L\le R_{i,m}
		<320(1+\lambda^{-1/2})L.
		\]
		The bounds
		\[
		20(1+\lambda^{-1/2})L\le \frac{R_{i,m}}8
		<40(1+\lambda^{-1/2})L
		\]
		put both \(R_{i,m}/8\) and \(s_L\) in the interval of
		\eqref{eq:5-confinement-uniform-polynomial}. The scale \(R_{i,m}\) is
		admissible in Lemma \ref{lem:5-canonical-polynomial-center} because
		\[
		R_{i,m}\ge160(1+\lambda^{-1/2})10^3\rho_i>200\rho_i,
		\qquad
		R_{i,m}<320(1+\lambda^{-1/2})\sigma.
		\]
		Thus \eqref{eq:5-canonical-scale-approximation},
		\eqref{eq:5-canonical-coherence}, and
		\eqref{eq:5-confinement-uniform-polynomial} give
		\[
		\|\widetilde P_i-P_i\|
		\le\|\widetilde P_i-u_{x_i,R_{i,m}/8}\|
		+\|u_{x_i,R_{i,m}/8}-P_{i,R_{i,m}}\|
		+\|P_{i,R_{i,m}}-P_i\|
		\le C\varepsilon^{1/2}.
		\]
		Consequently,
		\begin{equation}\label{eq:5-confinement-u-Pi}
			\begin{aligned}
				\|u_{x_i,s_L}-P_i\|^2
				&\le2\|u_{x_i,s_L}-\widetilde P_i\|^2
				+2\|\widetilde P_i-P_i\|^2\\
				&\le C\varepsilon.
			\end{aligned}
		\end{equation}
		
		For an orthonormal basis \(e_\alpha\) of \(W(P)\),
		\eqref{eq:5-polynomial-coherence-with-P} and the gradient identity give
		\[
		\|e_\alpha\cdot\nabla P_i\|
		\le(\gamma+C\varepsilon^{1/2})\|\nabla P_i\|
		\le2\gamma\|\nabla P_i\|.
		\]
		If \(W(P)=\{0\}\), put \(\bar P_i=P_i\) and \(P_i^\perp=0\).
		Otherwise, since \(\dim W(P)\le n-2\), Lemma
		\ref{lem:almost-invariant-polynomial} and
		\eqref{eq:5-gamma-consequences} yield
		\(P_i=\bar P_i+P_i^\perp\), where \(\bar P_i\)
		is invariant along \(W(P)\) and
		\(\|P_i^\perp\|\le C d^{n-2}\gamma\le C\varepsilon^{1/4}\).
		The same bound is immediate when \(W(P)=\{0\}\).
		Since \(\|P_i\|=1\), the reverse triangle inequality gives
		\[
		\bigl|\|\bar P_i\|-1\bigr|
		\le\|P_i^\perp\|\le\tfrac14,
		\qquad \|\bar P_i\|\ge\tfrac34.
		\]
		Its normalization \(Q_i=\bar P_i/\|\bar P_i\|\) satisfies
		\begin{align*}
			\|P_i-Q_i\|
			&\le\|P_i^\perp\|
			+\left\|\bar P_i-\frac{\bar P_i}{\|\bar P_i\|}\right\|\\
			&=\|P_i^\perp\|+\bigl|\|\bar P_i\|-1\bigr|
			\le2\|P_i^\perp\|,
		\end{align*}
		and hence
		\[
		\|P_i-Q_i\|^2
		\le C d^{2(n-2)}\gamma^2
		\le C\varepsilon^{1/2}.
		\]
		When \(n=2\), \(W(P)=\{0\}\) and we simply take \(Q_i=P_i\).
		Thus, since \(3r_L=s_L\),
		\[
		\|u_{x_i,3r_L}-Q_i\|^2
		\le2\|u_{x_i,s_L}-P_i\|^2+2\|P_i-Q_i\|^2
		\le C(\varepsilon+\varepsilon^{1/2})
		\le C\varepsilon^{1/2}.
		\]
		Put \(\varepsilon'=C\varepsilon^{1/2}\), with the fixed constant in the
		last display. By \eqref{eq:5-epsilon-choice},
		\[
		\varepsilon\le\varepsilon'\le
		\min\{\varepsilon_1/4,\varepsilon_4(\tau_1)\}.
		\]
		The defining conditions of \(r_6\) are monotone in their error
		parameter, while
		\[
		r_L\le C\sigma\le4Cr_8
		<4\cdot10^8C_3r_8\le r_6(\varepsilon)\le r_6(\varepsilon').
		\]
		Thus all scale and error hypotheses of Lemma
		\ref{lem:singular-confinement} hold with \(\tau=\tau_1\) and error
		\(\varepsilon'\).
		Also
		\[
		|A_{x_i}^{-1}(y-x_i)|\le\lambda^{-1/2}L
		<\frac{20}{9}(1+\lambda^{-1/2})L=r_L.
		\]
		Thus \(y\in x_i+r_LA_{x_i}B_1\), and that lemma gives
		\[
		\operatorname{dist}(y-x_i,A_{x_i}W(P))
		\le\sqrt\Lambda\,\tau_1r_L.
		\]
		We must express this estimate in the fixed plane \(W\) in the original coordinates.
		For a unit vector \(v\in A_{x_i}W(P)\), write \(v=A_{x_i}w\), where
		\(w\in W(P)\) and \(|w|\le\lambda^{-1/2}\). Thus
		\[
		\operatorname{dist}(v,W)
		\le |(A_{x_i}-A_{x_{i_0}})w|
		\le\lambda^{-1/2}\|A_{x_i}-A_{x_{i_0}}\|.
		\]
		Interchanging the two matrices gives the reverse estimate. Since the
		image subspaces have the same dimension, these estimates imply
		\begin{equation}\label{eq:5-confinement-plane-angle}
			\|\pi_{A_{x_i}W(P)}-\pi_W\|
			\le C(\lambda,\Lambda)\|A_{x_i}-A_{x_{i_0}}\|
			\le C\omega(2\sigma)\le C\varepsilon^{1/2}.
		\end{equation}
		Here we used the Lipschitz dependence of the positive square root on
		uniformly positive definite matrices, followed by \eqref{eq:5-rmaster}.
		For every vector \(z\),
		\[
		\operatorname{dist}(z,W)
		\le\operatorname{dist}(z,A_{x_i}W(P))
		+\|\pi_{A_{x_i}W(P)}-\pi_W\|\,|z|.
		\]
		Applying this with \(z=y-x_i\), using \(r_L=\frac{20}{9}
		(1+\lambda^{-1/2})L\), gives
		\[
		\operatorname{dist}(y-x_i,W)
		\le\left[\frac{20\sqrt\Lambda}{9}(1+\lambda^{-1/2})
		\tau_1+C\varepsilon^{1/2}\right]L
		\le\frac{\tau}{10^4}L,
		\]
		as required.
	\end{proof}
	
	We formulate the covering step for a subset of the singular set. This
	allows each child ball to carry precisely the part of its parent set
	which it covers.
	
	\begin{lemma}\label{lem:5-covering}
		Let \(2\le d\le d_{\mathcal D}\), \(0<\rho\le r_8\), and
		\(E\subset S(u)\cap B_\rho(x_0)\). Assume
		\begin{equation}\label{eq:5-good-ball}
			D^u(y,s)\le d+\varepsilon
			\qquad(y\in E,\ 0<s\le4C_3\rho).
		\end{equation}
		For every \(0<r\le\rho\), there are balls \(B_{s_i}(x_i)\), with
		\(x_i\in E\), and subsets \(E_i\subset E\cap B_{s_i}(x_i)\), such that
		\begin{equation}\label{eq:5-covering-inclusion}
			E=\bigcup_i E_i,\qquad
			\sum_i s_i^{n-2}\le\mathcal A\rho^{n-2}.
		\end{equation}
		Each child either has \(s_i=r\), or satisfies
		\begin{equation}\label{eq:5-covering-nonterminal-scale}
			r<s_i\le\theta\rho,\qquad
			D^u(y,s)\le d-1+\varepsilon
			\quad(y\in E_i,\ 0<s\le4C_3s_i).
		\end{equation}
	\end{lemma}
	
	\begin{proof}
		If \(E=\varnothing\), the empty family has union \(E\) and total mass
		zero. Hence assume \(E\ne\varnothing\). If \(r\ge\theta\rho\),
		a maximal \(r\)-separated subset of \(E\) has cardinality
		\[
		N|B_{r/2}|\le |B_{2\rho}|,
		\qquad N\le C(n)(\rho/r)^n.
		\]
		The corresponding radius-\(r\) balls cover \(E\), and
		\[
		Nr^{n-2}\le C(n)\rho^nr^{-2}
		\le C(n)\theta^{-2}\rho^{n-2}\le\mathcal A\rho^{n-2}.
		\]
		We may therefore assume
		\begin{equation}\label{eq:5-target-small}
			0<r<\theta\rho.
		\end{equation}
		All stopping radii and all sets in this proof are defined on \(E\).
		For \(x\in E\), put
		\begin{equation}\label{eq:5-stopping-radius}
			\begin{aligned}
				r'_x&=\sup\{0<s\le4C_3\rho:
				D^u(x,s)\le d-\varepsilon\},\\
				r_x&=\max\{r,\min\{r'_x,\rho\}\},
			\end{aligned}
		\end{equation}
		where the supremum of the empty set is zero.
		Continuity of the spherical integrals shows that a positive \(r'_x\)
		is attained. Consequently, if \(r_x>r\), there is
		\(\widehat r_x\in[r_x,4C_3\rho]\) with
		\begin{equation}\label{eq:5-stopping-drop-scale}
			D^u(x,\widehat r_x)\le d-\varepsilon.
		\end{equation}
		If \(r_x<\rho\), then \(r'_x\le r_x\), and hence
		\begin{equation}\label{eq:5-stopping-lower-pinching}
			d-\varepsilon<D^u(x,s)\le d+\varepsilon
			\qquad(r_x<s\le4C_3\rho).
		\end{equation}
		Decompose \(E=E_a\cup E_b\), where
		\begin{equation}\label{eq:5-Sa}
			E_a=\{x\in E:r_y\ge r_x/7
			\text{ for all }y\in E\cap B_{1000r_x}(x)\},
			\qquad E_b=E\setminus E_a.
		\end{equation}
		
		\smallskip
		\noindent\textit{The primary covering.}
		The Vitali lemma gives points \(x_i\in E_a\), with \(\rho_i=r_{x_i}\),
		such that \(B_{200\rho_i}(x_i)\) are disjoint and
		\begin{equation}\label{eq:5-Sa-Vitali}
			E_a\subset\bigcup_i B_{1000\rho_i}(x_i).
		\end{equation}
		This family is finite because \(\rho_i\ge r\).
		Set \(\sigma=4\rho\), and split the indices into
		\[
		I_1=\left\{i:\rho_i\le
		\frac{\sigma}{10^5C_3}\right\},
		\qquad I_2=I\setminus I_1.
		\]
		The disjoint balls lie in \(B_{201\rho}(x_0)\), so
		\(\sum_i \rho_i^n\le C(n)\rho^n\).
		For the large indices, \(\rho_i>\sigma/(10^5C_3)\).
		Hence
		\begin{equation}\label{eq:5-primary-large}
			\sum_{I_2}\rho_i^{n-2}
			\le\left(\frac{10^5C_3}{\sigma}\right)^2
			\sum_i \rho_i^n
			\le C(n)C_3^2\rho^{n-2}.
		\end{equation}
		The same disjointness also gives
		\(\#I_2\le C(n)C_3^n\).
		For the small indices, \(\rho_i<\rho\), and
		\eqref{eq:5-stopping-lower-pinching} gives the hypotheses of Lemma
		\ref{lem:5-coherent-packing} on
		\([2\rho_i,C_3\sigma]\).
		Its scale condition is satisfied because
		\(\sigma=4\rho\le4r_8\).
		Thus
		\begin{equation}\label{eq:5-primary-packing}
			\sum_i \rho_i^{n-2}\le C(n,\lambda,\Lambda)
			C_3^2\rho^{n-2}.
		\end{equation}
		Fix the graph and subspace \(W\) supplied for the small family, when
		that family is nonempty.
		
		For each \(i\), consider \(E\cap B_{1000\rho_i}(x_i)\).
		If \(\theta \rho_i\le r\), choose an \(r\)-cover centered in this set.
		Since \(r\le \rho_i\), its number of members is at most
		\(C(n)(\rho_i/r)^n\), and its mass satisfies
		\begin{equation}\label{eq:5-primary-terminal-mass}
			N_i r^{n-2}\le C(n)\rho_i^n r^{-2}
			\le C(n)\theta^{-2}\rho_i^{n-2}.
		\end{equation}
		Every radius in this family is exactly \(r\).
		If \(\theta \rho_i>r\), every point \(y\) of this set satisfies
		\(r_y\ge \rho_i/7>r\), by the definition of \(E_a\).
		Apply the definite-drop conclusion of Lemma
		\ref{lem:integer-pinching-drop} at \(\widehat r_y\) from
		\eqref{eq:5-stopping-drop-scale}. It gives
		\[
		D^u(y,s)\le d-1+\varepsilon
		\qquad\left(0<s\le\frac{\varepsilon \rho_i}{56}\right).
		\]
		Since \(4C_3\theta \rho_i=\varepsilon \rho_i/2500\),
		a cover centered in \(E\cap B_{1000\rho_i}(x_i)\), with radius
		\(\theta \rho_i\), has the required degree bound on its assigned subsets.
		It has at most \(C(n)\theta^{-n}\) members, so
		\begin{equation}\label{eq:5-primary-child-mass}
			N_i(\theta \rho_i)^{n-2}
			\le C(n)\theta^{-2}\rho_i^{n-2}.
		\end{equation}
		For this construction the child receives only its intersection with
		\(E\cap B_{1000\rho_i}(x_i)\); no assertion is made about singular points
		outside that parent set. Summing the two alternatives and using
		\eqref{eq:5-primary-packing} bounds their total mass by
		\(C C_3^2\theta^{-2}\rho^{n-2}\).
		
		\smallskip
		\noindent\textit{The secondary covering.}
		Every \(y\in E_b\) is close to a primary center:
		\begin{equation}\label{eq:5-secondary-near-primary}
			|y-x_i|<2200r_y,\qquad \rho_i\le r_y
			\quad\text{for some }i.
		\end{equation}
		Indeed, starting at \(y_0=y\), whenever \(y_m\notin E_a\) choose
		\(y_{m+1}\in E\cap B_{1000r_{y_m}}(y_m)\) with
		\(r_{y_{m+1}}<r_{y_m}/7\). The procedure must end at a point
		\(y_N\in E_a\): otherwise \(r_{y_m}<7^{-m}r_y<r\) for all sufficiently
		large \(m\), whereas \eqref{eq:5-stopping-radius} gives
		\(r_{y_m}\ge r\). Moreover,
		\begin{equation}\label{eq:5-chain-distance}
			|y-y_N|\le1000\sum_{m=0}^{N-1}r_{y_m}
			<\frac{7000}{6}r_y<1200r_y,
			\qquad r_{y_N}<r_y/7.
		\end{equation}
		Choose a primary ball \(B_{1000\rho_i}(x_i)\) containing \(y_N\).
		The defining property of \(x_i\in E_a\) gives
		\(r_{y_N}\ge \rho_i/7\), and hence \(\rho_i<r_y\). Therefore
		\(|y-x_i|<1200r_y+1000\rho_i<2200r_y\), as claimed.
		
		For
		\(y\in E_b\setminus\bigcup_i B_{1000\rho_i}(x_i)\), define
		\[
		t_y=10^{-6}\min_i|y-x_i|.
		\]
		Then \(0<t_y\le r_y/400\). A second Vitali selection gives disjoint
		balls \(B_{t_j}(y_j)\) such that
		\begin{equation}\label{eq:5-secondary-Vitali}
			E_b\setminus\bigcup_i B_{1000\rho_i}(x_i)
			\subset\bigcup_j B_{5t_j}(y_j).
		\end{equation}
		Choose a nearest primary center \(x_{i(j)}\), so that
		\begin{equation}\label{eq:5-secondary-radius-distance}
			|y_j-x_{i(j)}|=10^6t_j,\qquad
			0<t_j\le2\cdot10^{-6}\rho.
		\end{equation}
		We will use the following consequence of the construction:
		\begin{equation}\label{eq:5-rz-lower}
			r_z\ge t_j\qquad(z\in E\cap B_{6t_j}(y_j)).
		\end{equation}
		If \(z\in E_a\), choose a primary ball containing \(z\).
		The distance from \(y_j\) to its center is at least \(10^6t_j\);
		hence \(\rho_i\ge900t_j\), and \(r_z\ge \rho_i/7>t_j\).
		If \(z\in E_b\), use \eqref{eq:5-secondary-near-primary} at \(z\):
		\(10^6t_j\le6t_j+2200r_z\). This proves
		\eqref{eq:5-rz-lower}.
		
		Split the secondary indices according to whether \(i(j)\) belongs to
		\(I_2\) or \(I_1\).
		For the former class,
		\begin{equation}\label{eq:5-secondary-large-radius-lower}
			t_j>\frac{\sigma}{10^8C_3},
		\end{equation}
		because \(y_j\notin B_{1000\rho_{i(j)}}(x_{i(j)})\).
		Disjointness and \eqref{eq:5-secondary-radius-distance} give
		\(\sum_j t_j^n\le C(n)\rho^n\), and consequently
		\begin{align*}
			\sum_{i(j)\in I_2}t_j^{n-2}
			&\le\left(\frac{10^8C_3}{\sigma}\right)^2
			\sum_{i(j)\in I_2}t_j^n\\
			&\le C(n)C_3^2\rho^{n-2}.
		\end{align*}
		
		For \(i(j)\in I_1\), put \(k=\dim W\le n-2\).
		First suppose \(k\le n-3\). Group the radii into
		\(2^{-m-1}\rho<t_j\le2^{-m}\rho\), \(m\ge0\).
		By \eqref{eq:5-secondary-radius-distance}, every corresponding ball
		lies in the \(C2^{-m}\rho\)-neighborhood of the fixed Lipschitz graph
		inside \(B_{3\rho}(x_0)\). Orthogonal projection to \(W\), followed by
		integration in \(W^\perp\), bounds the volume of this neighborhood by
		\[
		C\rho^k(2^{-m}\rho)^{n-k}.
		\]
		If \(N_m\) denotes the number in the fixed group, disjointness and
		\(t_j>2^{-m-1}\rho\) give
		\[
		N_m c(n)(2^{-m}\rho)^n
		\le C\rho^k(2^{-m}\rho)^{n-k},
		\qquad N_m\le C2^{mk}.
		\]
		Consequently,
		\[
		\sum_{j\text{ in group }m}t_j^{n-2}
		\le N_m(2^{-m}\rho)^{n-2}
		\le C2^{-m(n-2-k)}\rho^{n-2}.
		\]
		Since \(n-2-k\ge1\), summation yields
		\[
		\sum_{i(j)\in I_1}t_j^{n-2}
		\le C\rho^{n-2}\sum_{m\ge0}2^{-m(n-2-k)}
		\le C\rho^{n-2}.
		\]
		
		Suppose instead that \(k=n-2\). For each such \(j\),
		\[
		1000\rho_{i(j)}\le |y_j-x_{i(j)}|\le2\rho<\sigma,\qquad
		[|y_j-x_{i(j)}|/100,C_3|y_j-x_{i(j)}|]
		\subset[2\rho_{i(j)},C_3\sigma].
		\]
		Lemma \ref{lem:5-coherent-confinement} applies and yields
		\[
		|\pi_{W^\perp}(y_j-x_{i(j)})|\le100\tau t_j.
		\]
		For two such centers write
		\[
		y_j-y_\ell=(y_j-x_{i(j)})+(x_{i(j)}-x_{i(\ell)})
		+(x_{i(\ell)}-y_\ell).
		\]
		The first and third terms have perpendicular components bounded by
		\(100\tau t_j\) and \(100\tau t_\ell\). For the middle term,
		\eqref{eq:5-fixed-cone} and the triangle inequality give
		\[
		|\pi_{W^\perp}(x_{i(j)}-x_{i(\ell)})|
		\le10^{-12}\bigl(|y_j-y_\ell|+10^6(t_j+t_\ell)\bigr).
		\]
		Consequently,
		\begin{equation}\label{eq:5-secondary-fixed-cone}
			\begin{aligned}
				|\pi_{W^\perp}(y_j-y_\ell)|
				&\le10^{-12}|y_j-y_\ell|
				+(100\tau+10^{-6})(t_j+t_\ell)\\
				&\le\frac1{100}|y_j-y_\ell|.
			\end{aligned}
		\end{equation}
		Indeed, disjointness gives \(t_j+t_\ell\le|y_j-y_\ell|\), and the choice
		of \(\tau\) makes
		\[
		10^{-12}+100\tau+10^{-6}\le10^{-2}.
		\]
		It follows from \eqref{eq:5-secondary-fixed-cone} that
		\[
		|\pi_W(y_j-y_\ell)|
		\ge\sqrt{1-10^{-4}}\,|y_j-y_\ell|
		>\tfrac12(t_j+t_\ell).
		\]
		If \(k\ge1\), the balls \(B^W_{t_j/4}(\pi_Wy_j)\) are therefore
		pairwise disjoint and lie in a ball of radius \(C\rho\); hence
		\[
		c(k)4^{-k}\sum_jt_j^k\le C(k)\rho^k,
		\qquad
		\sum_jt_j^{n-2}=\sum_jt_j^k\le C\rho^k=C\rho^{n-2}.
		\]
		If \(k=0\), the same projection inequality implies that there is at most
		one center, and its contribution is at most \(C\rho^{n-2}\). In both
		cases we obtain, together with the other
		classes,
		\begin{equation}\label{eq:5-secondary-packing}
			\sum_j t_j^{n-2}\le C C_3^2\rho^{n-2}.
		\end{equation}
		
		It remains to replace the secondary balls by children. For
		\(\theta t_j>r\), cover \(E\cap B_{5t_j}(y_j)\) by
		\(C(n)\theta^{-n}\) balls of radius \(\theta t_j\).
		For every point \(z\) assigned to such a ball,
		\[
		|z-y_j|\le5t_j+\theta t_j<6t_j,
		\]
		so \eqref{eq:5-rz-lower} gives \(r_z\ge t_j>r\). The definite-drop
		lemma at \(\widehat r_z\ge r_z\) therefore gives
		\[
		D^u(z,s)\le d-1+\varepsilon
		\qquad\left(0<s\le\frac{\varepsilon\widehat r_z}{8}\right).
		\]
		Since
		\[
		4C_3\theta t_j=\frac{\varepsilon t_j}{2500}
		<\frac{\varepsilon\widehat r_z}{8},
		\]
		the required degree bound holds through the child scale.
		Their total mass is at most
		\[
		C(n)\theta^{-2}\sum_jt_j^{n-2}
		\le C C_3^2\theta^{-2}\rho^{n-2}.
		\]
		
		The indices with \(\theta t_j\le r\) are treated collectively.
		For those associated with large primary centers,
		\eqref{eq:5-secondary-large-radius-lower} gives
		\(r\ge c\theta\rho/C_3\).
		Also
		\[
		\#\{j:i(j)\in I_2\}
		\left(\frac{\sigma}{10^8C_3}\right)^n
		\le\sum_{i(j)\in I_2}t_j^n\le C\rho^n,
		\]
		so their number is at most \(C C_3^n\).
		Cover each \(E\cap B_{5t_j}(y_j)\) by at most
		\(C(n)(1+(t_j/r)^n)\) balls of radius \(r\). Summing gives
		\begin{align*}
			C\left(C_3^n r^{n-2}
			+r^{-2}\sum_jt_j^n\right)
			&\le C C_3^n\rho^{n-2}+C r^{-2}\rho^n\\
			&\le C C_3^n\rho^{n-2}
			+C C_3^2\theta^{-2}\rho^{n-2}\\
			&\le C C_3^n\theta^{-2}\rho^{n-2}.
		\end{align*}
		For the terminal indices associated with small primary centers,
		\(t_j\le r/\theta\). Every point of \(B_{5t_j}(y_j)\) is within
		\((10^6+5)t_j\le Cr/\theta\) of \(x_{i(j)}\), which belongs to the
		fixed graph. Thus their union is contained in
		\begin{equation}\label{eq:5-secondary-terminal-small-tube}
			B_{3\rho}(x_0)\cap B_{Cr/\theta}(\operatorname{graph}f).
		\end{equation}
		We cover this union collectively. In the \(k\) directions of \(W\),
		a mesh of size \(cr\) has at most \(C(\rho/r)^k\) points. In the
		remaining directions, each fiber has length at most \(Cr/\theta\)
		and needs at most \(C\theta^{-(n-k)}\) mesh points. The Lipschitz bound
		on \(f\) makes the resulting mesh a cover by balls of radius \(r/2\).
		Consequently the number of balls is at most
		\[
		C(\rho/r)^k\theta^{-(n-k)}
		\]
		balls of radius \(r/2\). Discard the balls not meeting the part of \(E\)
		under consideration and recenter the others at points of that part;
		they then have radius exactly \(r\). Their mass is bounded by
		\[
		C\theta^{-(n-k)}(\rho/r)^k r^{n-2}
		=C\theta^{-2}
		\left(\frac r{\theta\rho}\right)^{n-2-k}\rho^{n-2}
		\le C\theta^{-2}\rho^{n-2}.
		\]
		This estimate includes \(k=0\) and uses \eqref{eq:5-target-small}.
		
		Combining the primary terminal and nonterminal masses with the three
		secondary classes gives
		\[
		\sum_i s_i^{n-2}
		\le C\bigl(C_3^2+C_3^n\bigr)\theta^{-2}\rho^{n-2}
		\le\mathcal A\rho^{n-2},
		\]
		by \eqref{eq:5-theta}. This proves
		\eqref{eq:5-covering-inclusion}. Assign to each child just the portion of the
		corresponding primary or secondary set which it covers. All preceding
		degree bounds then hold on its assigned subset \(E_i\).
		The nonterminal radii are \(\theta \rho_i\) or \(\theta t_j\), lie above
		\(r\), and are at most \(\theta\rho\). Thus
		\eqref{eq:5-covering-nonterminal-scale} also holds.
	\end{proof}
	
	Now we prove the interior quantitative singular set estimate.
	
	\medskip
	
	\begin{proof}[Proof of Theorem \ref{thm:main-volume-estimate}]
		Here \(q_0,\varepsilon_3,c_5\) are the constants supplied by
		Proposition \ref{prop:low-doubling-regularity}.
		If \(d_{\mathcal D}=1\), fix \(y\in S(u)\). The spherical doubling
		bound gives
		\(D^u(y,4C_3r_8)\le1+\varepsilon\). The definition of \(r_8\) gives
		\[
		4C_3r_8\le r_4(\varepsilon),\qquad
		2C_3r_8\le\min\{r_4(1/10),(10\sqrt\Lambda)^{-1}\},\qquad
		A_{\mathcal D}^{1/2}q_0^{-\varepsilon_1}
		\eta(2C_3r_8)\le c_5\varepsilon.
		\]
		Proposition \ref{prop:low-doubling-regularity} makes \(y\) regular,
		which is impossible for \(y\in S(u)\). Thus \(S(u)=\varnothing\).
		We therefore assume \(d_{\mathcal D}\ge2\).
		Fix \(0<r<r_8\). A maximal
		\(r_8\)-separated family in \(S(u)\) gives a cover by
		balls of radius \(r_8\), centered in \(S(u)\). The balls of radius
		\(r_8/2\) about the same centers are disjoint and lie in \(B_2\), so
		\[
		N_0|B_{r_8/2}|\le|B_2|,
		\qquad N_0\le C(n)r_8^{-n},
		\qquad N_0r_8^{n-2}\le C(n)r_8^{-2}.
		\]
		Associate each ball with
		its intersection with \(S(u)\). The spherical doubling bound gives
		\eqref{eq:5-good-ball} with \(d=d_{\mathcal D}\) on each associated set,
		and the above inequality gives the \((n-2)\)-mass bound. Next we show this process clearly.
		
		Let \(\mathcal F_0\) be this
		initial family, with its associated subsets, and let
		\(\mathcal T_0=\varnothing\). At stage \(j\), the family
		\(\mathcal F_j\) consists of triples \((x_i,s_i,E_i)\) with
		\(s_i>r\), \(E_i\subset S(u)\cap B_{s_i}(x_i)\), and
		\begin{equation}\label{eq:5-generation-active-property}
			D^u(y,s)\le d_{\mathcal D}-j+\varepsilon
			\quad(y\in E_i,\ 0<s\le4C_3s_i).
		\end{equation}
		The family \(\mathcal T_j\) consists of the terminal triples already
		obtained, all with radius exactly \(r\). We maintain
		\begin{equation}\label{eq:5-generation-cover}
			S(u)\subset
			\bigcup_{(x_i,s_i,E_i)\in\mathcal F_j\cup\mathcal T_j}E_i.
		\end{equation}
		Apply Lemma \ref{lem:5-covering} to every nonempty set carried by
		\(\mathcal F_j\), with degree \(d_{\mathcal D}-j\). Define
		\[
		\mathcal T_{j+1}=\mathcal T_j\cup
		\{\text{new children of radius }r\},
		\qquad
		\mathcal F_{j+1}=\{\text{new nonterminal children}\},
		\]
		retaining each child's assigned subset. This preserves
		\eqref{eq:5-generation-cover} and gives
		\eqref{eq:5-generation-active-property} with \(j+1\).
		The subsets are essential here: every application uses only points
		in the set carried by the parent ball.
		
		Writing the sums over the corresponding triples, the covering mass
		estimate gives
		\begin{align}
			\sum_{\mathcal F_{j+1}}s_i^{n-2}
			+\sum_{\mathcal T_{j+1}}s_i^{n-2}
			&\le\mathcal A\sum_{\mathcal F_j}s_i^{n-2}
			+\sum_{\mathcal T_j}s_i^{n-2}\notag\\
			&\le\mathcal A\left(
			\sum_{\mathcal F_j}s_i^{n-2}
			+\sum_{\mathcal T_j}s_i^{n-2}\right).
			\label{eq:5-generation-mass-recursion}
		\end{align}
		We used \(\mathcal A\ge1\) in the second line. Since the corresponding
		initial sum is at most \(Cr_8^{-2}\), induction in
		\eqref{eq:5-generation-mass-recursion} gives
		\[
		\sum_{\mathcal F_j}s_i^{n-2}
		+\sum_{\mathcal T_j}s_i^{n-2}
		\le C\mathcal A^jr_8^{-2}.
		\]
		At stage \(d_{\mathcal D}-1\), every remaining associated set has
		upper doubling index \(1+\varepsilon\) at all required scales.
		Indeed, for \((x_i,s_i,E_i)\in\mathcal F_{d_{\mathcal D}-1}\), fix
		\(y\in E_i\). Then
		\eqref{eq:5-generation-active-property} gives
		\(D^u(y,4C_3s_i)\le1+\varepsilon\). Since \(s_i\le r_8\),
		\eqref{eq:5-rmaster}, the definition of \(r_7\), and the smallness
		condition in \(r_6(\varepsilon)\) give
		\[
		4C_3s_i\le r_4(\varepsilon),\qquad
		2C_3s_i\le\min\{r_4(1/10),(10\sqrt\Lambda)^{-1}\},\qquad
		A_{\mathcal D}^{1/2}q_0^{-\varepsilon_1}
		\eta(2C_3s_i)\le c_5\varepsilon.
		\]
		Proposition \ref{prop:low-doubling-regularity} makes \(y\) regular,
		contradicting \(E_i\subset S(u)\). Thus
		\begin{equation}\label{eq:5-final-active-balls-empty}
			\mathcal F_{d_{\mathcal D}-1}=\varnothing
		\end{equation}
		after empty sets are discarded. The terminal balls cover \(S(u)\),
		and, since each terminal radius is exactly \(r\), their number \(N(r)\)
		satisfies
		\begin{equation}\label{eq:5-covering-number-precise}
			N(r)r^{n-2}\le
			C r_8^{-2}\mathcal A^{d_{\mathcal D}}.
		\end{equation}
		For countable covers the same bound follows by monotone convergence;
		in particular it forces finiteness of the terminal family.
		If \(\operatorname{dist}(z,S(u))<r\), choose \(y\in S(u)\) with
		\(|z-y|<r\), and then a terminal ball \(B_r(x_i)\) containing \(y\).
		Thus \(|z-x_i|<2r\), and
		\[
		\{z\in B_1:\operatorname{dist}(z,S(u))<r\}
		\subset\bigcup_{i=1}^{N(r)}B_{2r}(x_i).
		\]
		Hence
		\begin{equation}\label{eq:5-main-volume-estimate}
			\begin{aligned}
				|\{z\in B_1:\operatorname{dist}(z,S(u))<r\}|
				&\le N(r)|B_{2r}|\\
				&\le C[N(r)r^{n-2}]r^2\\
				&\le C r_8^{-2}\mathcal A^{d_{\mathcal D}}r^2.
			\end{aligned}
		\end{equation}
		For \(r\ge r_8\), since \(\mathcal A\ge1\),
		\[
		|\{z\in B_1:\operatorname{dist}(z,S(u))<r\}|
		\le|B_1|\le C r_8^{-2}\mathcal A^{d_{\mathcal D}}r^2.
		\]
		
		We now specify the dependence of the covering scale. Recall
		\(\delta\) from \eqref{eq:5-final-delta}. With the structural constants
		there chosen sufficiently small and large, respectively, we may take
		\begin{equation}\label{eq:5-explicit-r8}
			r_8
			=\frac{c}{(1+\mathcal D)^2}
			\min\left\{\rho_\omega(\delta),
			\left(\frac{\delta}{1+M}\right)^{1/(2-n/p)},1\right\}.
		\end{equation}
		The choice satisfies all smallness requirements because
		\begin{equation}\label{eq:5-structural-constant-bounds}
			\begin{gathered}
				K\le C(1+M)^C,\qquad
				K4^{3C_0(1+\mathcal D)}
				\le C(1+M)^C e^{C(1+\mathcal D)},\qquad
				d_{\mathcal D}\le C(1+\mathcal D),\\
				A_{\mathcal D}^C100^{d_{\mathcal D}+1}
				(1+d_{\mathcal D})^C\le C e^{C(1+\mathcal D)},\qquad
				\varepsilon^{-1}\le C e^{C(1+\mathcal D)}.
			\end{gathered}
		\end{equation}
		In particular, by increasing the constant in the exponent of \(\delta\)
		and decreasing its prefactor, we obtain
		\begin{equation}\label{eq:5-delta-thresholds}
			\begin{aligned}
				8\delta&\le
				c\left(\frac{\varepsilon}{4000}\right)^{
					5(n+2)C_0(n,\lambda,\Lambda)(1+\mathcal D)}
				\bigl(K4^{3C_0(1+\mathcal D)}\bigr)^{-(n+2)/2},\\
				8\delta&\le
				\frac{c\varepsilon}
				{K A_{\mathcal D}^C100^{d_{\mathcal D}+1}
					(1+d_{\mathcal D})^C},\\
				8\delta&\le c\tau A_{\mathcal D}^{-C}\varepsilon^{1/2}.
			\end{aligned}
		\end{equation}
		The first line is the smallness required for definite drop and for
		almost monotonicity with error \(\varepsilon/100\).
		The second gives the common harmonic comparison and all relative
		normalizations, including the dyadic polynomial estimate.
		The last gives the changes of metric and of planes.
		Together with \eqref{eq:5-structural-constant-bounds}, these
		inequalities show that the choice
		\(\delta=c(1+M)^{-C}e^{-C(1+\mathcal D)^2}\) in
		\eqref{eq:5-final-delta} satisfies all three thresholds, with constants
		independent of the number of scales or covering generations.
		
		All smallness expressions defining \(r_2,r_7\), and
		\(r_6(\varepsilon)\) are nondecreasing in the trial radius. It is
		therefore enough to check them at \(4\cdot10^8C_3r_8\) and, for the
		occurrence of \(\eta\) in \(r_6\), at the enlarged argument
		\(C(1+d_{\mathcal D})4\cdot10^8C_3r_8\). The first line of
		\eqref{eq:5-delta-thresholds} gives the comparison, almost-monotonicity,
		and definite-drop requirements. The second gives the enlarged
		\(\eta\)-condition in \(r_6\), and the third gives the explicit
		coefficient-modulus condition in the second line of
		\eqref{eq:5-rmaster}.
		
		Every radius and every argument of the modulus in these
		applications is bounded by
		\[
		C(n,\lambda,\Lambda)(1+d_{\mathcal D})^2r_8.
		\]
		Indeed, the largest radius required directly by
		\eqref{eq:5-rmaster} is
		\[
		4\cdot10^8C_3r_8
		\le C(1+d_{\mathcal D})r_8.
		\]
		The definition of \(r_6\) enlarges this once more by
		\(2^{12}(1+\lambda^{-1/2})(1+d_{\mathcal D})\). Thus its largest
		argument in \(\eta\), and hence in the Dini integral, is at most
		\(C(1+d_{\mathcal D})^2r_8\). The direct harmonic
		comparison uses the same bound, including its factor \(8\sqrt\Lambda\).
		Since \(d_{\mathcal D}\le C(1+\mathcal D)\),
		\eqref{eq:5-explicit-r8} makes every such argument, denoted by
		\(t\), satisfy
		\[
		t<\rho_\omega(\delta),\qquad
		t\le\left(\frac{\delta}{1+M}\right)^{1/(2-n/p)},\qquad t\le c.
		\]
		Here the constant \(c\) in the choice of \(r_8\) has been
		decreased only in terms of the structural data. By the definition of
		\(\rho_\omega\) and monotonicity,
		\[
		\omega(t)+\int_0^t\frac{\omega(s)}s\,ds\le\delta,
		\qquad Mt^{2-n/p}\le\delta
		\]
		at the largest required argument \(t\). They also lie below every
		fixed geometric radius in the preceding comparison results.
		Combining these inequalities with \eqref{eq:5-delta-thresholds}
		verifies \eqref{eq:5-rmaster} and all its stated consequences.
		
		By \eqref{eq:5-A5-exponential-D},
		\[
		\mathcal A^{d_{\mathcal D}}
		\le\exp\bigl(C(1+\mathcal D)^2\bigr).
		\]
		Insert this and \eqref{eq:5-explicit-r8} into
		\eqref{eq:5-main-volume-estimate}. Since
		\(\rho_\omega(\delta)\le1\) and \(\delta\le1\),
		\[
		r_8^{-2}
		\le C(1+\mathcal D)^4\rho_\omega(\delta)^{-2}
		\left(\frac{1+M}{\delta}\right)^{2/(2-n/p)}.
		\]
		Equation \eqref{eq:5-final-delta} gives
		\[
		\left(\frac{1+M}{\delta}\right)^{\!2/(2-n/p)}
		=c^{-2/(2-n/p)}(1+M)^{2(1+C)/(2-n/p)}
		\exp\!\left[\frac{2C}{2-n/p}(1+\mathcal D)^2\right].
		\]
		Since \((1+\mathcal D)^4\le C\exp[C(1+\mathcal D)^2]\), we obtain
		\[
		r_8^{-2}\mathcal A^{d_{\mathcal D}}
		\le C(1+M)^C\exp\bigl[C(1+\mathcal D)^2\bigr]
		\rho_\omega(\delta)^{-2},
		\]
		which inserted into \eqref{eq:5-main-volume-estimate} proves
		\eqref{eq:5-final-clean-D-M}.
		In particular,
		\[
		\limsup_{r\downarrow0}r^{-2}
		|\{z\in B_1:\operatorname{dist}(z,S(u))<r\}|
		\le C(1+M)^C\exp\bigl[C(1+\mathcal D)^2\bigr]
		\rho_\omega(\delta)^{-2}<\infty.
		\]
		Thus the upper Minkowski content is finite and
		\(\dim_{\mathrm{Mink}}S(u)\le n-2\).
	\end{proof}
	
	\begin{proof}[Proof of Corollary \ref{cor:main-holder}]
		Under \eqref{eq:5-holder-modulus-corollary}, we may apply
		Theorem \ref{thm:main-volume-estimate} with the fixed modulus
		\(t^\alpha\); its constants then depend only on
		\(n,p,\lambda,\Lambda,\alpha\). For this modulus,
		\[
		\Omega_{t^\alpha}(s)
		=\left(1+\frac1\alpha\right)s^\alpha,
		\qquad
		\rho_{t^\alpha}(\delta)
		=\left(\frac{\alpha\delta}{1+\alpha}\right)^{1/\alpha}.
		\]
		By \eqref{eq:5-final-delta},
		\[
		\rho_{t^\alpha}(\delta)^{-2}
		\le C(1+M)^C
		\exp\bigl[C(1+\mathcal D)^2\bigr].
		\]
		Substitution into \eqref{eq:5-final-clean-D-M} proves
		\eqref{eq:5-final-holder-clean}.
	\end{proof}
	
	\medskip
	
	\section{Boundary estimates}\label{sec:doubling-verification}
	
	Throughout this section, $\Omega$ is a bounded connected
	$C^{1,\mathrm{Dini}}$ domain, $p>n$, and
	\begin{equation}\label{eq:boundary-schrodinger-equation}
		\Delta u+Vu=0\quad\text{in }\Omega,
		\qquad u\in W^{1,2}_0(\Omega),\quad u\not\equiv0,
	\end{equation}
	which naturally extends trivially to a solution in whole space $\mathbb{R}^n$.
	We put
	\begin{equation}\label{eq:boundary-potential-bound}
		M=\|V\|_{L^p(\Omega)}.
	\end{equation}
	The constants in this section depend only on $n,p$ and the quantitative
	$C^{1,\mathrm{Dini}}$ character of $\Omega$, unless otherwise stated.
	These quantitative domain data include the radius and number of boundary
	charts, their uniform graph bounds and Dini modulus, and a diameter bound.
	In boundary charts, the modulus $\omega_{\partial\Omega}$ controls the
	gradient of the defining function, after multiplication by a geometric
	constant, and satisfies
	\begin{equation}\label{eq:boundary-dini-modulus}
		\int_0^{r_\Omega}\frac{\omega_{\partial\Omega}(t)}t\,dt<\infty.
	\end{equation}
	Define
	\begin{equation}\label{eq:schrodinger-ball-mass}
		I_u(y,r)=\int_{\Omega\cap B_r(y)}u^2,
		\qquad y\in\overline\Omega,
	\end{equation}
	and
	\begin{equation}\label{eq:schrodinger-ball-doubling}
		D_\Omega^u(y,r)=
		\log_4\frac{\vint_{\Omega\cap B_{2r}(y)}u^2}
		{\vint_{\Omega\cap B_r(y)}u^2}.
	\end{equation}
	The masses are positive by unique continuation for $\Delta+V$;
	see \cite{JerisonKenig1985}.
	
	\subsection{Weighted frequency and three-ball inequalities}
	
	Fix a center $y$ for which $\Omega\cap B_R(y)$ is star-shaped with
	respect to $y$. For $\alpha\ge2$, write
	\begin{equation}\label{eq:schrodinger-weight}
		\psi_{y,r}(x)=r^2-|x-y|^2
	\end{equation}
	and set
	\begin{align}
		\mathscr H_y(r)&=\int_{\Omega\cap B_r(y)}u^2\psi_{y,r}^{\alpha-1},
		\label{eq:weighted-H-schrodinger}\\
		\mathscr E_y(r)&=\int_{\Omega\cap B_r(y)}|\nabla u|^2\psi_{y,r}^{\alpha},
		\label{eq:weighted-E-schrodinger}\\
		\mathscr P_y(r)&=-\int_{\Omega\cap B_r(y)}Vu^2\psi_{y,r}^{\alpha},
		\label{eq:weighted-P-schrodinger}\\
		\mathscr N_y(r)&=\frac{\mathscr E_y(r)}{\mathscr H_y(r)}.
		\label{eq:weighted-frequency-schrodinger}
	\end{align}
	\begin{lemma}\label{lem:weighted-identities-schrodinger}
		For almost every $0<r<R$,
		\begin{equation}\label{eq:H-prime-schrodinger}
			\frac{\mathscr H_y'(r)}{\mathscr H_y(r)}
			=\frac{2\alpha+n-2}{r}
			+\frac{\mathscr N_y(r)}{\alpha r}
			+\frac{\mathscr P_y(r)}{\alpha r\mathscr H_y(r)},
		\end{equation}
		and
		\begin{equation}\label{eq:N-prime-preliminary}
			\mathscr N_y'(r)\ge
			-\frac{1}{4\alpha r\mathscr H_y(r)}
			\int_{\Omega\cap B_r(y)}|Vu|^2\psi_{y,r}^{\alpha+1}.
		\end{equation}
	\end{lemma}
	
	\begin{proof}
		The full-ball calculation is that of
		\cite[Lemmas~3.1 and~3.2]{Davey2026Frequency}; we record the boundary
		terms. Testing the equation with $u\psi_{y,r}^{\alpha}$ gives
		\[
		\mathscr E_y+\mathscr P_y
		=2\alpha\int_{\Omega\cap B_r(y)}
		u(x-y)\cdot\nabla u\,\psi_{y,r}^{\alpha-1}.
		\]
		Since the weight vanishes on the spherical boundary, differentiation
		under the integral sign gives
		\[
		\mathscr H_y'(r)=2(\alpha-1)r
		\int_{\Omega\cap B_r(y)}u^2\psi_{y,r}^{\alpha-2}.
		\]
		The boundary integral of
		$u^2(x-y)\psi_{y,r}^{\alpha-1}$ vanishes on both parts of the
		boundary. Hence
		\begin{align*}
			0={}&2\int_{\Omega\cap B_r(y)}u(x-y)\cdot\nabla u\,
			\psi_{y,r}^{\alpha-1}+n\mathscr H_y(r)\\
			&-2(\alpha-1)\int_{\Omega\cap B_r(y)}
			u^2|x-y|^2\psi_{y,r}^{\alpha-2}.
		\end{align*}
		Using $|x-y|^2=r^2-\psi_{y,r}$ in the last integral and the preceding
		testing identity yields
		\[
		r\mathscr H_y'(r)
		=(2\alpha+n-2)\mathscr H_y(r)
		+\frac{\mathscr E_y(r)+\mathscr P_y(r)}{\alpha},
		\]
		which is \eqref{eq:H-prime-schrodinger}.
		
		We first justify the Rellich calculation at the regularity of the
		lemma. If $n>2$, then $p>n$ implies
		\[
		\frac{2p}{p-2}\le\frac{2n}{n-2},
		\qquad
		\|Vu\|_{L^2}\le\|V\|_{L^p}
		\|u\|_{L^{2p/(p-2)}}<\infty;
		\]
		for $n=2$ the same conclusion follows from
		$W^{1,2}\hookrightarrow L^q$ for every finite $q$. Thus
		$\Delta u=-Vu\in L^2$. For a smooth Dirichlet solution, apply the
		divergence theorem to the vector field
		\[
		|\nabla u|^2(x-y)\psi_{y,r}^{\alpha}
		-2((x-y)\cdot\nabla u)\nabla u\,\psi_{y,r}^{\alpha}.
		\]
		Expanding its divergence and using
		$\mathscr E_y'(r)=2\alpha r\int|\nabla u|^2
		\psi_{y,r}^{\alpha-1}$ gives, after discarding the nonnegative
		boundary term on \(\partial\Omega\cap B_r(y)\),
		\begin{align}\label{eq:Rellich-weighted-schrodinger}
			\mathscr E_y'(r)
			\ge{}&\frac{2\alpha+n-2}{r}\mathscr E_y(r)
			+\frac{4\alpha}{r}\int_{\Omega\cap B_r(y)}
			((x-y)\cdot\nabla u)^2\psi_{y,r}^{\alpha-1}\notag\\
			&-\frac2r\int_{\Omega\cap B_r(y)}
			((x-y)\cdot\nabla u)\Delta u\,\psi_{y,r}^{\alpha}.
		\end{align}
		Indeed, on \(\partial\Omega\cap B_r(y)\) the discarded term is
		$r^{-1}\int ((x-y)\cdot\nu)(\partial_\nu u)^2
		\psi_{y,r}^{\alpha}\ge0$, while the weight vanishes on the spherical
		boundary.
		
		Differentiate $\mathscr E_y/\mathscr H_y$, discard the nonnegative
		boundary term, and complete the square. Then, we obtain
		\begin{align*}
			\mathscr N_y'(r)\ge{}&
			\frac{4\alpha}{r\mathscr H_y}
			\int_{\Omega\cap B_r(y)}
			\left((x-y)\cdot\nabla u-\frac{\psi_{y,r}\Delta u}{4\alpha}\right)^2
			\psi_{y,r}^{\alpha-1}\\
			&-\frac{\mathscr E_y(\mathscr E_y+\mathscr P_y)}
			{\alpha r\mathscr H_y^2}
			-\frac1{4\alpha r\mathscr H_y}
			\int_{\Omega\cap B_r(y)}|\Delta u|^2\psi_{y,r}^{\alpha+1}.
		\end{align*}
		Since
		\[
		\int_{\Omega\cap B_r(y)}u
		\left((x-y)\cdot\nabla u-\frac{\psi_{y,r}\Delta u}{4\alpha}\right)
		\psi_{y,r}^{\alpha-1}
		=\frac{2\mathscr E_y+\mathscr P_y}{4\alpha},
		\]
		Cauchy--Schwarz bounds the first two terms from below by
		\[
		\frac{(2\mathscr E_y+\mathscr P_y)^2
			-4\mathscr E_y(\mathscr E_y+\mathscr P_y)}
		{4\alpha r\mathscr H_y^2}
		=\frac{\mathscr P_y^2}{4\alpha r\mathscr H_y^2}\ge0.
		\]
		Using $\Delta u=-Vu$ proves \eqref{eq:N-prime-preliminary}.
	\end{proof}
	
	\begin{lemma}\label{lem:weighted-potential-schrodinger}
		For $\alpha\ge2$, the following estimates hold with $C=C(n,p)$:
		\begin{equation}\label{eq:P-estimate-schrodinger}
			|\mathscr P_y(r)|\le\tfrac12\mathscr E_y(r)
			+C(\alpha+M^{\frac{2p}{2p-n}}r^2)\mathscr H_y(r),
		\end{equation}
		and
		\begin{equation}\label{eq:V2-estimate-schrodinger}
			\int_{\Omega\cap B_r(y)}|Vu|^2\psi_{y,r}^{\alpha+1}
			\le C M^{\frac{2p}{3p-2n}}r^{\frac{2(2p-n)}{3p-2n}}
			\bigl(M^{\frac{4p}{3p-2n}}r^{\frac{4(2p-n)}{3p-2n}}
			+\alpha+\mathscr N_y(r)\bigr)\mathscr H_y(r).
		\end{equation}
	\end{lemma}
	
	\begin{proof}
		We adapt the constant-coefficient calculation in
		\cite[Lemma~4.3]{Davey2026Frequency}. The zero extension of
		$u\psi_{y,r}^{\alpha/2}$ belongs to $W^{1,2}(\mathbb R^n)$.
		H\"older's inequality and Gagliardo--Nirenberg interpolation at
		exponent $2p/(p-1)$ give, also when $n=2$,
		\begin{equation}\label{eq:V-form-interpolation}
			|\mathscr P_y(r)|
			\le CMr^{2-n/p}\mathscr H_y(r)^{1-n/(2p)}
			\left(\int_{\Omega\cap B_r(y)}
			|\nabla(u\psi_{y,r}^{\alpha/2})|^2\right)^{n/(2p)}.
		\end{equation}
		Expansion of the gradient and integration by parts give the exact identity
		\[
		\int_{\Omega\cap B_r(y)}|\nabla(u\psi_{y,r}^{\alpha/2})|^2
		=\frac{\alpha}{2(\alpha-1)}\mathscr E_y
		-\frac{\alpha-2}{2(\alpha-1)}\mathscr P_y
		+\frac{n\alpha^2}{2(\alpha-1)}\mathscr H_y.
		\]
		In particular,
		\begin{equation}\label{eq:weighted-gradient-first-bound}
			\int_{\Omega\cap B_r(y)}|\nabla(u\psi_{y,r}^{\alpha/2})|^2
			\le\mathscr E_y+\tfrac12|\mathscr P_y|+n\alpha\mathscr H_y.
		\end{equation}
		For $0<\zeta\le1/4$, Young's inequality therefore yields
		\[
		|\mathscr P_y|
		\le\frac\zeta2\left(\mathscr E_y+\tfrac12|\mathscr P_y|
		+n\alpha\mathscr H_y\right)
		+C\zeta^{-n/(2p-n)}M^{\frac{2p}{2p-n}}r^2\mathscr H_y.
		\]
		Absorbing the term containing $|\mathscr P_y|$ proves
		\begin{equation}\label{eq:refined-P-estimate-recentering}
			|\mathscr P_y(r)|\le\zeta\mathscr E_y(r)
			+C\left(\zeta\alpha+
			\zeta^{-\frac{n}{2p-n}}M^{\frac{2p}{2p-n}}r^2\right)
			\mathscr H_y(r).
		\end{equation}
		Taking a fixed $\zeta$ proves \eqref{eq:P-estimate-schrodinger}.
		The factor $\zeta\alpha$ in the refined estimate will be important when
		the center is changed below.
		
		For the second estimate, extend
		$u\psi_{y,r}^{(\alpha+1)/2}$ by zero. Formula
		\eqref{eq:weighted-gradient-first-bound}, with
		$\alpha+1$ in place of $\alpha$, and
		\eqref{eq:P-estimate-schrodinger} imply
		\[
		\|u\psi_{y,r}^{(\alpha+1)/2}\|_2^2\le r^4\mathscr H_y,
		\qquad
		\|\nabla(u\psi_{y,r}^{(\alpha+1)/2})\|_2^2
		\le Cr^2\bigl(\mathscr E_y+(\alpha+M^{\frac{2p}{2p-n}}r^2)\mathscr H_y\bigr).
		\]
		H\"older's inequality and Gagliardo--Nirenberg interpolation at exponent
		$2p/(p-2)$ now give
		\begin{align*}
			\int_{\Omega\cap B_r(y)}|Vu|^2\psi_{y,r}^{\alpha+1}
			&\le CM^2\|u\psi_{y,r}^{(\alpha+1)/2}\|_2^{2(1-n/p)}
			\|\nabla(u\psi_{y,r}^{(\alpha+1)/2})\|_2^{2n/p}\\
			&\le CM^2r^{4-2n/p}\mathscr H_y
			\bigl(\mathscr N_y+\alpha+M^{\frac{2p}{2p-n}}r^2\bigr)^{n/p}.
		\end{align*}
		Young's inequality gives
		\[
		\begin{aligned}
			&M^2r^{4-2n/p}
			\bigl(\mathscr N_y+\alpha+M^{\frac{2p}{2p-n}}r^2\bigr)^{n/p}\\
			&\qquad\le C
			M^{\frac{2p}{3p-2n}}r^{\frac{2(2p-n)}{3p-2n}}
			\left(M^{\frac{4p}{3p-2n}}r^{\frac{4(2p-n)}{3p-2n}}
			+\mathscr N_y+\alpha+M^{\frac{2p}{2p-n}}r^2\right).
		\end{aligned}
		\]
		Since
		\[
		M^{\frac{2p}{2p-n}}r^2
		\le1+M^{\frac{4p}{3p-2n}}r^{\frac{4(2p-n)}{3p-2n}}
		\]
		and $\alpha\ge2$, this proves
		\eqref{eq:V2-estimate-schrodinger}.
	\end{proof}
	
	\begin{proposition}\label{prop:modified-frequency-schrodinger}
		There are constants $C_1,C_2$, depending only on $n,p$, such that
		\begin{equation}\label{eq:modified-frequency-schrodinger}
			\mathscr F_y(r)=
			\left[\mathscr N_y(r)+\frac{C_2}{\alpha}
			M^{\frac{2p}{3p-2n}}r^{\frac{2(2p-n)}{3p-2n}}
			\bigl(M^{\frac{4p}{3p-2n}}r^{\frac{4(2p-n)}{3p-2n}}+\alpha\bigr)\right]
			\exp\left(\frac{C_1}{\alpha}
			M^{\frac{2p}{3p-2n}}r^{\frac{2(2p-n)}{3p-2n}}\right)
		\end{equation}
		is nondecreasing on $(0,R)$.
	\end{proposition}
	
	\begin{proof}
		The preceding lemmas give
		\begin{equation}\label{eq:modified-frequency-differential}
			\mathscr N_y'(r)\ge
			-\frac{C}{\alpha r}M^{\frac{2p}{3p-2n}}r^{\frac{2(2p-n)}{3p-2n}}
			\bigl(M^{\frac{4p}{3p-2n}}r^{\frac{4(2p-n)}{3p-2n}}
			+\alpha+\mathscr N_y(r)\bigr).
		\end{equation}
		Direct differentiation gives
		\begin{align}
			\frac{d}{dr}M^{\frac{2p}{3p-2n}}r^{\frac{2(2p-n)}{3p-2n}}
			&=\frac{2(2p-n)}{(3p-2n)r}
			M^{\frac{2p}{3p-2n}}r^{\frac{2(2p-n)}{3p-2n}},
			\label{eq:modified-frequency-power-derivative}\\
			\frac{d}{dr}\left[
			M^{\frac{2p}{3p-2n}}r^{\frac{2(2p-n)}{3p-2n}}
			\left(M^{\frac{4p}{3p-2n}}r^{\frac{4(2p-n)}{3p-2n}}+\alpha\right)\right]
			&=\frac{2(2p-n)\alpha}{(3p-2n)r}
			M^{\frac{2p}{3p-2n}}r^{\frac{2(2p-n)}{3p-2n}}\notag\\
			&\quad+\frac{6(2p-n)}{(3p-2n)r}
			M^{\frac{6p}{3p-2n}}r^{\frac{6(2p-n)}{3p-2n}}\notag\\
			&\ge\frac{c}{r}M^{\frac{2p}{3p-2n}}r^{\frac{2(2p-n)}{3p-2n}}
			\left(\alpha+M^{\frac{4p}{3p-2n}}r^{\frac{4(2p-n)}{3p-2n}}\right).
			\label{eq:modified-frequency-correction-derivative}
		\end{align}
		Since both terms in the brackets in
		\eqref{eq:modified-frequency-schrodinger} are nonnegative,
		combining
		\eqref{eq:modified-frequency-differential} and
		\eqref{eq:modified-frequency-correction-derivative}, and then choosing
		$C_2\ge2C/c$, gives
		\[
		\begin{aligned}
			&\frac{d}{dr}\left[\mathscr N_y(r)+\frac{C_2}{\alpha}
			M^{\frac{2p}{3p-2n}}r^{\frac{2(2p-n)}{3p-2n}}
			\left(M^{\frac{4p}{3p-2n}}r^{\frac{4(2p-n)}{3p-2n}}+\alpha\right)\right]\\
			&\quad\ge-\frac{C}{\alpha r}M^{\frac{2p}{3p-2n}}r^{\frac{2(2p-n)}{3p-2n}}
			\left[\mathscr N_y(r)+\frac{C_2}{\alpha}
			M^{\frac{2p}{3p-2n}}r^{\frac{2(2p-n)}{3p-2n}}
			\left(M^{\frac{4p}{3p-2n}}r^{\frac{4(2p-n)}{3p-2n}}+\alpha\right)\right]\\
			&\qquad+\frac{cC_2-C}{\alpha r}M^{\frac{2p}{3p-2n}}r^{\frac{2(2p-n)}{3p-2n}}
			\left(\alpha+M^{\frac{4p}{3p-2n}}r^{\frac{4(2p-n)}{3p-2n}}\right)\\
			&\quad\ge-\frac{C}{\alpha r}M^{\frac{2p}{3p-2n}}r^{\frac{2(2p-n)}{3p-2n}}
			\left[\mathscr N_y(r)+\frac{C_2}{\alpha}
			M^{\frac{2p}{3p-2n}}r^{\frac{2(2p-n)}{3p-2n}}
			\left(M^{\frac{4p}{3p-2n}}r^{\frac{4(2p-n)}{3p-2n}}+\alpha\right)\right].
		\end{aligned}
		\]
		By \eqref{eq:modified-frequency-power-derivative},
		\begin{align*}
			\mathscr F_y'(r)
			&\ge\frac{\frac{2(2p-n)}{3p-2n}C_1-C}{\alpha r}
			M^{\frac{2p}{3p-2n}}r^{\frac{2(2p-n)}{3p-2n}}\\
			&\quad\times\left[\mathscr N_y(r)+\frac{C_2}{\alpha}
			M^{\frac{2p}{3p-2n}}r^{\frac{2(2p-n)}{3p-2n}}
			\left(M^{\frac{4p}{3p-2n}}r^{\frac{4(2p-n)}{3p-2n}}+\alpha\right)\right]
			\exp\left(\frac{C_1}{\alpha}
			M^{\frac{2p}{3p-2n}}r^{\frac{2(2p-n)}{3p-2n}}\right)\ge0
		\end{align*}
		once $C_1\ge C(3p-2n)/(2(2p-n))$. This proves the monotonicity.
	\end{proof}
	
	For the rest of this section, choose
	\begin{equation}\label{eq:alpha-star-schrodinger}
		\alpha_1 =2+M^{\frac{2p}{3p-2n}}.
	\end{equation}
	The reason we do this choice is such $\alpha_1$ gives sharp exponent in three balls inequality and thus in vanishing order, see for example \cite{Davey2026Frequency}.
	All weighted integrals below use $\alpha=\alpha_1$. Next is the three balls inequality on star-shaped domain, which plays key role in later doubling index estimate in $C^{1,Dini}$ domain.
	
	\begin{lemma}\label{lem:fixed-ratio-three-ball-schrodinger}
		For fixed $0<a_1<a_2<a_3<1$, there are $C>0$ and $0<\vartheta<1$,
		depending only on $n,p,a_1,a_2,a_3$, such that
		\begin{equation}\label{eq:fixed-ratio-three-ball-schrodinger}
			I_u(y,a_2r)\le e^{C\alpha_1}
			I_u(y,a_1r)^\vartheta I_u(y,a_3r)^{1-\vartheta}
		\end{equation}
		whenever $0<r\le1$ and $\Omega\cap B_r(y)$ is star-shaped with
		respect to $y$. The constants are uniform when the radius ratios range
		over a compact subset of $0<a_1<a_2<a_3<1$.
	\end{lemma}
	
	\begin{proof}
		Let \(C_1,C_2\) be the constants supplied by Proposition
		\ref{prop:modified-frequency-schrodinger}.
		For $0<t\le1$, the choice of $\alpha_1$ implies
		$M^{\frac{2p}{3p-2n}}t^{\frac{2(2p-n)}{3p-2n}}\le\alpha_1$. Consequently,
		\[
		e^{-C_1}\mathscr F_y(t)-C\alpha_1^2
		\le\mathscr N_y(t)\le\mathscr F_y(t).
		\]
		Equations \eqref{eq:H-prime-schrodinger} and
		\eqref{eq:P-estimate-schrodinger} give
		\begin{equation}\label{eq:three-ball-H-comparison}
			\frac{2\alpha_1+n-2}{t}
			+\frac{c\mathscr F_y(t)}{\alpha_1t}-\frac{C\alpha_1}{t}
			\le\frac{\mathscr H_y'(t)}{\mathscr H_y(t)}
			\le\frac{2\alpha_1+n-2}{t}
			+\frac{C\mathscr F_y(t)}{\alpha_1t}+\frac{C\alpha_1}{t}.
		\end{equation}
		Set
		\[
		t_0=\frac{a_2+a_3}{2}r,\qquad
		\ell_1=\log\frac{t_0}{a_1r},\qquad
		\ell_2=\log\frac{a_3r}{t_0}.
		\]
		Integrating the upper bound in \eqref{eq:three-ball-H-comparison} over
		$[a_1r,t_0]$ and using
		$\mathscr F_y(t)\le\mathscr F_y(t_0)$ gives
		\begin{equation}\label{eq:three-ball-upper-integrated}
			\log\frac{\mathscr H_y(t_0)}{\mathscr H_y(a_1r)}
			\le(2\alpha_1+n-2)\ell_1
			+\frac{C\ell_1}{\alpha_1}\mathscr F_y(t_0)
			+C\alpha_1\ell_1.
		\end{equation}
		Similarly, the lower bound on $[t_0,a_3r]$ and
		$\mathscr F_y(t)\ge\mathscr F_y(t_0)$ give
		\begin{equation}\label{eq:three-ball-lower-integrated}
			\log\frac{\mathscr H_y(a_3r)}{\mathscr H_y(t_0)}
			\ge(2\alpha_1+n-2)\ell_2
			+\frac{c\ell_2}{\alpha_1}\mathscr F_y(t_0)
			-C\alpha_1\ell_2.
		\end{equation}
		Solving \eqref{eq:three-ball-lower-integrated} for
		$\mathscr F_y(t_0)/\alpha_1$ and substituting into
		\eqref{eq:three-ball-upper-integrated} yields
		\[
		\log\frac{\mathscr H_y(t_0)}{\mathscr H_y(a_1r)}
		\le\frac{C\ell_1}{c\ell_2}
		\log\frac{\mathscr H_y(a_3r)}{\mathscr H_y(t_0)}
		+C\alpha_1.
		\]
		With $\vartheta=(1+C\ell_1/(c\ell_2))^{-1}$, this is
		\[
		\mathscr H_y(t_0)
		\le e^{C\alpha_1}
		\mathscr H_y(a_1r)^\vartheta
		\mathscr H_y(a_3r)^{1-\vartheta}.
		\]
		Finally,
		\[
		\begin{split}
			\mathscr H_y(t_0)
			&\ge\bigl(((a_2+a_3)^2/4-a_2^2)r^2\bigr)^{\alpha_1-1}
			I_u(y,a_2r),\\
			\mathscr H_y(a_ir)&\le(a_i^2r^2)^{\alpha_1-1}I_u(y,a_ir),
			\qquad i=1,3.
		\end{split}
		\]
		After substitution, the powers of $r$ cancel and all fixed radius
		factors are bounded by $e^{C\alpha_1}$. This proves
		\eqref{eq:fixed-ratio-three-ball-schrodinger}.
		The logarithmic interval lengths and the displayed positive radius
		differences are uniformly controlled on compact sets of radius ratios.
	\end{proof}
	
	\subsection{Doubling at the boundary}
	
	In this subsection, we prove the global uniform upper bound of doubling index. This is achieved by the summable radius series acting by Dini function, together with the propagation of smallness argument.
	
	The next lemma gives the star-shaped condition needed by previous subsection.
	
	\begin{lemma}\label{lem:dini-boundary-star-shaped}
		There is a geometric radius $r_9>0$ such that, for $y\in\overline\Omega$
		and $0<r\le r_9$, one can choose $y_r\in\overline\Omega$ with
		\begin{equation}\label{eq:center-shift-size}
			|y_r-y|\le Cr\omega_{\partial\Omega}(Cr)\le r/32
		\end{equation}
		so that $\Omega\cap B_{8r}(y_r)$ is star-shaped with respect to $y_r$.
	\end{lemma}
	
	\begin{proof}
		If $B_{8r}(y)\subset\Omega$, take $y_r=y$. Otherwise, choose
		$\xi\in\partial\Omega$ with $|\xi-y|\le8r$. Translate and rotate a
		boundary chart about $\xi$ so that $\xi=0$, the inward normal is $e_n$,
		and
		\[
		\Omega=\{x_n>\varphi(x')\},\qquad
		\varphi(0)=0,\qquad \nabla\varphi(0)=0.
		\]
		In these coordinates put
		$\tau=\sup_{|x'|\le32r}|\nabla\varphi(x')|$.
		After decreasing $r_9$, $\tau\le C\omega_{\partial\Omega}(Cr)\le1/512$.
		Take $y_r=y+16\tau r e_n$. For $z\in\Omega\cap B_{8r}(y_r)$ and
		$0<t<1$, put $w_t=ty_r+(1-t)z$. The Lipschitz bound gives
		\[
		\varphi(ty'+(1-t)z')
		-t\varphi(y')-(1-t)\varphi(z')
		\le2t(1-t)\tau|y'-z'|\le16t(1-t)\tau r.
		\]
		Since $y_n\ge\varphi(y')$ and $z_n>\varphi(z')$, it follows that
		\[
		\begin{aligned}
			(w_t)_n-\varphi(w_t')
			&=t(y_n+16\tau r)+(1-t)z_n
			-\varphi(ty'+(1-t)z')\\
			&\ge16t\tau r-16t(1-t)\tau r
			=16t^2\tau r\ge0.
		\end{aligned}
		\]
		Thus the segment from $y_r$ to $z$ lies above the graph. Convexity of
		the ball proves the assertion. Allowing $y_r$ on the boundary includes
		the flat case $\tau=0$.
	\end{proof}
	
	Now we use the Dini condition and three balls inequality Lemma \ref{lem:fixed-ratio-three-ball-schrodinger} to obtain doubling index estimate at a fixed point $y$.
	
	\begin{lemma}\label{lem:dyadic-dini-recentering}
		There is $r_{10}>0$, depending only on $n,p,\Omega$, such that
		\begin{equation}\label{eq:dyadic-dini-recentering-conclusion}
			\sup_{0<s\le r}\log\frac{I_u(y,s)}{I_u(y,s/2)}
			\le C\left(\alpha_1+
			\log\frac{I_u(y,4r)}{I_u(y,2r)}\right),
			\qquad y\in\overline\Omega,\quad 0<r\le r_{10}.
		\end{equation}
	\end{lemma}
	
	\begin{proof}
		Define
		\[
		\tau_2=\min\left\{\frac{2(2p-n)}{3p-2n},2-\frac np\right\},
		\qquad
		\varepsilon(t)=C\omega_{\partial\Omega}(Ct)+Ct^{\tau_2}.
		\]
		Choose
		$r_{10}\le\min\{r_9/2,1/8\}$, where $r_9$ is supplied by Lemma
		\ref{lem:dini-boundary-star-shaped}, so that
		$C\varepsilon(t)\le1/8$ for $0<t\le2r_{10}$.
		This choice is independent of $M$.
		For $0<t\le2r_{10}$, Lemma~\ref{lem:dini-boundary-star-shaped}
		supplies $y_t$.
		
		Apply \eqref{eq:refined-P-estimate-recentering} with
		$\zeta=t^{2-n/p}/4$. For $t/8\le s\le t$, the elementary bounds
		\[
		\begin{split}
			\alpha_1^{-1}M^{\frac{2p}{3p-2n}}s^{\frac{2(2p-n)}{3p-2n}}
			&\le t^{\frac{2(2p-n)}{3p-2n}},\\
			\alpha_1^{-2}M^{\frac{2p}{3p-2n}}s^{\frac{2(2p-n)}{3p-2n}}
			\bigl(M^{\frac{4p}{3p-2n}}s^{\frac{4(2p-n)}{3p-2n}}+\alpha_1\bigr)
			&\le C\alpha_1\varepsilon(t),\\
			\zeta+\alpha_1^{-1}\zeta^{-\frac n{2p-n}}
			M^{\frac{2p}{2p-n}}s^2
			&\le C\alpha_1\varepsilon(t)
		\end{split}
		\]
		follow from \eqref{eq:alpha-star-schrodinger}. To see that the bounds
		are uniform also when $M<1$, split into two cases. If $M\le1$, every
		positive power of $M$ is at most one and $\alpha_1\ge2$. If $M\ge1$,
		\begin{align*}
			\frac{2p}{2p-n}\frac{2p-n}{3p-2n}
			&=\frac{2p}{3p-2n},\qquad
			\frac{2p}{2p-n}-\frac{2p}{3p-2n}
			\le\frac{2p}{3p-2n},\\
			\left(2-\frac np\right)\frac n{2p-n}&=\frac np,
			\qquad M^{\frac{2p}{3p-2n}}\le\alpha_1.
		\end{align*}
		Thus
		$\alpha_1^{-1}M^{\frac{2p}{3p-2n}}s^{\frac{2(2p-n)}{3p-2n}}
		\le t^{\frac{2(2p-n)}{3p-2n}}$ and
		\[
		\alpha_1^{-1}\zeta^{-n/(2p-n)}M^{\frac{2p}{2p-n}}s^2
		\le C\alpha_1t^{2-n/p}.
		\]
		The middle bound follows from
		$M^{\frac{2p}{3p-2n}}s^{\frac{2(2p-n)}{3p-2n}}
		\le\alpha_1t^{\frac{2(2p-n)}{3p-2n}}$ and
		$M^{\frac{4p}{3p-2n}}s^{\frac{4(2p-n)}{3p-2n}}
		\le\alpha_1^2t^{\frac{4(2p-n)}{3p-2n}}$. Since
		$\tau_2\le\frac{2(2p-n)}{3p-2n}$ and
		$\tau_2\le2-n/p$, all three right-hand sides are
		bounded by the displayed multiples of $\varepsilon(t)$.
		The first bound controls the exponential in
		\eqref{eq:modified-frequency-schrodinger}, and the second controls
		its additive term after division by $\alpha_1$. More explicitly,
		\begin{align*}
			\mathscr N_{y_t}(s)
			={}&\exp\left(-\frac{C_1}{\alpha_1}
			M^{\frac{2p}{3p-2n}}s^{\frac{2(2p-n)}{3p-2n}}\right)\mathscr F_{y_t}(s)\\
			&-\frac{C_2}{\alpha_1}M^{\frac{2p}{3p-2n}}s^{\frac{2(2p-n)}{3p-2n}}
			\left(M^{\frac{4p}{3p-2n}}s^{\frac{4(2p-n)}{3p-2n}}+\alpha_1\right).
		\end{align*}
		Combining these bounds with
		\eqref{eq:H-prime-schrodinger} and
		\eqref{eq:modified-frequency-schrodinger} gives
		\begin{equation}\label{eq:recentering-H-differential}
			\begin{split}
				\frac{\mathscr H_{y_t}'(s)}{\mathscr H_{y_t}(s)}
				&\le\frac{2\alpha_1+n-2}{s}
				+\frac{1+C\varepsilon(t)}{\alpha_1s}\mathscr F_{y_t}(s)
				+\frac{C\alpha_1\varepsilon(t)}s,\\
				\frac{\mathscr H_{y_t}'(s)}{\mathscr H_{y_t}(s)}
				&\ge\frac{2\alpha_1+n-2}{s}
				+\frac{1-C\varepsilon(t)}{\alpha_1s}\mathscr F_{y_t}(s)
				-\frac{C\alpha_1\varepsilon(t)}s.
			\end{split}
		\end{equation}
		
		The multiplicative error in this comparison is independent of $M$;
		the potential enters only through the factor $\alpha_1$ in the
		additive error. This distinction permits iteration over infinitely
		many scales.
		
		Set $\rho_1=t/4-|y_t-y|$, $\rho_2=t/2+|y_t-y|$, and
		$\rho_3=t-|y_t-y|$. Since $|y_t-y|/t\le1/32$, these radii belong
		to $[t/8,t]$. Put
		\[
		\ell_1=\log\frac{\rho_2}{\rho_1},\qquad
		\ell_2=\log\frac{\rho_3}{\rho_2}.
		\]
		Both are bounded above and below by positive absolute constants.
		The mean value theorem, applied to the logarithms on
		$0\le |y_t-y|/t\le1/32$, gives
		\begin{equation}\label{eq:log-radius-ratio-error}
			|\ell_1-\log2|+|\ell_2-\log2|
			\le C\frac{|y_t-y|}t,
			\qquad
			\left|\frac{\ell_1}{\ell_2}-1\right|
			\le C\frac{|y_t-y|}t\le C\varepsilon(t).
		\end{equation}
		Integrate the upper bound in \eqref{eq:recentering-H-differential}
		over $[\rho_1,\rho_2]$ and the lower bound over $[\rho_2,\rho_3]$.
		Monotonicity of $\mathscr F_{y_t}$ gives
		\begin{align*}
			\frac1{\ell_1}
			\log\frac{\mathscr H_{y_t}(\rho_2)}{\mathscr H_{y_t}(\rho_1)}
			&\le 2\alpha_1+n-2
			+\frac{1+C\varepsilon(t)}{\alpha_1}\mathscr F_{y_t}(\rho_2)
			+C\alpha_1\varepsilon(t),\\
			\frac1{\ell_2}
			\log\frac{\mathscr H_{y_t}(\rho_3)}{\mathscr H_{y_t}(\rho_2)}
			&\ge 2\alpha_1+n-2
			+\frac{1-C\varepsilon(t)}{\alpha_1}\mathscr F_{y_t}(\rho_2)
			-C\alpha_1\varepsilon(t).
		\end{align*}
		Eliminating $\mathscr F_{y_t}(\rho_2)$ and using
		\eqref{eq:log-radius-ratio-error} yield
		\[
		\log\frac{\mathscr H_{y_t}(\rho_2)}{\mathscr H_{y_t}(\rho_1)}
		\le(1+C\varepsilon(t))
		\log\frac{\mathscr H_{y_t}(\rho_3)}{\mathscr H_{y_t}(\rho_2)}
		+C\alpha_1\varepsilon(t).
		\]
		For every $R>|y_t-y|$, the triangle inequality and positivity of
		the weights give
		\[
		\mathscr H_y(R-|y_t-y|)
		\le\mathscr H_{y_t}(R)
		\le\mathscr H_y(R+|y_t-y|).
		\]
		For example,
		$R^2-(|x-y|+|y_t-y|)^2
		\ge(R-|y_t-y|)^2-|x-y|^2$
		on $B_{R-|y_t-y|}(y)$; the reverse comparison follows by interchanging
		the centers. At the three chosen radii these inequalities give
		\[
		\begin{gathered}
			\mathscr H_y(t/2)\le\mathscr H_{y_t}(\rho_2),
			\qquad \mathscr H_{y_t}(\rho_1)\le\mathscr H_y(t/4),\\
			\mathscr H_{y_t}(\rho_3)\le\mathscr H_y(t),
			\qquad \mathscr H_y(t/2)\le\mathscr H_{y_t}(\rho_2).
		\end{gathered}
		\]
		Substitution proves
		\begin{equation}\label{eq:adjacent-increment-comparison}
			\log\frac{\mathscr H_y(t/2)}{\mathscr H_y(t/4)}
			\le(1+C\varepsilon(t))
			\log\frac{\mathscr H_y(t)}{\mathscr H_y(t/2)}
			+C\alpha_1\varepsilon(t).
		\end{equation}
		These integrals are now all centered at $y$.
		
		Fix $0<r\le r_{10}$ and put $t_j=2^{1-j}r$. By
		\eqref{eq:adjacent-increment-comparison},
		\begin{equation}\label{eq:dyadic-recurrence}
			\log\frac{\mathscr H_y(t_{j+1})}{\mathscr H_y(t_{j+2})}
			\le(1+C\varepsilon(t_j))
			\log\frac{\mathscr H_y(t_j)}{\mathscr H_y(t_{j+1})}
			+C\alpha_1\varepsilon(t_j).
		\end{equation}
		The Dini condition gives
		\[
		\sum_{j=0}^\infty\varepsilon(2^{1-j}r)
		\le C\int_0^{Cr}\frac{\omega_{\partial\Omega}(s)}s\,ds
		+Cr^{\tau_2}\le C.
		\]
		In particular,
		\[
		\prod_{j=0}^\infty(1+C\varepsilon(t_j))
		\le\exp\left(C\sum_{j=0}^\infty\varepsilon(t_j)\right)\le C.
		\]
		Iteration of \eqref{eq:dyadic-recurrence} gives
		\begin{align*}
			\log\frac{\mathscr H_y(t_{j+1})}{\mathscr H_y(t_{j+2})}
			\le{}&\prod_{i=0}^{j}(1+C\varepsilon(t_i))
			\log\frac{\mathscr H_y(2r)}{\mathscr H_y(r)}\\
			&+C\alpha_1\sum_{i=0}^{j}\varepsilon(t_i)
			\prod_{k=i+1}^{j}(1+C\varepsilon(t_k)).
		\end{align*}
		Thus every dyadic weighted increment below $2r$ is bounded by
		\begin{equation}\label{eq:dyadic-iteration-explicit-potential}
			C\left(\alpha_1+
			\log\frac{\mathscr H_y(2r)}{\mathscr H_y(r)}\right).
		\end{equation}
		
		To bound the initial increment, fix $0<q\le2r_{10}$ and use the
		three-ball inequality at the shifted center $y_q$, with radii
		$q/2-|y_q-y|$, $q+|y_q-y|$, and $2q-|y_q-y|$.
		After division by the outer scale $4q$, these ratios remain in a fixed
		compact set. The center inclusions are
		\[
		B_{q/2-|y_q-y|}(y_q)\subset B_{q/2}(y),\qquad
		B_q(y)\subset B_{q+|y_q-y|}(y_q),\qquad
		B_{2q-|y_q-y|}(y_q)\subset B_{2q}(y).
		\]
		Lemma~\ref{lem:fixed-ratio-three-ball-schrodinger}, applied at $y_q$
		with outer scale $4q$, therefore gives a fixed $0<\vartheta<1$ such that
		\begin{align*}
			I_u(y,q)
			&\le I_u(y_q,q+|y_q-y|)\\
			&\le e^{C\alpha_1}I_u(y_q,q/2-|y_q-y|)^\vartheta
			I_u(y_q,2q-|y_q-y|)^{1-\vartheta}\\
			&\le e^{C\alpha_1}I_u(y,q/2)^\vartheta
			I_u(y,2q)^{1-\vartheta}.
		\end{align*}
		Taking logarithms in the preceding inequality gives
		\[
		\vartheta\log\frac{I_u(y,q)}{I_u(y,q/2)}
		\le C\alpha_1+(1-\vartheta)
		\log\frac{I_u(y,2q)}{I_u(y,q)}.
		\]
		Hence
		\begin{equation}\label{eq:one-step-ordinary-doubling}
			\log\frac{I_u(y,q)}{I_u(y,q/2)}
			\le C\left(\alpha_1+
			\log\frac{I_u(y,2q)}{I_u(y,q)}\right).
		\end{equation}
		On the other hand,
		\[
		\log\frac{\mathscr H_y(2r)}{\mathscr H_y(r)}
		\le\log\frac{I_u(y,2r)}{I_u(y,r/2)}+C\alpha_1.
		\]
		Two applications of \eqref{eq:one-step-ordinary-doubling} bound this
		by the right-hand side of
		\eqref{eq:dyadic-dini-recentering-conclusion}.
		Finally,
		\[
		\log\frac{I_u(y,s)}{I_u(y,s/2)}
		\le\log\frac{\mathscr H_y(2s)}{\mathscr H_y(s/2)}.
		\]
		At $s=2^{-j}r$, the right-hand side is exactly
		\[
		\log\frac{\mathscr H_y(2^{-j+1}r)}
		{\mathscr H_y(2^{-j}r)}
		+\log\frac{\mathscr H_y(2^{-j}r)}
		{\mathscr H_y(2^{-j-1}r)},
		\]
		and both increments are bounded by
		\eqref{eq:dyadic-iteration-explicit-potential}.  If
		$2^{-j-1}r<s\le2^{-j}r$, monotonicity of $I_u(y,\cdot)$ gives
		\[
		\frac{I_u(y,s)}{I_u(y,s/2)}
		\le\frac{I_u(y,2^{-j}r)}{I_u(y,2^{-j-2}r)}.
		\]
		The last logarithm is
		\[
		\log\frac{I_u(y,2^{-j}r)}{I_u(y,2^{-j-2}r)}
		=\log\frac{I_u(y,2^{-j}r)}{I_u(y,2^{-j-1}r)}
		+\log\frac{I_u(y,2^{-j-1}r)}{I_u(y,2^{-j-2}r)}.
		\]
		For the two ratios on the right, respectively,
		\[
		\begin{aligned}
			\log\frac{I_u(y,2^{-j}r)}{I_u(y,2^{-j-1}r)}
			&\le
			\log\frac{\mathscr H_y(2^{-j+1}r)}{\mathscr H_y(2^{-j}r)}
			+\log\frac{\mathscr H_y(2^{-j}r)}{\mathscr H_y(2^{-j-1}r)},\\
			\log\frac{I_u(y,2^{-j-1}r)}{I_u(y,2^{-j-2}r)}
			&\le
			\log\frac{\mathscr H_y(2^{-j}r)}{\mathscr H_y(2^{-j-1}r)}
			+\log\frac{\mathscr H_y(2^{-j-1}r)}{\mathscr H_y(2^{-j-2}r)}.
		\end{aligned}
		\]
		All four weighted increments occur among the neighboring dyadic
		increments controlled by
		\eqref{eq:dyadic-iteration-explicit-potential}. This proves the lemma.
	\end{proof}
	
	\begin{definition}\label{def:fixed-scale-boundary-doubling}
		Let $r_{10}$ be the radius supplied by Lemma
		\ref{lem:dyadic-dini-recentering}, decreased if necessary, and set
		\begin{equation}\label{eq:fixed-scale-radius}
			r_{11}=2r_{10}.
		\end{equation}
		Define
		\begin{equation}\label{eq:fixed-scale-boundary-doubling}
			\mathcal D_0
			=\sup_{y\in\overline\Omega}D_\Omega^u(y,r_{11}).
		\end{equation}
	\end{definition}
	 
	 Combing Lemma \ref{lem:dyadic-dini-recentering} and  propagation of smallness in Lemma \ref{lem:uniform-mass-propagation} below, we can now prove the upper bound of doubling index in \eqref{eq:uniform-schrodinger-doubling}.
	 
	\begin{theorem}\label{thm:small-scale-schrodinger-doubling-fixed-scale}
		Under \eqref{eq:boundary-schrodinger-equation}--\eqref{eq:boundary-potential-bound},
		\begin{equation}\label{eq:small-scale-schrodinger-doubling-fixed-scale}
			\sup_{y\in\overline\Omega}\sup_{0<s\le r_{11}/4}
			D_\Omega^u(y,s)
			\le C\bigl(1+M^{\frac{2p}{3p-2n}}
			+\mathcal D_0\bigr).
		\end{equation}
		Using Lemma \ref{lem:uniform-mass-propagation}, we finally reach
		\begin{equation}\label{eq:uniform-schrodinger-doubling}
			\sup_{y\in\overline\Omega}\sup_{0<s\le r_{11}/4}
			D_\Omega^u(y,s)
			\le C\left(1+M^{\frac{2p}{3p-2n}}\right).
		\end{equation}
	\end{theorem}
	
	\begin{proof}
		For $0<s\le r_{11}/4=r_{10}/2$, apply
		Lemma~\ref{lem:dyadic-dini-recentering} with $r=r_{10}$ and $2s$ in
		place of its inner radius. Since $2s\le r_{10}$,
		\begin{align*}
			\log\frac{I_u(y,2s)}{I_u(y,s)}
			&\le C\left(\alpha_1+
			\log\frac{I_u(y,4r_{10})}{I_u(y,2r_{10})}\right)\\
			&=C\left(\alpha_1+
			\log\frac{I_u(y,2r_{11})}{I_u(y,r_{11})}\right).
		\end{align*}
		Uniform volume density gives
		\[
		cr^n\le|\Omega\cap B_r(y)|\le Cr^n
		\qquad(0<r\le2r_{11}),
		\]
		so the definition of $\mathcal D_0$ implies
		\[
		\log\frac{I_u(y,2r_{11})}{I_u(y,r_{11})}
		\le(\log4)\mathcal D_0+C.
		\]
		Using the same density estimate at $s$ and $2s$, we conclude that
		\begin{align*}
			D_\Omega^u(y,s)
			&=\frac1{\log4}\left[
			\log\frac{I_u(y,2s)}{I_u(y,s)}
			+\log\frac{|\Omega\cap B_s(y)|}{|\Omega\cap B_{2s}(y)|}
			\right]\\
			&\le C\bigl(\alpha_1+1+\mathcal D_0\bigr)
			\le C\left(1+M^{\frac{2p}{3p-2n}}+\mathcal D_0\right).
		\end{align*}
		Taking the supremum over $y$ and $s$ proves the theorem.
	\end{proof}
	
	\begin{lemma}\label{lem:uniform-mass-propagation}
		For the global Dirichlet solution,
		\begin{equation}\label{eq:fixed-scale-doubling-global-polynomial}
			\mathcal D_0
			\le C\left(1+M^{\frac{2p}{3p-2n}}\right).
		\end{equation}
	\end{lemma}
	
	\begin{proof}
		With $r_{11}$ defined in \eqref{eq:fixed-scale-radius}, normalize
		$\int_\Omega u^2=1$ and put
		$\rho=r_{11}/512$. A finite cover of $\overline\Omega$
		by $N_\Omega$ balls of radius $\rho/2$ contains a ball with mass at least
		$1/N_\Omega$, because
		\[
		1=\int_\Omega u^2
		\le\sum_{k=1}^{N_\Omega}I_u(x_k,\rho/2).
		\]
		Connectedness gives, from any prescribed center $y$ to
		that ball, a chain $x_0=y,\ldots,x_J$ with $J\le N_\Omega$ and
		$|x_{j+1}-x_j|<\rho$. The cover and its intersection graph depend only
		on the domain.
		
		Apply Lemma~\ref{lem:dini-boundary-star-shaped} at $x_j$ with radius
		$8\rho$, obtaining $\widetilde x_j$ with
		$|\widetilde x_j-x_j|\le\rho/4$ and a star-shaped intersection at radius
		$64\rho$. The inclusions
		\[
		B_{\rho/2}(\widetilde x_j)\subset B_\rho(x_j),
		\qquad B_\rho(x_{j+1})\subset B_{4\rho}(\widetilde x_j)
		\]
		and the three-ball inequality at radii $\rho/2,4\rho,32\rho$ give
		\begin{align*}
			I_u(x_{j+1},\rho)
			&\le I_u(\widetilde x_j,4\rho)\\
			&\le e^{C\alpha_1}I_u(\widetilde x_j,\rho/2)^{\vartheta_0}
			I_u(\widetilde x_j,32\rho)^{1-\vartheta_0}\\
			&\le e^{C\alpha_1}I_u(x_j,\rho)^{\vartheta_0},
			\qquad 0<\vartheta_0<1,
		\end{align*}
		where the last line uses $I_u(\widetilde x_j,32\rho)\le
		\int_\Omega u^2=1$. The interpolation exponent
		$\vartheta_0$ is fixed, since the three radii divided by $64\rho$
		are $1/128,1/16,1/2$. In particular, it is independent of $M$ and $u$.
		Since $0<I_u(x_j,\rho)\le1$, taking logarithms gives
		\begin{equation}\label{eq:fixed-scale-logarithmic-one-link}
			-\log I_u(x_j,\rho)
			\le\vartheta_0^{-1}
			\bigl(-\log I_u(x_{j+1},\rho)+C\alpha_1\bigr).
		\end{equation}
		The terminal ball has mass at least $1/N_\Omega$. Iterating backwards,
		\begin{equation}\label{eq:fixed-scale-logarithmic-iteration}
			-\log I_u(y,\rho)
			\le\vartheta_0^{-J}\log N_\Omega
			+C\alpha_1\sum_{j=1}^{J}\vartheta_0^{-j}
			\le C\alpha_1.
		\end{equation}
		The last constant is geometric because $J\le N_\Omega$. Consequently,
		\[
		I_u(y,r_{11})\ge I_u(y,\rho)
		\ge e^{-C\alpha_1}.
		\]
		Since $I_u(y,2r_{11})\le1$, while uniform volume density gives
		$cr_{11}^n\le|\Omega\cap B_{r_{11}}(y)|$ and
		$|\Omega\cap B_{2r_{11}}(y)|\le Cr_{11}^n$, it follows that
		\begin{align*}
			4^{D_\Omega^u(y,r_{11})}
			&=\frac{I_u(y,2r_{11})}{I_u(y,r_{11})}
			\frac{|\Omega\cap B_{r_{11}}(y)|}
			{|\Omega\cap B_{2r_{11}}(y)|}\\
			&\le Ce^{C\alpha_1}.
		\end{align*}
		Taking $\log_4$, then the supremum over $y$, and using
		$\alpha_1=2+M^{2p/(3p-2n)}$ proves
		\eqref{eq:fixed-scale-doubling-global-polynomial}.
	\end{proof}
	
	\medskip
	
	\subsection{Boundary flattening and reflection}
	From this subsection, we are using reflection to reach interior problem \eqref{eq:main-equation}.
	Fix $x_0\in\partial\Omega$. Translate and rotate so that $x_0=0$ and,
	in a boundary chart,
	\begin{equation}\label{eq:boundary-graph-for-flattening}
		\Omega=\{(x',x_n):x_n>\varphi(x')\},
		\qquad\varphi(0)=0,\quad\nabla\varphi(0)=0.
	\end{equation}
	The inward unit normal is
	\[
	\mathbf n(y')=
	\frac{(-\nabla\varphi(y'),1)}{\sqrt{1+|\nabla\varphi(y')|^2}}.
	\]
	Choose a smooth nonnegative kernel $\rho$, supported in $B_1^{n-1}$,
	with integral one, and write $\rho_s(z')=s^{1-n}\rho(z'/s)$.
	Following \cite{KenigZhao2025}, define
	\begin{equation}\label{eq:boundary-flattening-map}
		G(y',s)=(y',\varphi(y'))+s(\rho_s*\mathbf n)(y'),
		\qquad s>0,
	\end{equation}
	with $G(y',0)=(y',\varphi(y'))$.
	
	\begin{lemma}\label{lem:flattening-map-properties}
		There is a geometric radius $r_0>0$ such that $G$ is a $C^1$
		bi-Lipschitz diffeomorphism on $B_{16r_0}^+$, maps its interior into
		$\Omega$, and satisfies
		\begin{align}
			L^{-1}|Y-Z|\le |G(Y)-G(Z)|&\le L|Y-Z|,
			\label{eq:G-bilipschitz}\\
			0<J_-\le J_G(Y):=\det DG(Y)&\le J_+,
			\label{eq:G-jacobian-bounds}\\
			|DG(Y)-DG(Z)|&\le C\omega_{\partial\Omega}(C|Y-Z|).
			\label{eq:DG-dini-modulus}
		\end{align}
		After reducing the radius,
		\begin{equation}\label{eq:G-image-inclusions}
			\Omega\cap B_{cr_0}\subset G(B_{r_0}^+),
			\qquad G(B_{4r_0}^+)\subset\Omega\cap B_{Cr_0}.
		\end{equation}
	\end{lemma}
	
	\begin{proof}
		Put $\mathbf n_s=\rho_s*\mathbf n$. Change of variables gives
		\begin{align}
			\partial_iG(y',s)
			&=(e_i,\partial_i\varphi(y'))
			+\int_{B_1}\partial_i\rho(q)
			\bigl(\mathbf n(y'-sq)-\mathbf n(y')\bigr)\,dq,
			\label{eq:G-tangential-derivative-kernel}\\
			\partial_sG(y',s)
			&=\int_{B_1}\rho(q)\mathbf n(y'-sq)\,dq
			+\int_{B_1}\bigl(-(n-1)\rho(q)-q\cdot\nabla\rho(q)\bigr)
			\bigl(\mathbf n(y'-sq)-\mathbf n(y')\bigr)\,dq.
			\label{eq:G-normal-derivative-kernel}
		\end{align}
		Indeed, every derivative of $\rho_s$ has integral zero, and the same
		identities before the change of variables read
		\begin{align*}
			s\partial_i\mathbf n_s(y')
			&=s\int\partial_i\rho_s(z')
			\bigl(\mathbf n(y'-z')-\mathbf n(y')\bigr)\,dz',\\
			s\partial_s\mathbf n_s(y')
			&=s\int\partial_s\rho_s(z')
			\bigl(\mathbf n(y'-z')-\mathbf n(y')\bigr)\,dz'.
		\end{align*}
		The support of the kernels is $|z'|\le s$, and their $L^1$ derivative
		norms are bounded by $C/s$. Hence
		\[
		|s\nabla_{y'}\mathbf n_s(y')|
		+|s\partial_s\mathbf n_s(y')|
		\le C\omega_{\partial\Omega}(Cs).
		\]
		Consequently,
		\begin{equation}\label{eq:DG-boundary-limit}
			\partial_iG(y',0)=(e_i,\partial_i\varphi(y')),
			\qquad\partial_sG(y',0)=\mathbf n(y').
		\end{equation}
		For $Y=(y',s)$ and $Z=(z',t)$, each difference of the integral terms
		in \eqref{eq:G-tangential-derivative-kernel} is bounded by
		\begin{align*}
			&C\int_{B_1}\bigl(
			|\mathbf n(y'-sq)-\mathbf n(z'-tq)|
			+|\mathbf n(y')-\mathbf n(z')|\bigr)\,dq\\
			&\qquad\le C\omega_{\partial\Omega}
			(C(|y'-z'|+|s-t|))
			\le C\omega_{\partial\Omega}(C|Y-Z|).
		\end{align*}
		The same estimate applied to the two integrals in
		\eqref{eq:G-normal-derivative-kernel}, together with
		$|\nabla\varphi(y')-\nabla\varphi(z')|
		\le C\omega_{\partial\Omega}(C|y'-z'|)$, gives
		\[
		|DG(Y)-DG(Z)|
		\le C\omega_{\partial\Omega}(C|Y-Z|).
		\]
		this is the derivative estimate in
		\cite[Section~4, equations~(4.1)--(4.3)]{KenigZhao2025}.
		The formula
		\[
		\widehat G(y',s)=(y',\varphi(y'))
		+s(\rho_{|s|}*\mathbf n)(y')
		\]
		extends $G$ to a $C^1$ map across $s=0$.
		Since $DG(0)=I$, choose $r_0$ so that
		$\sup_{B_{16r_0}}|DG-I|\le1/4$.  The fundamental theorem of calculus
		along line segments then gives
		\[
		\frac34|Y-Z|\le|G(Y)-G(Z)|\le\frac54|Y-Z|,
		\]
		and all singular values of $DG$ lie in $[3/4,5/4]$. Thus
		$(3/4)^n\le J_G\le(5/4)^n$ and the quantitative inverse theorem makes
		$G$ a bi-Lipschitz diffeomorphism on the smaller half-ball.
		
		Finally,
		\[
		\partial_s(G_n-\varphi(G'))(y',0)
		=\sqrt{1+|\nabla\varphi(y')|^2}>0.
		\]
		After one further reduction of $r_0$, continuity makes this derivative
		at least $1/2$ for $|s|\le16r_0$. Hence
		\[
		G_n(y',s)-\varphi(G'(y',s))
		=\int_0^s\partial_t(G_n-\varphi(G'))(y',t)\,dt\ge s/2,
		\]
		so $G$ maps the open half-ball into $\Omega$; the same integral for
		$s<0$ is negative, so $\widehat G(y',s)$ lies below the graph.
		Moreover, for $|Y|=r_0$,
		$|\widehat G(Y)|\ge3r_0/4$. Since $\widehat G$ is injective and open,
		the boundary of $\widehat G(B_{r_0})$ is
		$\widehat G(\partial B_{r_0})$; consequently
		$B_{3r_0/4}\subset\widehat G(B_{r_0})$. If
		$x\in\Omega\cap B_{3r_0/4}$, its unique preimage therefore has
		positive last coordinate by the sign just proved. Thus
		$\Omega\cap B_{3r_0/4}\subset G(B_{r_0}^+)$. Finally,
		$|G(Y)|\le(5/4)|Y|$ for $Y\in B_{4r_0}^+$, so
		$G(B_{4r_0}^+)\subset\Omega\cap B_{5r_0}$. This proves
		\eqref{eq:G-image-inclusions} with fixed $c$ and $C$.
	\end{proof}
	
	\begin{proposition}\label{prop:equation-after-flattening}
		Put $\bar u=u\circ G$. Then
		\begin{equation}\label{eq:flattened-equation}
			\operatorname{div}(\bar a\nabla \bar u)+\bar V \bar u=0
			\quad\text{in }B_{16r_0}^+,
		\end{equation}
		with zero trace on the flat boundary, where
		\begin{equation}\label{eq:flattened-coefficient}
			\bar a=J_GDG^{-1}DG^{-T},
			\qquad \bar V=J_G(V\circ G).
		\end{equation}
		The matrix is symmetric, uniformly elliptic, and satisfies
		\begin{equation}\label{eq:flattened-coefficient-modulus}
			|\bar a(Y)-\bar a(Z)|
			\le C\omega_{\partial\Omega}(C|Y-Z|),
			\qquad \|\bar V\|_{L^p(B_{16r_0}^+)}\le CM.
		\end{equation}
		On $s=0$,
		\begin{equation}\label{eq:flattened-boundary-block}
			\bar a(y',0)=\sqrt{1+|\nabla\varphi|^2}
			\begin{pmatrix}
				(I_{n-1}+\nabla\varphi\otimes\nabla\varphi)^{-1}&0\\
				0&1
			\end{pmatrix}.
		\end{equation}
	\end{proposition}
	
	\begin{proof}
		The principal-coefficient transformation is that of
		\cite[Section~4]{KenigZhao2025}. To verify the equation with the
		potential, let $\eta$ be a test function on the half-ball, vanishing
		on its flat face. The chain rule gives
		\[
		\nabla u(G(Y))=DG(Y)^{-T}\nabla \bar u(Y),
		\qquad
		\nabla(\eta\circ G^{-1})(G(Y))=DG(Y)^{-T}\nabla\eta(Y).
		\]
		Changing variables in the weak equation therefore gives
		\[
		\int_{B_{16r_0}^+}J_GDG^{-1}DG^{-T}
		\nabla \bar u\cdot\nabla\eta
		-\int_{B_{16r_0}^+}J_G(V\circ G)\bar u\eta=0.
		\]
		Tests with zero trace on the flat face are admissible by zero extension
		and density. This proves \eqref{eq:flattened-equation}, and the trace
		of $\bar u$ is zero because $G$ maps the flat face to $\partial\Omega$.
		The transformed potential satisfies
		\begin{align*}
			\|\bar V\|_{L^p(B_{16r_0}^+)}^p
			&=\int_{B_{16r_0}^+}J_G(Y)^p|V(G(Y))|^p\,dY\\
			&=\int_{G(B_{16r_0}^+)}
			J_G(G^{-1}(x))^{p-1}|V(x)|^p\,dx
			\le J_+^{p-1}M^p.
		\end{align*}
		For every vector $\xi$,
		\[
		J_-\|DG\|^{-2}|\xi|^2
		\le \bar a\xi\cdot\xi
		=J_G|DG^{-T}\xi|^2
		\le J_+\|DG^{-1}\|^2|\xi|^2,
		\]
		so the bi-Lipschitz and Jacobian bounds give uniform ellipticity.
		The identities
		\[
		DG(Y)^{-1}-DG(Z)^{-1}
		=DG(Y)^{-1}\bigl(DG(Z)-DG(Y)\bigr)DG(Z)^{-1}
		\]
		and
		\[
		|\det DG(Y)-\det DG(Z)|\le C|DG(Y)-DG(Z)|.
		\]
		Expanding the three factors in
		\[
		(\det DG(Y))DG(Y)^{-1}DG(Y)^{-T}
		-(\det DG(Z))DG(Z)^{-1}DG(Z)^{-T}
		\]
		therefore gives
		\[
		|\bar a(Y)-\bar a(Z)|
		\le C|DG(Y)-DG(Z)|
		\le C\omega_{\partial\Omega}(C|Y-Z|).
		\]
		Finally, the boundary columns in \eqref{eq:DG-boundary-limit} yield
		\[
		DG(y',0)^TDG(y',0)
		=\operatorname{diag}(I_{n-1}+\nabla\varphi\otimes\nabla\varphi,1),
		\]
		with $J_G(y',0)=\sqrt{1+|\nabla\varphi(y')|^2}$.
		Taking the inverse Gram matrix proves \eqref{eq:flattened-boundary-block}.
	\end{proof}
	
	Let $R=\operatorname{diag}(1,\ldots,1,-1)$. Extend $\bar u$
	oddly and $\bar V$ evenly across $s=0$, and extend the matrix by
	\begin{equation}\label{eq:reflected-coefficient}
		\bar a(Y)=
		\begin{cases}
			\bar a(Y),&Y_n\ge0,\\
			R\,\bar a(RY)R,&Y_n<0.
		\end{cases}
	\end{equation}
	We retain the same notation for the extended functions and matrix.
	
	\begin{proposition}\label{prop:reflected-equation}
		The odd extension belongs to $W^{1,2}(B_{16r_0})$ and satisfies
		\begin{equation}\label{eq:reflected-full-ball-equation}
			\operatorname{div}(\bar a\nabla \bar u)
			+\bar V\bar u=0
			\quad\text{in }B_{16r_0}.
		\end{equation}
		Its coefficients are symmetric and uniformly elliptic, with
		\begin{equation}\label{eq:reflected-coefficient-dini}
			|\bar a(Y)-\bar a(Z)|
			\le C\omega_{\partial\Omega}(C|Y-Z|),
			\qquad\|\bar V\|_{L^p(B_{16r_0})}\le CM.
		\end{equation}
	\end{proposition}
	
	\begin{proof}
		The zero trace permits odd reflection in $W^{1,2}$.
		For $Y_n>0$, the extension formulas give
		\[
		\nabla \bar u(RY)
		=-R\nabla \bar u(Y),\qquad
		\bar V(RY)=\bar V(Y).
		\]
		For $\eta\in C_c^\infty(B_{16r_0})$, changing variables by
		$R$ in the lower half-ball yields
		\begin{align*}
			&\int_{B_{16r_0}}\bar a
			\nabla \bar u\cdot\nabla\eta
			-\int_{B_{16r_0}}\bar V\bar u\eta\\
			&\qquad=\int_{B_{16r_0}^+}\bar a\nabla \bar u\cdot
			\nabla(\eta-\eta\circ R)
			-\int_{B_{16r_0}^+}\bar V \bar u
			(\eta-\eta\circ R)=0.
		\end{align*}
		The difference of test functions has zero trace on the flat face.
		This proves the weak equation across that face and shows why the
		potential must be extended evenly. In particular,
		\[
		\|\bar V\|_{L^p(B_{16r_0})}^p
		=2\|\bar V\|_{L^p(B_{16r_0}^+)}^p\le CM^p.
		\]
		Orthogonal conjugation preserves ellipticity. The diagonal blocks of the
		matrix are even and the mixed entries are odd; the latter vanish on the
		flat face by \eqref{eq:flattened-boundary-block}. If $Y_nZ_n<0$, put
		$Y_0=(Y',0)$ and $Z_0=(Z',0)$. Then
		\begin{align*}
			|\bar a(Y)-\bar a(Z)|
			&\le|\bar a(Y)-\bar a(Y_0)|
			+|\bar a(Y_0)-\bar a(Z_0)|
			+|\bar a(Z_0)-\bar a(Z)|\\
			&\le C\bigl[
			\omega_{\partial\Omega}(C|Y-Y_0|)
			+\omega_{\partial\Omega}(C|Y_0-Z_0|)
			+\omega_{\partial\Omega}(C|Z-Z_0|)\bigr]\\
			&\le C\omega_{\partial\Omega}(C|Y-Z|).
		\end{align*}
		The last inequality uses
		$|Y-Y_0|,|Y_0-Z_0|,|Z-Z_0|\le|Y-Z|$ and monotonicity of the
		modulus. On one side of the face the same estimate is
		\eqref{eq:flattened-coefficient-modulus}. This proves
		\eqref{eq:reflected-coefficient-dini}.
	\end{proof}
	
	\subsection{Transfer of doubling and proof of the boundary theorem}
	Let $r_0,L,J_-,J_+$ be the constants supplied by Lemma
	\ref{lem:flattening-map-properties}, and let $\mathcal D_0$ be defined in
	\eqref{eq:fixed-scale-boundary-doubling}.
	
	\begin{proposition}\label{prop:doubling-transfer-reflection}
		After decreasing the geometric radius $r_0$, for every
		$B_{2s}(z)\subset B_{4r_0}$,
		\begin{equation}\label{eq:reflected-euclidean-doubling-fixed-scale}
			\log_4\frac{\vint_{B_{2s}(z)}(\bar u)^2}
			{\vint_{B_s(z)}(\bar u)^2}
			\le C\bigl(1+M^{\frac{2p}{3p-2n}}
			+\mathcal D_0\bigr).
		\end{equation}
		In particular, for the global Dirichlet solution,
		\begin{equation}\label{eq:reflected-doubling-explicit-V}
			\sup_{B_{2s}(z)\subset B_{4r_0}}
			\log_4\frac{\vint_{B_{2s}(z)}(\bar u)^2}
			{\vint_{B_s(z)}(\bar u)^2}
			\le C\left(1+M^{\frac{2p}{3p-2n}}\right).
		\end{equation}
	\end{proposition}
	
	\begin{proof}
		Put $z^+=(z',|z_n|)$ and $x=G(z^+)$. Since the squared solution is even,
		folding the ball into the upper half-space gives
		\[
		\int_{B_s(z^+)\cap\{Y_n>0\}}(\bar u)^2
		\le\int_{B_s(z)}(\bar u)^2
		\le2\int_{B_s(z^+)\cap\{Y_n>0\}}(\bar u)^2.
		\]
		The bi-Lipschitz and Jacobian bounds therefore imply
		\begin{equation}\label{eq:coordinate-original-mass-comparison}
			cI_u(x,s/L)\le\int_{B_s(z)}(\bar u)^2
			\le CI_u(x,Ls).
		\end{equation}
		Indeed, the signed extension $\widehat G$ is bi-Lipschitz on
		$B_{16r_0}$ and $z^+\in B_{4r_0}$. Hence
		\[
		\operatorname{dist}\bigl(x,\widehat G(\partial B_{16r_0})\bigr)
		\ge L^{-1}(16r_0-|z^+|)\ge12L^{-1}r_0.
		\]
		Since $s\le2r_0$, every $w\in\Omega\cap B_{s/L}(x)$ belongs to the
		image of $\widehat G$. Its unique preimage $W$ has positive last
		coordinate by the sign property proved in
		Lemma~\ref{lem:flattening-map-properties}, and
		$|W-z^+|\le L|w-x|<s$. Conversely,
		$|G(W)-x|\le L|W-z^+|$ for
		$W\in B_s(z^+)\cap\{Y_n>0\}$. Therefore
		\[
		\Omega\cap B_{s/L}(x)
		\subset G(B_s(z^+)\cap\{Y_n>0\})
		\subset\Omega\cap B_{Ls}(x),
		\]
		which proves \eqref{eq:coordinate-original-mass-comparison} after the
		Jacobian bounds are applied.
		
		Choose $m=\lceil\log_2(2L^2)\rceil$ and reduce $r_0$ so that
		$2^{m+1}r_0/L\le r_{11}/4$, where $r_{11}$ is defined in
		\eqref{eq:fixed-scale-radius}. Apply
		\eqref{eq:small-scale-schrodinger-doubling-fixed-scale} at the finitely
		many radii $s/L,2s/L,\ldots,2^{m-1}s/L$. Since $s\le2r_0$,
		all these radii lie in the range of that estimate. Uniform volume
		density and telescoping the average ratios give
		\begin{align*}
			\frac{I_u(x,2^ms/L)}{I_u(x,s/L)}
			={}&\frac{|\Omega\cap B_{2^ms/L}(x)|}{|\Omega\cap B_{s/L}(x)|}
			\prod_{j=0}^{m-1}
			\frac{\vint_{\Omega\cap B_{2^{j+1}s/L}(x)}u^2}
			{\vint_{\Omega\cap B_{2^js/L}(x)}u^2}\\
			\le{}&C2^{mn}
			4^{Cm(1+M^{2p/(3p-2n)}+\mathcal D_0)}.
		\end{align*}
		Now $2Ls\le2^ms/L$, so monotonicity of the unnormalized mass and
		\eqref{eq:coordinate-original-mass-comparison} imply
		\[
		\frac{\vint_{B_{2s}(z)}(\bar u)^2}
		{\vint_{B_s(z)}(\bar u)^2}
		\le C\frac{I_u(x,2Ls)}{I_u(x,s/L)}.
		\]
		Taking logarithms proves
		\eqref{eq:reflected-euclidean-doubling-fixed-scale}. Substitution of
		\eqref{eq:fixed-scale-doubling-global-polynomial} gives
		\eqref{eq:reflected-doubling-explicit-V}.
	\end{proof}
	
	\begin{proof}[Proof of Theorem~\ref{thm:boundary-volume-estimate}]
		Let $r_0$ be the geometric radius from Lemma
		\ref{lem:flattening-map-properties}, and set $U(Y)=\bar u(r_0Y)$.
		Its equation in $B_2$ has leading
		matrix $\bar a(r_0Y)$ and potential
		$r_0^2\bar V(r_0Y)$. By the preceding propositions,
		\begin{align}
			\|r_0^2\bar V(r_0\,\cdot)\|_{L^p(B_2)}
			&\le Cr_0^{2-n/p}M,
			\label{eq:rescaled-boundary-potential-bound}\\
			\widetilde\omega(t)&=C\omega_{\partial\Omega}(Cr_0t),
			\label{eq:rescaled-boundary-modulus}\\
			\sup_{B_{2s}(z)\subset B_2}
			\log_4\frac{\vint_{B_{2s}(z)}U^2}{\vint_{B_s(z)}U^2}
			&\le C\left(1+M^{\frac{2p}{3p-2n}}\right).
			\label{eq:rescaled-boundary-doubling-global}
		\end{align}
		Reduce $r_0$ so that $\widetilde\omega(1)\le1$.
		Theorem~\ref{thm:main-volume-estimate} now applies to $U$.
		
		For $Y_n>0$, the chain rule gives
		\[
		\nabla U(Y)=r_0DG(r_0Y)^T\nabla u(G(r_0Y)).
		\]
		Invertibility of $DG$ therefore gives
		\begin{equation}\label{eq:singular-set-under-flattening}
			U(Y)=0,\ \nabla U(Y)=0
			\quad\Longleftrightarrow\quad
			G(r_0Y)\in S_\Omega(u).
		\end{equation}
		Choose $c>0$ sufficiently small that
		$\Omega\cap B_{4cr_0}(x_0)\subset G(r_0B_{1/2}^+)$.
		If $x\in\Omega\cap B_{2cr_0}(x_0)$ and
		$q\in S_\Omega(u)\cap B_{4cr_0}(x_0)$ satisfy $|x-q|<r$, then
		\[
		\left|r_0^{-1}G^{-1}(x)-r_0^{-1}G^{-1}(q)\right|
		\le\frac{Lr}{r_0}.
		\]
		Both preimages belong to $B_{1/2}^+$, and the second is singular for $U$.
		After reducing $c$, $Lr/r_0\le1$. Since
		$\det D(G(r_0\,\cdot))=r_0^nJ_G(r_0\,\cdot)\le J_+r_0^n$,
		change of variables gives
		\begin{align*}
			&\left|\left\{x\in\Omega\cap B_{2cr_0}(x_0):
			\operatorname{dist}(x,S_\Omega(u)\cap B_{4cr_0}(x_0))<r\right\}\right|\\
			&\qquad\le J_+r_0^n
			\left|\left\{Y\in B_1:
			\operatorname{dist}(Y,S(U))<\frac{Lr}{r_0}\right\}\right|.
		\end{align*}
		Let $\delta_2$ be the parameter supplied by the interior theorem for
		the rescaled equation.  The bounds
		\eqref{eq:rescaled-boundary-potential-bound} and
		\eqref{eq:rescaled-boundary-doubling-global} give
		\[
		\delta_2\ge c(1+r_0^{2-n/p}M)^{-C}
		\exp\left[-C\left(1+
		M^{\frac{2p}{3p-2n}}\right)^2\right],
		\]
		where $c,C>0$ are structural constants obtained from
		\eqref{eq:5-final-delta} after applying Theorem
		\ref{thm:main-volume-estimate} to the rescaled equation and inserting
		\eqref{eq:rescaled-boundary-potential-bound} and
		\eqref{eq:rescaled-boundary-doubling-global}. By decreasing $c$ and
		increasing $C$ in
		\eqref{eq:intro-boundary-delta-global}, we may arrange that
		\[
		0<\widetilde{\delta}
		\le c(1+r_0^{2-n/p}M)^{-C}
		\exp\left[-C\left(1+
		M^{\frac{2p}{3p-2n}}\right)^2\right]
		\le\delta_2.
		\]
		Since
		$\rho_{\widetilde\omega}$ is nondecreasing,
		$\rho_{\widetilde\omega}(\widetilde{\delta})\le
		\rho_{\widetilde\omega}(\delta_2)$. The interior theorem at radius
		$Lr/r_0$ consequently yields
		\[
		\begin{split}
			&\left|\left\{x\in\Omega\cap B_{2cr_0}(x_0):
			\operatorname{dist}\bigl(x,S_\Omega(u)\cap B_{4cr_0}(x_0)\bigr)<r
			\right\}\right|\\
			&\quad\le
			\frac{Cr_0^{n-2}r^2}
			{\widetilde{\delta}\,\rho_{\widetilde\omega}(\widetilde{\delta})^2},
		\end{split}
		\]
		where we also used
		\[
		(1+r_0^{2-n/p}M)^C
		\exp\left[C\left(1+M^{\frac{2p}{3p-2n}}\right)^2\right]
		\le\frac C{\widetilde{\delta}}.
		\]
		Restriction to the smaller balls
		proves \eqref{eq:intro-boundary-volume-estimate}.
		
		Put $\Sigma=S_\Omega(u)\cap B_{cr_0}(x_0)$ and let
		$x_1,\ldots,x_N$ be a maximal $r$-separated subset of $\Sigma$.
		For $0<r\le cr_0$, the sets
		$\Omega\cap B_{r/4}(x_i)$ are pairwise disjoint, lie in the
		$r/4$-neighborhood of $\Sigma$ inside $B_{2cr_0}(x_0)$, and satisfy
		$|\Omega\cap B_{r/4}(x_i)|\ge cr^n$. Applying the tube estimate at
		radius $r/4$ gives
		\[
		Ncr^n
		\le\frac{Cr_0^{n-2}r^2}
		{\widetilde{\delta}\,\rho_{\widetilde\omega}(\widetilde{\delta})^2},
		\qquad
		N\le\frac{Cr_0^{n-2}r^{2-n}}
		{\widetilde{\delta}\,\rho_{\widetilde\omega}(\widetilde{\delta})^2}.
		\]
		Maximality gives $\Sigma\subset\bigcup_{i=1}^NB_r(x_i)$, and hence
		\[
		\left|\{x:\operatorname{dist}(x,\Sigma)<r\}\right|
		\le\sum_{i=1}^N|B_{2r}(x_i)|
		\le\frac{Cr_0^{n-2}r^2}
		{\widetilde{\delta}\,\rho_{\widetilde\omega}(\widetilde{\delta})^2}.
		\]
		This is the asserted upper Minkowski-content bound and implies
		$\overline{\dim}_{\mathrm M}\Sigma\le n-2$.
	\end{proof}
	
	\raggedbottom

\end{document}